\documentclass[11pt,reqno,draft]{amsart}
\usepackage{amssymb,verbatim}
\def\sh{\sinh}
\def\ch{\cosh}

 \let\Bbb\mathbb \let\Cal\mathcal
\let\le\leqslant \let\ge\geqslant \let\phi\varphi
 \let\TIL\widetilde \let\OVER\overline
\let\opn\operatorname \let\opl\operatornamewithlimits
\def\CC {{ \Bbb C {\,} }} \def\RR {{ \Bbb R {\,} }} \def\ZZ {{ \Bbb Z {\,} }}
\newcommand{\diag}{\operatorname{diag}}
\def\vac{{1\!\!\Bbb I}}

\makeatletter
\def\part{\@startsection{part}{0}%
\z@{\linespacing\@plus\linespacing}{.5\linespacing}%
{\normalfont\scshape\centering}}                
\def\section{\@startsection{section}{1}%
\z@{.7\linespacing\@plus\linespacing}{.5\linespacing}%
{\normalfont\bfseries\centering}}    
\def\l@part{\@tocline{1}{6pt plus 1pt}{0pc}{}{\hfil\scshape}} 
\def\l@section{\@tocline{1}{4pt}{0pc}{}{}} 
\def\l@subsection{\@tocline{2}{0pt}{0.4pc}{5pc}{}}
\def\l@subsubsection{\@tocline{3}{0pt}{0.4pc}{7pc}{}}
\makeatother
\theoremstyle{definition}
\newtheorem*{DEF*}{Definition} 
\newtheorem{DEF}{Definition} 

\theoremstyle{remark}
\newtheorem*{REM*}{Remark} 
\newtheorem*{NOTE*}{Note}

\theoremstyle{plain}   
\newtheorem*{COR*}{Corollary}
\newtheorem{COR}{Corollary}

\theoremstyle{remark}
\newtheorem{REM}{Remark} 

\theoremstyle{plain}  
\newtheorem{PROP}{Proposition}[section]  
\newtheorem{THM}{Theorem}[section]       

\numberwithin{equation}{section}         

\renewcommand{\sectionname}{\bfseries\S} 

\makeatletter
\def\@setaddresses{\par
  \nobreak \begingroup
\footnotesize
  \def\author##1{\nobreak\addvspace\bigskipamount}%
  \def\\{\unskip, \ignorespaces}%
  \interlinepenalty\@M
  \def\address##1##2{\begingroup
    \par\addvspace\bigskipamount\indent
    \@ifnotempty{##1}{(\ignorespaces##1\unskip) }%
    {\scshape\ignorespaces##2}\par\endgroup}%
  \def\curraddr##1##2{\begingroup
    \@ifnotempty{##2}{\nobreak\indent{\itshape Current address}%
      \@ifnotempty{##1}{, \ignorespaces##1\unskip}\/:\space
      ##2\par}\endgroup}%
  \def\email##1##2{\begingroup
    \@ifnotempty{##2}{\nobreak\indent{\itshape \EMAILADDRESS}%
      \@ifnotempty{##1}{, \ignorespaces##1\unskip}\/:\space
      \rmfamily##2\par}\endgroup}%
  \def\urladdr##1##2{\begingroup
    \@ifnotempty{##2}{\nobreak\indent{\itshape URL}%
      \@ifnotempty{##1}{, \ignorespaces##1\unskip}\/:\space
      \ttfamily##2\par}\endgroup}%
  \addresses
  \endgroup
}
\makeatother
\def\EMAILADDRESS{E-mail}
\address{A.~M.~Vershik, St.~Petersburg Department of Steklov Institute
of Mathematics, Fontanka 27, 191023 St.~Petersburg, Russia}
\email{vershik@pdmi.ras.ru}
\address{M.~I.~Graev,  Institute for System Studies, 36-1
Nakhimovsky pr.,
117218 Moscow, Russia}
\email{graev\_36@mtu-net.ru}
\title{}
\begin{document}



\begin{center}
\sc
The structure of complementary series and special representations
\\
of the groups $O(n,1)$ and $U(n,1)$.
\\ 
\bigskip
\rm A.~M.~Vershik,\footnote{Supported by the grants
NSh.--4329.2006.1 and CRDF RUM1-2662-ST-04.}
M.~I.~Graev\footnote{Supported by the grant RFBR 01-04-00363.}
\end{center}
\begin{abstract}{}
We give a survey of several models of irreducible 
complementary series representations  and their limits, special
representations, for the groups
$SU(n,1)$ and $SO(n,1)$, including new ones. These groups, 
whose geometrical meaning is well known, exhaust the list of simple
Lie groups for which the identity representation is
not isolated in the space of irreducible unitary representations
(i.e., which do not have the Kazhdan property) and hence
there exist irreducible unitary representations of these
groups --- so-called ``special representations'' --- for which the first
cohomology of the group with coefficients in these representations 
is nontrivial. By technical reasons, it is more convenient to consider
the groups $O(n,1)$ and $U(n,1)$. Most part of the paper is devoted to the group
$U(n,1)$.

The main emphasis
is on the so-called commutative models of special and complementary series
representations: in these models, the maximal unipotent subgroup
is represented by multiplicators in the case of
$O(n,1)$, and by the canonical model of the Heisenberg representations in the case of
$U(n,1)$. Earlier, these models were studied only for the group
$SL(2,\RR)$. They are especially important for realization of 
nonlocal representations of current groups, which will be considered
elsewhere.

We substantially use the ``density'' of the irreducible representations
 under study of $SO(n,1)$: their restrictions to the maximal parabolic
subgroup $ {{P}} $ are equivalent irreducible representations. Conversely, 
in order to extend an irreducible representation of
${{P}}$ to a representation of $SO(n,1)$, we must additionally define only one involution.
For the group $U(n,1)$, the situation is similar but slightly more
complicated.

UDC: 517.5.
\end{abstract}

\maketitle
        \setcounter{tocdepth}{1}
{\footnotesize\tableofcontents}

 \section*{Introduction}

\subsection{Special representations}
An irreducible unitary nonidentity representation
$T_g$ of a group $G$ in a Hilbert space $L$ is called
{\it special} if it can be extended to an indecomposable reducible representation
$\hat T$ in a one-dimensional extension
$L \oplus \{\xi \}$ of the space $L$: ${\hat T}_g (h,\xi) =(T_g h
+b(g),\xi)$, where the map $b: G\to L$ satisfies the condition
$$
b(g_1g_2) = b(g_1) + T_{g_1} b(g_2) \quad \text{for any $g_1,g_2 \in G$.\footnotemark}
$$
\footnotetext{Note that the term ``special representation'' in
\cite{NO} has another meaning.}
Such a map is called a $1$-cocycle, and the 
representation is indecomposable if and only if
this cocycle is nontrivial, i.e., the function
$b(g)$ cannot be written in the form $b(g) = T_g \xi - \xi $ for any
$\xi \in L$. Given a representation
$T$ of the group $G$, we define, in the usual way, the additive 
first cohomology group
$H^1(G,T)$ of $G$ with values in this representation:
the quotient of the additive group of cocycles by
the subgroup of trivial (i.e., cohomologous to zero) cocycles.
Thus the fact that an irreducible representation $T$ is special means
that the first cohomology with values in this representation is nontrivial.

Cocycles with values in the identity representation of a group are
additive one-dimensional characters of this group; 
hence for groups $G$ having such nontrivial
characters, it is natural to consider that the identity representation is
also special. In this case, the cohomology group of the identity 
representation is precisely
the group of additive characters of $G$.

In \cite{VK}, it is proved that every special nonidentity representation of
a locally compact compactly generated group cannot be separated 
(in the Fell--Jacobson topology)
from the identity representation; in other words, every open neighborhood
of a special representation has a nonempty intersection with every 
open neighborhood of the identity representation. Thus for groups
having a nonidentity special representation, the identity representation
is not isolated in the set of all unitary irreducible representations;
in other words, such groups do not have the Kazhdan property $(T)$.
Starting from the pioneering paper \cite{Kazdan-FA} and subsequent presentations
(see, e.g., \cite{DK}), many studies have been devoted to groups with the Kazhdan 
property (one of the recent books is \cite{BHV}); but we are interested in
semisimple groups that do not have this property. In \cite{VK},
the set of nonidentity irreducible unitary representations
that cannot be separated from the identity representation
was called the {\it core},
and these representations themselves, by a clear reason, were called
{\it infinitesimal representations}. 
As one can easily check on the example
of the group of Euclidean motions, not every infinitesimal representation
is special.

It may happen that a special representation belongs to the closure
of the identity representation. This is equivalent to the fact that 
a special representation contains almost invariant vectors
($\forall\,\varepsilon >0\,\,
\forall\, g_1,g_2, \dots g_k \in G\, \exists\, h\in L,\, ||T_{g_i}h-h||<\varepsilon$).

In \cite{VK}, it was also conjectured that every locally compact group
without property~$(T)$ has a special representation (possibly, the identity one).
This was proved by Y.~Shalom
\cite{S} (see also \cite{BHV} and references therein), even in a stronger
form: every such group has a special representation, and, moreover, not only 
the ordinary cohomology group with values in this representation
is nontrivial, but the strict
cohomology group is also nontrivial. (The group of strict cohomologies
is the quotient of the group of cocycles by
the {\it closure} of the subgroup of trivial cocycles,
see \cite{BHV}). However, the problem of computing the cohomology and
explicitly describing special representations, and even that of finding the 
number of such representations, 
is far from being solved. It is known that for Abelian and nilpotent groups,
only the identity representation is special. Even for countable solvable groups,
the situation is more complicated. The cohomology group was studied 
in a number of papers, see, e.g.,
\cite{Del} and the books \cite{I, Gui-1980}; for the cohomology of the group
of automorphisms of a tree, see
\cite{KV}.

The squared norm of a $1$-cocycle is a conditionally positive definite
function on the group. For semisimple groups of rank~1, such  functions were described in
\cite{V-G-2} and, slightly later and independently, in
\cite{FH}. In this case, the fact that a cocycle is not cohomologous to zero 
if and only if its norm, regarded  as
a function on the group, is unbounded, follows from the results of
\cite{V-G-2}. In the general form, for arbitrary unitary representations, 
it was apparently first proved in
\cite{HV} (see also \cite{BHV}).

Among the simple Lie groups, only
$SO(n,1)$ and $SU(n,1)$ have special representations, see
\cite{V-G-2}. The groups $SO(n,1)$, $n>2$, have exactly one special representation,
and the groups $SU(n,1)$ and $SO(2,1)$ have exactly two special representations
\cite[pp.~168--169]{Gui-1980}. They are the subject of this paper. As to
groups of rank higher than~1 and the groups
$Sp(n,1)$, for which the identity representation is isolated in the space
of irreducible unitary representations, similar opportunities appear for them
provided that we consider nonunitary representations 
with invariant bilinear, but not
positive definite, form; in this wider class, special representations exist
for every simple Lie group, but it seems that this opportunity
is not yet investigated in detail.

Though  for groups of rank~$1$, special representations as such were of course
known before the paper \cite{V-G-2}, only that paper made clear their link 
to cohomology, as well as to unbounded conditionally positive definite functions
and the canonical state from \cite{V-G-1};
see comments on these two papers in
\cite{Gui-exposee486, G.Segal}.

One of the important applications of nonidentity special  representations
is as follows: using such a representation, one can construct an irreducible nonlocal
unitary representation of the corresponding current group
$G^X$, i.e., the group of bounded measurable maps 
$X \to G$ with pointwise multiplication, where $X$ is a measurable space;
this question is considered in the series of papers
\cite{V-G-1,V-G-2,V-G-4,V-G-5,89-1982,V-G-6,V-G-2005,V-G-indag},
\cite{125-1993,131-1994,Ber-1}. However, we do not consider these applications
in this paper.\footnote{For representations of infinite-dimensional groups, see
also \cite{Araki,Araki-Woods,Part-S-1,Part-S-2,Part-S-LN,Streater-1,Streater-3},
\cite{69-1975,70-1975,
113-1990,O-15}, and the monographs \cite{I,Gui-1972}.
See also the monographs \cite{Gui-1980,V-G-7} and the papers \cite{Ber-4,TS-V}.
}

This paper is devoted to irreducible unitary representations of the groups
$O(n,1)$ and $U(n,1)$, i.e., the groups of linear transformations in 
$\RR^{n+1}$ and $\CC^{n+1}$ preserving, respectively, the quadratic and Hermitian forms
of signature $(n,1)$. Note that, up to central subgroups,
they are the groups of motions of the real and complex Lobachevsky spaces
$O(n,1)/O(n) \times O(1)$ and $U(n,1)/U(n) \times U(1)$.\footnote{
Replacing the groups $SO(n,1)$ and $SU(n,1)$ by their extensions
$O(n,1)$ and $U(n,1)$, while not causing any substantial changes, 
simplifies constructions.
}

We study only the special representations of the groups
$O(n,1)$ and $U(n,1)$ and the closely related complementary series representations,
leaving aside the general representation theory of these groups,
which is presented in an extensive literature.
We place emphasis on the properties that distinguish these representations 
and their realizations from other representations; in particular,
\par\hangindent=\parindent
{\it the restriction of a special representation of the group
$O(n,1)$ or $U(n,1)$ to the maximal parabolic subgroup
${{P}}$ remains a special irreducible representation.}
\par\noindent
Thus the study of the special representations of
$O(n,1)$ and  $U(n,1)$ reduces to studying the special representations
of their maximal parabolic subgroups.

The structures of the complementary series representations
and the special representations in the case
of $U(n,1)$ are richer and more beautiful than the analogous structures in the case
of $O(n,1)$. We consider these groups separately; the case of
$O(n,1)$ is presented more briefly, and the main part of the paper
concerns the group $U(n,1)$. For representations of the group 
$O(n,1)$, see also \cite{V-G-2005,V-G-indag}.

The theorems presented in the paper describe various models of the representations
under study, including new ones. 

\subsection{The special representations of the groups 
$O(n,1)$ and $U(n,1)$ as limits of the complementary 
series  representations}\label{sect:0.2}
The special representations of
$O(n,1)$ and $U(n,1)$ can be obtained by passing to the limit from 
the so-called complementary series representations, which are
of independent interest.

Let $G$ be one of these groups and ${{P}}$ be its
maximal parabolic subgroup (which is unique up to conjugation). 
Complementary series representations of the 
group $G$ are its irreducible unitary representations induced from 
one-dimensional \emph{\bf nonunitary} representations of the subgroup ${{P}}$.
These representations depend on a real parameter
$\lambda $, which lies in the interval $(0,n{-}1)$ in the case of $O(n,1)$,
and in the interval $(0,2n)$ in the case of $U(n,1)$. Their description in the case of 
$G=SL(2,K)$ can be found in \cite{IG} for $K=\RR$ and $\CC$, and in
\cite{G-G-P} for a non-Archimedean field $K$.

Complementary series representations can be realized, for example, in
Hilbert spaces $L^{\lambda }$ of functions on the sphere
$S = G/{{P}}$. Other realizations are described below.

At the endpoints of the interval $(0,n{-}1)$ in the case of $O(n,1)$
and, respectively, $(0,2n)$ in the case of $U(n,1)$,
the space $L^{\lambda }$ becomes reducible.

In the case of $O(n,1)$, for $\lambda = 0$ the norm $\|f\|_{\lambda }$
in the space $L^\lambda$ degenerates on a subspace $L$
of codimension~$1$. The special representation $T$ of
$O(n,1)$ is realized in this subspace $L$ with the norm
$$
\|f\|=\lim\limits_{\lambda \,\to\,  0} \frac{ d\|f\|_{\lambda } }{ d \lambda }.
$$
The associated nontrivial $1$-cocycle is given by the formula
$b(g) = T_g \xi - \xi $, where $\xi \notin L$.

For $\lambda = n{-}1$, the norm $\|f\|_{\lambda }$
degenerates on a one-dimensional subspace $L_0$. The special representation
is realized in the quotient by
$L_0$ with the norm
$$
\|f\|_0 = \lim\limits_{\lambda\, \to \, n-1} \|f\|_{\lambda } .
$$
The special irreducible unitary representations of the group 
$O(n,1)$ corresponding to the endpoints of the interval
$(0,n{-}1)$ are equivalent.

\smallskip

In the case of $U(n,1)$, for $\lambda = 0$ the norm
$\|f\|_{\lambda }$ also degenerates on a subspace $L$ of codimension~$1$.
In contrast to the case of $O(n,1)$, the norm
$\|f\|=\lim\limits_{\lambda \,\to\,  0} \frac{ d\|f\|_{\lambda } }{ d \lambda }$
in $L$ for ${{n>1}}$ degenerates on some infinite-dimensional
subspace $L_1 \subset L$. The Hilbert space
$L/L_1$ (with the inherited norm) splits into the direct sum of two
irreducible nonequivalent invariant subspaces
$H_+$ and $H_-$. The representations of the group $U(n,1)$
in the subspaces $H_+$ and $H_-$ are special representations.

For $\lambda = 2n$ with ${{n>1}}$, the norm $\|f\|_{\lambda }$
degenerates on a subspace $L$ of infinite codimension. The norm
$\|f\|=\lim\limits_{\lambda \,\to\,  2n} \frac{ d\|f\|_{\lambda } }{ d \lambda }$
on $L$ degenerates on a one-dimensional subspace
$L_0 \subset L$. The representation of the group $U(n,1)$
in the Hilbert space $L/L_0$ splits into the direct sum of two
nonequivalent special irreducible representations.

\smallskip

The special unitary representations of
$U(n,1)$ corresponding to the endpoints of the interval
$(0,2n)$ are equivalent.

\subsection{Matrix realizations of the groups
$O(n,1)$, $U(n,1)$ and of some subgroups of these groups}

Matrix realizations of the groups $O(n,1)$ and $U(n,1)$
are determined by the choice of a quadratic (respectively, Hermitian) form 
of signature $(n,1)$ in the space $\RR^{n+1}$ (respectively, $\CC^{n+1}$). This choice is in turn
equivalent to fixing a Cartan subgroup of the corresponding group.

We realize $O(n,1)$ and $U(n,1)$ as the groups of linear transformations
in $\RR^{n+1}$ and $\CC^{n+1}$ preserving, respectively, the quadratic form
$2x_1 x_{n+1}+x_2^2+ \ldots +x_n^2$ and the Hermitian form
$x_1\OVER x_{n+1}+x_{n+1}\OVER x_1+|x_2|^2+ \ldots +|x_n|^2.$

In this realization, elements of the groups can be written as block
matrices $g=\|g_{ij}\|$ of order $3$, where the diagonal blocks
$g_{11}$, $g_{22}$, and $g_{33}$ are square matrices of orders $1$,
$n{-}1$, and $1$, respectively. The maximal unipotent subgroups
$Z \subset O(n,1)$ and $H \subset U(n,1)$ are, respectively, the groups of
matrices of the form
$$
z=\begin{pmatrix}{}1&0&0\\
-{\gamma}^*&e&0\\ -\frac{{\gamma} {\gamma}^*}{2}&{\gamma}&1\end{pmatrix},
\,\, {\gamma} \in\RR^{n-1}, \quad \text{and} \quad
h=\begin{pmatrix}{}1&0&0\\
-z^*&e&0\\ it-\frac{z z^*}{2}&z&1\end{pmatrix},
\,\, t\in\RR, \,\,  z \in\CC^{n-1},
$$
where $e$ is the unit matrix of order $n{-}1$,
${\gamma}{\gamma}^* = {\gamma}_1^2 + \ldots + {\gamma}_{n-1}^2$,
$zz^* = |z_1|^2 + \ldots + |z_{n-1}|^2$.

It is convenient to write elements of $H$ as pairs
$(t,z)$, $t\in\RR$,
$z\in\CC^{n-1}$, with the multiplication law
$$
(t_1,z_1)\,(t_2,z_2)=(t_1+t_2-\opn{Im} z_1z_2^*,z_1+z_2).
$$

{\it The normalizer of these subgroups is the subgroup
${{P}}$ of lower block triangular matrices, i.e., the maximal parabolic
subgroup.} It is essential that each of the groups
$O(n,1)$ and $U(n,1)$ is algebraically generated by the elements of 
${{P}}$ and one element
$$
s=\left(\begin{array}{ccc}0&0&1\\ 0&e&0\\ 1&0&0\end{array}\right).
$$
Thus representations of the groups
$O(n,1)$ and $U(n,1)$ are completely determined by the operators
corresponding to the elements of the subgroup
${{P}}$ and this element $s$.

The group ${{P}}$ is the semidirect product
${{P}}=Z\leftthreetimes D$ and
${{P}}=H\leftthreetimes D$, respectively, where
$D$ is the subgroup of block-diagonal matrices of the form
$$
d = \diag (\epsilon^{-1}, u, \epsilon),  \quad
\epsilon \in\RR^*,\quad u\in O(n{-}1),\quad \text{in the case of}\quad O(n,1).
$$
$$
d = \diag (\bar\epsilon^{-1}, u, \epsilon),  \quad
\epsilon \in\CC^*,\quad u\in U(n{-}1),\quad \text{in the case of}\quad U(n,1).
$$
The subgroup $D$ can be written as the direct product
$D = D_0 \times D_1$, where $D_0$ is the maximal compact subgroup of
$D$. In the case of $U(n,1)$, $D_0$ is the subgroup of block-diagonal matrices 
$d=\diag(\epsilon ,u, \epsilon )$, $|\epsilon |=1$, $u\in U(n{-}1)$, and
$D_1\cong\RR^*_+$ is the subgroup of block-diagonal matrices
of the form $d=\diag(r^{-1},e,r)$, $r>0$.
In the case of $O(n,1)$, the subgroups
$D_0$ and $D_1$ are defined in a similar way.

{\it By ${{P}}_0$ we denote a subgroup of 
${{P}}$ of codimension $1$, namely,
${{P}}_0 = Z \leftthreetimes D_0$ and ${{P}}_0 = H \leftthreetimes D_0$,
respectively.}

\subsection{The difference in the structures of 
complementary series  representations and special representations of the groups
$O(n,1)$ and $U(n,1)$}

The difference in the structures of representations of
the groups $O(n,1)$ and $U(n,1)$ is predetermined 
by the structure of their maximal unipotent subgroups
$Z$ and $H$, respectively.

The subgroup $Z \subset O(n,1)$ is commutative and isomorphic to the additive 
group $\RR^{n-1}$; thus all irreducible representations of $Z$ are one-dimensional.
Hence in the case of $O(n,1)$ there exists a model of representations in which
the operators corresponding to the unipotent subgroup act as multiplicators.
We call it the commutative model associated with the unipotent subgroup.

The subgroup $H \subset U(n,1)$ is isomorphic to the Heisenberg group of 
dimension $2n{-}1$, so that it has a one-parameter family of 
infinite-dimensional irreducible unitary representations. In this case, the analog
of the commutative model for $O(n,1)$ is a direct 
integral of multiples of
irreducible infinite-dimensional subspaces of the Heisenberg group. 
Thus the theory of complementary series representations
and special representations of
$U(n,1)$ merges with the representation theory of the Heisenberg group.

Note that the Heisenberg subgroup $H \subset  U(n,1)$,
in contrast to the maximal unipotent subgroup of
$O(n,1)$, has the following remarkable property. Every 
infinite-dimensional irreducible unitary representation of $H$
can be (uniquely) extended to a unitary representation of the subgroup
${{P}}_0 \subset {{P}}$ of codimension $1$:
${{P}} = {{P}}_0 \leftthreetimes D_1$, $D_1 \cong \RR^*_+$.

It is known that the groups
$O(n,1)$ and $U(n,1)$ are algebraically generated by the elements of the subgroup
${{P}}$ and one element $s$ of order $2$. Thus if the space of a unitary
representation of $U(n,1)$ is decomposed into a direct integral
of infinite-dimensional subspaces invariant and irreducible with respect
to the Heisenberg subgroup $H$, then, in order to describe this representation,
it suffices to define only the action of the operators of the one-parameter
subgroup $D_1$ and the operator corresponding to the element $s$.

\subsection{Restrictions of representations of
$O(n,1)$ and $U(n,1)$ to the maximal compact and parabolic subgroups}

The difference between $O(n,1)$ and $U(n,1)$ shows itself also in the
restrictions of complementary series representations and special representations
to the maximal compact subgroup $K$ and the maximal parabolic subgroup
${{P}}$.

In both cases, the restrictions of complementary series  representations
to $K$ are pairwise equivalent
and contain the identity representation with multiplicity one;
the spectrum of these restrictions is simple only in the case of $O(n,1)$.
When passing to the special representations, in the case of
$O(n,1)$, the identity representation vanishes; in the case of
$U(n,1)$, $n>1$, besides the identity representation, a certain infinite-dimensional
representation of $K$
also vanishes. Each of the remaining irreducible representations of
$K$ occurs in one of the two special
representations with multiplicity one.

In both cases, the restrictions of all
complementary series  representations to $P$ are equivalent to one and the same representation of
${{P}}$ induced from the identity representation of its
block-diagonal subgroup.

In the case of $O(n,1)$, $n>2$, this representation of 
${{P}}$ is irreducible and equivalent to the restriction to
${{P}}$ of the special representation. Thus on the subgroup
${{P}}$ these representations coincide. On the whole group
$O(n,1)$ they differ by the action of the single operator
$T_s$.

In the case of $U(n,1)$, this representation of
${{P}}$ decomposes into a direct sum of pairwise nonequivalent
irreducible representations. {\it The restrictions 
of the two special representations to $P$ are irreducible and occur
in this decomposition as two direct summands.}

Observe the difference between the restrictions
of complementary series representations and special representations of
the group $U(n,1)$  to the Heisenberg subgroup  $H$.

The decomposition of every complementary series  representation 
contains every infinite-dimensional representation of the Heisenberg group
with infinite multiplicity.

{\it The sum of the two special representations of
$U(n,1)$ decomposes into the direct integral of all pairwise nonequivalent
representations of the Heisenberg group. }

\par\medskip

In the remaining part of the Introduction we briefly describe
the main results of the paper.

\subsection{Realization of the complementary series representations 
and the special representation for the group $O(n,1)$}

In \S\,\ref{sect:2} we describe
six models of the complementary series representations and
six models of the special representation of the group
$O(n,1)$; these models differ by the choice of coordinate system,
functional space, and the endpoints of the interval corresponding
to the complementary series.

More exactly, the invariant bilinear functional that determines the
Hilbert space of the complementary series representation
$T^{\lambda }$ can be defined for any value of the parameter
$\lambda \in [0,n{-}1]$. The values 
$\lambda $ and  $(n{-}1){-}\lambda $
(that are distinct from the endpoints of the interval) give rise to equivalent
complementary series representations 
$T^{\lambda }$ and $T^{\,n-1- \lambda }$;
their realizations are different.

This difference in the realization of
$T^{\lambda }$ and $T^{\,n-1- \lambda }$ is especially vivid
when one compares the cases $\lambda =0$ and $\lambda = n{-}1$,
when these representations become reducible and nonequivalent. Starting from
$T^{0}$ and $T^{\,n-1}$, we can obtain two models of the unique
special irreducible representation of the group
$O(n,1)$. These models are different. Namely, as observed above, for
$ \lambda = 0 $ the special representation is realized in an
invariant subspace of codimension $1$, while for
$\lambda = n{-}1$ it is realized in the quotient by
a one-dimensional invariant subspace.

\smallskip

Realization of representations also depends on the choice of
coordinate system on the sphere $S^{n-1} \subset \RR^n$.
Model $A$ uses spherical coordinates, while
model $B$ uses the local coordinates of the punctured sphere, which
determine a bijection between $\RR^{n-1}$ and the unipotent subgroup $Z$.

Thus we obtain four models for each complementary series  representation
and four models for the special representation of the group
$O(n,1)$. Two of these four models of the special representation are obtained
for $\lambda =0$, and the other two models, for $\lambda = n{-}1$.

Besides, one can pass to the Fourier transform in the space
of model $B$ realized in the local coordinates.
Then the operators corresponding to the unipotent subgroup
$Z$ will act as multiplicators, and we obtain two more models of
representations of $O(n,1)$
(models $C$, or commutative models), which are especially important
in applications to the representation theory of the current group
$O(n,1)^X$ \cite{V-G-2005,V-G-indag}.

In particular, this realization is convenient for obtaining an explicit formula
for the nontrivial 
$1$-cocycle $b:\, O(n,1)\to L^0$ associated with the special representation
(see Section~\ref{sect:b(g):O(n,1)}).

In the commutative model with $\lambda =0$, this $1$-cocycle
is given by the formula
$$
        b(g)=T_g \xi - \xi ,
$$
where
$$
\xi = |\xi|^{-\frac{1-n}{2} }\,
 K_{\frac{1-n}{2}}(\sqrt2\,|\xi|)).
$$
Here $K_{\rho}(x)$ is a Bessel function
($\xi$ is a vector that does not belong to $L^0$ and is invariant under
the maximal compact subgroup $K$ of the group $O(n,1)$).
In particular,
$$
\xi = e^{\,- \sqrt2\,|\xi |}\quad\mbox{for}\quad n=2.
$$

For each model of the complementary series representations, we explicitly describe the
spherical function and the vacuum vector
$\xi $, i.e., a vector from the representation space that is invariant
under the maximal commutative subgroup.
We also describe the embedding of the spaces of complementary series representations
into their tensor product. This embedding leads to one
of the constructions of representations of the current group  $O(n,1)^X$.

\subsection{Construction of an irreducible representation of the subgroup
${{P}} \subset U(n,1)$}

In \S\S\,\ref{sect:3}--\ref{sect:6} we consider representations of the group
$U(n,1)$ and its parabolic subgroup ${{P}}$.

We begin in \S\,\ref{sect:3} with a construction of a class of unitary representations
of the group ${{P}}$ based on the construction of irreducible
representations of the Heisenberg group $H$. It is known that every
infinite-dimensional unitary representation of
$H$ is determined by a real parameter $\rho \ne 0$ (the Planck constant).
We consider the Bargmann realization of these representations.  In this
realization, the representation with parameter
$\rho >0$ is defined in the Hilbert space $\Cal H(\rho )$ 
of entire analytic functions
$f(z)$ on $\CC^{n-1}$ with the norm
\begin{equation}{}\label{11-00}
\| f \|^{2} = |\rho |^{n-1}
\int_{\CC^{n-1}} |f(z)|^{2}\,e^{\, - |\rho| \,|z|^2}\,d\mu(z),
\end{equation}
where $|z|^2=zz^*=|z_1|^2+ \ldots +|z_{n-1}|^2$ and $d\mu(z)$ is the Lebesgue measure on
$\CC^{n-1}$. The operators corresponding to elements
$h=(t_0,z_0)\in H$ have the form
\begin{equation}{}\label{11-01}
T_{h} f(z) = e^{\,\rho  \,(it_0 - \frac12\,|z_0|^2 - zz_0^*)}\,f(z+z_0).
\end{equation}

Analogously, the representation with parameter $\rho <0$ is defined in the
Hilbert space $\Cal H(\rho )$ of entire antianalytic  functions on
$\CC^{n-1}$ with the norm \eqref{11-00}.
The operators corresponding to elements
$h=(t_0,z_0)\in H$ have the form
\begin{equation}{}\label{11-02}
T_{h} f(z) = e^{\,i\rho t_0-|\rho |\,( \frac12\,|z_0|^2 + z_0z^*)}\,f(z+z_0).
\end{equation}

Each of these representations can be extended to a unitary representation
of the group ${{P}}_0=H\leftthreetimes D_0$. Namely, the operators $T_d$
corresponding to elements $d=(\epsilon ,u, \epsilon  )\in D_0$ 
are given, for every $\rho $, by the formulas
\begin{equation}{}\label{11-03}
T_df(z)=f(\epsilon zu).
\end{equation}
To construct unitary representations of the group
${{P}}$, we associate with every  $\lambda \ge 0$ direct integrals of
the spaces $\Cal H(\rho )$:\label{sect:Hlambda}
$$
\Cal H^{\lambda }_+=\int_0^\infty \Cal H(\rho )\,
\rho ^{\lambda -1}\, d \rho \quad \text{and} \quad
\Cal H^{\lambda }_-=\int^0_{-\infty} \Cal H(\rho )\,
|\rho| ^{\lambda -1}\, d \rho.
$$
The representations of the group ${{P}}_0$ in these spaces can be extended to
unitary representations of the group ${{P}}={{P}}_0\leftthreetimes D_1$. 
Namely, the operators $T_d$ corresponding to elements
$d=\diag(r^{-1},e,r)\in D_1$ have the form 
\begin{equation}{}\label{11-04}
T_d f(\rho ,z)=f(r^2 \rho , r^{-1}z)\,r^{\lambda }.
\end{equation}

The representations of the group
${{P}}$ in the spaces
$\Cal H^{\lambda }_+$ (respectively, 
$\Cal H^{\lambda }_-$) defined in this way are irreducible and pairwise equivalent;
the representations in $\Cal H^{\lambda }_+$ and $\Cal H^{\lambda }_-$
are not equivalent.

In \S\,\ref{sect:6} it is shown that these representations of the group
${{P}}$ can be extended to special representations of the group $U(n,1)$.

\subsection{Representations of the group ${{P}}$ in the spaces 
$H^{\lambda }_{\pm}$}
\label{sect:0.8}

In  \S\,\ref{sect:4} we construct a family of reducible unitary 
representations of the subgroup
${{P}} \subset U(n,1)$ in spaces
$H^{\lambda }_{+}$ and $H^{\lambda }_{-}$. As shown
later in~\S\,\ref{sect:5}, the representations of the subgroup ${{P}}$
in the spaces $H^{\lambda }=H^{\lambda }_{+}\oplus H^{\lambda }_{-}$
can be extended to (irreducible)  complementary
series representations of the whole group $U(n,1)$.

The spaces $H^{\lambda }_{\pm}$ are defined as the direct integrals
$$
H^{\lambda }_+ = \int_0^\infty  H(\rho )\,
\rho ^{\lambda -1}\, d \rho \quad \text{and} \quad
H^{\lambda }_-=\int^0_{-\infty} H(\rho )\,
|\rho| ^{\lambda -1}\, d \rho,
$$
where $H(\rho )$ is the Hilbert space of
{\it all} functions $f(z)$ on $\CC^{n-1}$ with the norm \eqref{11-00}.

In $H(\rho )$, a unitary representation of the group
${{P}}_0$ acts; it is given by \eqref{11-01} and \eqref{11-03}
for $\rho >0$, and
by \eqref{11-02} and \eqref{11-03}  for $\rho <0$.

The space $H(\rho )$ decomposes into a direct sum of invariant 
${{P}}_0$-irreducible  pairwise nonequivalent subspaces:
$$
H(\rho )= \bigoplus _{m=0}^\infty H_m(\rho ),
$$
where $H_m(\rho )$ with $\rho >0$ and $\rho <0$ is cyclically generated by
the vectors
$\OVER z^k=\OVER z_1^{k_1} \ldots \OVER z_{n-1}^{k_{n-1}}$
and $z^k=z_1^{k_1} \ldots  z_{n-1}^{k_{n-1}}$, respectively, where $k_1 + \ldots + k_{n-1}=m$.
In particular, $H_0(\rho )=\Cal H(\rho )$.

The representation of the Heisenberg subgroup in each subspace
$H_m(\rho )$ is a multiple of the irreducible representation with parameter
$\rho $. The multiplicity is equal to the number of partitions of
$m$ into the sum of $n{-}1$ nonnegative integer summands.

Each space $H(\rho )$ has an orthogonal basis
$\{f_{pq} \mid p,q\in\ZZ^{n-1}_+\}$ compatible
with the action of the Heisenberg group $H$. For every fixed
$p$, the vectors $f_{pq}$, $q\in\ZZ^{n-1}_+$, form an orthogonal basis
in the subspace of an irreducible representation of the group
$H$. We describe the matrix elements of
the representations of the group ${{P}}_0$ in the same subspace 
in terms of this basis.

The representations of the group ${{P}}_0$ in the subspaces $H^{\lambda }_\pm$
can be extended, according to \eqref{11-04}, to unitary
representations of the group ${{P}}$. These spaces are reducible with respect to
${{P}}$ and decompose into a direct sum of $P$-irreducible pairwise
nonequivalent subspaces:
$$
H^{\lambda }_\pm = \bigoplus _{m=0}^\infty (H^{\lambda }_\pm)_m\,\,,
$$
where
$$
(H^{\lambda }_+)_m = \int_0^\infty H_m(\rho )\,\rho ^{\lambda -1}\,d \rho
\quad \text{and}\quad
(H^{\lambda }_-)_m= \int^0_{-\infty} H_m(\rho )\,|\rho|^{\lambda -1}\,d \rho .
$$

\subsection{Models of the complementary series representations of the group
$U(n,1)$}

In \S\,\ref{sect:5}, by analogy with the case of the group
$O(n,1)$, we construct six different models of the complementary 
series  representations of the group $U(n,1)$. We start from the known realization of these 
representations in a space of functions on the unitary unit sphere 
$S \subset \CC^n$ (model $A$), and then pass to their realization
in a space of functions on the Heisenberg group $H$, i.e., in the space
of functions $\phi (t,z)$ on $\RR \times \CC^{n-1}$
(model $B$).

In the next model (model $C$), the complementary series representations of the group  
$U(n,1)$ are realized as direct integrals of 
multiples of irreducible representations of the Heisenberg group.
They are obtained by passing from functions $\phi (t,z)$
in model $B$ to functions  $\psi (\rho ,z)$ on 
$\RR \times \CC^{n-1}$, which are given by the integral
$$
\psi(\rho, z ) = |\rho |^{-n+1}\,e^{\frac{|\rho |}{2}|z|^2}\,
        \int_{- \infty }^{+ \infty }e^{\,-i\rho t}\,
        \phi (t, z )\,d t.
$$
In this model, the space $L^{\lambda }$ of the complementary series  representation
is a direct integral of ${{P}}_0$-invariant subspaces $L(\rho )$:
$$
L^{\lambda } = \int_{- \infty }^{+ \infty }
        L(\rho )\,|\rho |^{\,\lambda -1}\,d \rho ,
$$
where $L(\rho )$ is the space of functions $\psi (z)$ on $\CC^{n-1}$ with the norm
\begin{equation*}{}
\|\psi\|^2_{\rho }=|\rho |^{2n-2}\int_{\CC^{n-1} \times \CC^{n-1}}
R_{\lambda }(\rho ,z_1,z_2)\,\psi(\rho ,z)\,
\OVER{\psi(\rho ,z)}\,d\mu(z_1)\,d\mu(z_2).
\end{equation*}
An explicit expression for $R_{\lambda }$ in terms of the Bessel function
$K_{\frac{1 - \lambda }{2}}(t)$ is given in the main text.

The operators of the subgroup ${{P}}_0$ in the space $L(\rho )$
are given by the same formulas as for the subspace 
$H(\rho )$ defined in Section~\ref{sect:0.8}.

The total number of models of each complementary series  representation
of the group $U(n,1)$ is equal to $3$ rather than $6$, because the representations
in the spaces $L^{\lambda }$
and $L^{2n-\lambda }$ are equivalent.

\par\smallskip
Let us formulate the main result of \S\,\ref{sect:5}.
It is proved that the representations of the subgroup ${{P}}_0$ 
in the spaces $L(\rho )$ and $H(\rho )$ are equivalent.
Thus there exists an isomorphism of Hilbert spaces
$$
J : H(\rho ) \to L(\rho )
$$
commuting with the operators of the group
${{P}}_0$. Since
$H(\rho ) = \bigoplus _{m=0}^{\infty } H_m(\rho )$,
where $H_m(\rho )$ are $P_0$-invariant
irreducible pairwise nonequivalent subspaces, it follows that 
on each of them the operator $J$ is a multiple of the identity operator:
$J = c_m \,\opn{id}$ on $H_m(\rho )$.
(In the main text, we explicitly calculate the coefficient $c_m$.)

The isomorphism $J : H(\rho ) \to L(\rho ) $ induces an isomorphism of Hilbert spaces
$$
H^{\lambda }\to L^{\lambda }
$$
commuting with the action of the whole group ${{P}} $
(and not only the subgroup ${{P}}_0 \subset {{P}}$).

In view of the isomorphism $ H^{\lambda }\to L^{\lambda }$, the action of the group
$U(n,1)$ can be transferred from the space $L^{\lambda }$
to the space $H^{\lambda }$.

\par\hangindent=\parindent
{\it Thus the space $H^{\lambda }$ can be regarded as a new model of
the space of a  complementary series representation.
}
\par\noindent
In this new model, the operators of the subgroup
${{P}}$ are given by the same formulas as in the space $L^{\lambda }$.
Hence, to describe a representation of the group
$U(n,1)$ in the new model, it suffices to know 
only the operator $T_s$.

\subsection{Models of the special representations of the group $U(n,1)$}

In \S\,\ref{sect:6} we construct six models of the special representations of the 
group $U(n,1)$. They are obtained from the models of the complementary series
representations by passing to the limit as
$\lambda \to 0$ and $\lambda \to 2n$.

The structure of the spaces $L^0$ and $L^{2n}$, which are the limits of
the spaces $L^{\lambda }$ of complementary series representations, was already discussed in
Section~\ref{sect:0.2}.

Let us describe the main model (model  $C$, $\lambda =0$).

The special representations of the group $U(n,1)$ are realized in the spaces
$$
\Cal L_+=\int_0^{\infty }\Cal L(\rho )\,\rho ^{-1}\,d \rho, \qquad
\Cal L_-=\int_{- \infty }^0\Cal L(\rho )\,|\rho |^{-1}\,d \rho.
$$

Here $\Cal L(\rho ) $ is the quotient of the space of functions $\psi(z)$ on
$\CC^{n-1}$ with the norm
$$
\|\psi\|^2_{\rho }=|\rho |^{2n-2}\int_{\CC^{n-1} \times \CC^{n-1}}
e^{|\rho |(a_{\rho }(z_1,z_2)-|z_1|^2-|z_2|^2)}\,
\psi(z_1)\,\OVER{\psi(z_2)}\,d\mu(z_1)\,d\mu(z_2),
$$
where
\begin{equation}{}\label{33-818}
a_{\rho }(z_1,z_2)=
\begin{cases}
 z_1z_2^* & \text{for $\rho >0$},\\[1ex]
 z_2z_1^* & \text{for $\rho <0$},
\end{cases}
\end{equation}
by the subspace of functions of norm $0$.

The spaces $\Cal L(\rho)$
are irreducible and pairwise nonequivalent
under the action of the Heisenberg group $H$. Thus the special
representations of the group $U(n,1)$ decompose into direct integrals of
irreducible pairwise nonequivalent representations
of the Heisenberg group.

There is another realization of the special representations of 
$U(n,1)$. Namely, it is proved that there exists an isomorphism of Hilbert spaces
$$
\Cal H^0_{\pm}\to\Cal L_{\pm}
$$
commuting with the action of the subgroup ${{P}}$, where $\Cal H^0_{\pm}$
are the spaces defined in~\S\,\ref{sect:Hlambda} of the Introduction.
\par\hangindent=\parindent
{\it
In view of this isomorphism, the action of the group
$U(n,1)$ can be transferred from the spaces $\Cal L_{\pm}$
to the spaces $\Cal H^0_{\pm}$.
}
\par\noindent
We call the realizations of the special representations of
$U(n,1)$ in
$\Cal H^0_{\pm}$ the Bargmann models. In these models,
the action of the operators of the subgroup
${{P}}$ is defined from the beginning, so that
in order to give a complete description of a representation of
$U(n,1)$, it suffices to determine only the action of the operator $T_s$.

It is establish that in the case of $\Cal H^0_+$ the action of  $T_s$ 
on functions of the form $\rho ^{|k|}\,e^{-a \rho }\,z^k$, where
$z^k=z_1^{k_1} \ldots z_{n-1}^{k_{n-1}}$,
$|k|=k_1+ \ldots +k_{n-1}$, and $\opn{Re} a>0$, is given by 
$$
T_s(\rho ^{|k|}\,e^{-a \rho }\,z^k)=a^{-|k|}\,
\rho ^{|k|}\,e^{-\frac{\rho }{a}}\,z^k.
$$
In particular, $T_s(e^{-a \rho })=e^{-\frac{\rho }{a}}$.
In a similar way $T_s$ is defined on the space $\Cal H^0_-$.

In the Bargmann model, the nontrivial $1$-cocycles
$b_{\pm}:U(n,1)\to \Cal H^0_{\pm}$ are given by 
$$
b_{\pm}(g)=T_g \xi _{\pm}-\xi _{\pm},\quad
\xi _{\pm}= \chi _{\pm}(\rho )\,e^{- |\rho| },
$$
where $\chi _{\pm}$ are the characteristic functions of the half-lines
$\rho >0$ and $\rho <0$, respectively. Since
$b_{\pm}(g)=0$ on the maximal compact subgroup, this function
is uniquely determined by its restriction to the subgroup
${{P}}_1$, a one-dimensional extension of the Heisenberg group $H$.

\subsection{Canonical representations}

In the construction of the special representations of the groups
$O(n,1)$ and $U(n,1)$, one can replace the complementary series with another
one-parameter series of unitary representations, the canonical representations.

Canonical representations were first introduced in
\cite{V-G-1,V-G-6} in connection with constructing representations of the current groups
$O(n,1)^X$ and $U(n,1)^X$.

In Appendix~1, we first give a general definition of canonical representations
for an arbitrary locally compact group $G$, and then explicitly construct
the canonical representations of the groups
$O(n,1)$ and $U(n,1)$ and  establish their relation to the complementary
series representations and the special representations.

In Appendix~2, we derive some formulas used in the paper.

\addcontentsline{toc}{part}{Preliminary information on the groups $O(n,1)$ and $U(n,1)$}

\section{The groups $O(n,1)$, $U(n,1)$ and some subgroups of these groups}\label{sect:1}

\subsection{The groups $O(n,1)$ and $U(n,1)$}

By $O(n,1)$ and $U(n,1)$ we denote the groups of all linear transformations
in $\RR^{n+1}$ and $\CC^{n+1}$ preserving, respectively, the quadratic and
Hermitian forms of signature $(n,1)$; this condition uniquely determines
the groups $O(n,1)$ and $U(n,1)$ up to isomorphism. In what follows, we use
two matrix models of these groups, which are determined by the choice of the
corresponding quadratic and Hermitian form.

\subsubsection{Model $a$}
The quadratic and Hermitian forms are, respectively,
$$
x_1^2+ \ldots +x_n^2-x_{n+1}^2,
$$
$$
|x_1|^2+ \ldots +|x_n|^2-|x_{n+1}|^2.
$$
In this model, elements of the groups are written as block matrices 
$$
g=\begin{pmatrix}\alpha  &\beta  \\ \gamma   &\delta \end{pmatrix},
$$
where the diagonal blocks $\alpha , \delta $ are matrices of orders $n$ and $1$,
respectively. In this notation, $O(n,1)$ and $U(n,1)$ are, respectively,
the groups of real and complex matrices $g$ of order
$n+1$ satisfying the condition
\begin{equation}{}\label{1-1}
g \sigma  g^*= \sigma ,\quad \text{where}\quad
\sigma  =\begin{pmatrix}e&0\\ 0&-1\end{pmatrix},
\end{equation}
i.e.,
$$
\alpha \alpha ^* - \beta \beta ^*=e, \quad \alpha \gamma ^* - \beta \delta ^*=0,
\quad  \delta \delta ^*- \gamma \gamma ^*=1.
$$
Here $e$ is the unit matrix of order $n$ and
$*$ denotes conjugation, i.e., 
$\alpha ^* = \alpha ^\top$ in the real case, and
$\alpha ^* = \OVER\alpha ^\top$ in the complex case.

\subsubsection{Model $b$}
The quadratic and Hermitian forms are, respectively,
$$
2x_1 x_{n+1}+x_2^2+ \ldots +x_n^2,
$$
$$
x_1\bar x_{n+1}+x_{n+1}\bar x_1+|x_2|^2+ \ldots +|x_n|^2.
$$
In this model, elements of the groups are written as block matrices
$$
g=\begin{pmatrix}{}g_{11}&g_{12}&g_{13}\\
g_{21}&g_{22}&g_{23}\\
g_{31}&g_{32}&g_{33}\end{pmatrix},
$$
where the diagonal blocks are matrices of orders $1$, $n{-}1$, and $1$, respectively.
In this notation, $O(n,1)$ and $U(n,1)$ are, respectively, the groups
of real and complex matrices $g$ of order $n+1$
satisfying the condition
\begin{equation}{}\label{1-2}
g s  g^*= s ,\quad \text{где}\quad
s  =\left(\begin{array}{ccc}0&0&1\\ 0&e&0\\ 1&0&0\end{array}\right),
\end{equation}
where $e$ is the unit matrix of order $n{-}1$.

\begin{REM*}{}
Model $a$ is convenient for describing the maximal compact subgroup
in $U(n,1)$.
Namely, the latter subgroup  in model $a$ consists of block-diagonal matrices
$g=\begin{pmatrix}u&0\\0&\epsilon \end{pmatrix}$, where
$u \in O(n)$, $\epsilon =\pm1$ in the case of $O(n,1)$;
and $u\in U(n) $, $\epsilon \in\CC^* $ in the case of $U(n,1)$.
Model $b$, mostly used in this paper, 
is convenient for describing the unipotent subgroup.
\end{REM*}

\subsubsection{Transition from model $a$ to model $b$ and vice versa}
Note that the matrices $\sigma $ and $s$ are related by 
$$
a^{-1}s a= \sigma , \quad \text{where} \quad
a  =\frac{1}{\sqrt 2}\,
\left(\begin{array}{ccr}1&0&-1\\ 0&e&0\\ 1&0&1\end{array}\right).
$$
It follows that the transition from model $a$ to model $b$ reduces 
to the transformation $g\to ag a^{-1}$. Let us give explicit formulas
for this transition. For this, let us write the blocks of matrices 
$g=\begin{pmatrix}\alpha  &\beta  \\ \gamma   &\delta \end{pmatrix}$
in model $a$ in the form
$$
\alpha =\begin{pmatrix}\alpha_{11}&\alpha_{12}\\\
\alpha_{21}&\alpha_{22} \end{pmatrix},\quad
\beta =\begin{pmatrix}\beta_1\\\beta_2 \end{pmatrix},\quad
\gamma =(\gamma _1, \gamma _2),
$$
where $\alpha _{22}$ is a matrix of order $n{-}1$ and $\beta_2$, $\gamma _2$
are, respectively, a column and a row of length $n{-}1$. Then the transition
from model $a$ to model $b$ is given by the following formula:
$$
\begin{pmatrix}\alpha  &\beta  \\ \gamma   &\delta \end{pmatrix}\to
\begin{pmatrix}\frac{1}{2}(\alpha _{11}+ \delta - \beta _1- \gamma _1)&
\frac{\alpha _{12}- \gamma _2}{\sqrt 2}&
\frac{1}{2}(\alpha _{11}- \delta + \beta _1- \gamma _1)\\
\frac{\alpha _{21}- \beta _2}{\sqrt 2}&\alpha _{22}&
\frac{\alpha _{21}+\beta_2}{\sqrt 2}\\
\frac{1}{2}(\alpha _{11}- \delta - \beta _1+ \gamma _1)&
\frac{\alpha _{12}+ \gamma _2}{\sqrt 2}&
\frac{1}{2}(\alpha _{11}+\delta + \beta _1+ \gamma _1)
\end{pmatrix}.
$$

Conversely, the transition from model $b$ to model $a$ is given by
$$
\begin{pmatrix}{}g_{11}&g_{12}&g_{13}\\
g_{21}&g_{22}&g_{23}\\
g_{31}&g_{32}&g_{33}\end{pmatrix} \longrightarrow 
\begin{pmatrix}\alpha  &\beta  \\ \gamma   &\delta \end{pmatrix},
$$
where
\begin{gather*}{}
\alpha =
\begin{pmatrix}\frac{1}{2}(g_{11} + g_{33}+g_{13}+g_{31})&
\frac{g_{12}+g_{32}}{\sqrt 2} \\
\frac{g_{21}+g_{23}}{\sqrt 2}&g_{22} \end{pmatrix},
 \quad
\beta =\begin{pmatrix}\frac{1}{2}(-g_{11}+ g_{33}+g_{13}-g_{31}\\
\frac{-g_{21}+g_{23}}{\sqrt 2}\end{pmatrix},
\\
  \gamma = (\frac{1}{2}(-g_{11}+ g_{33}-g_{13}+g_{31}),\,
\frac{-g_{12} + g_{32}}{\sqrt 2}),\quad
\delta =\frac{1}{2}(g_{11}+ g_{33}-g_{13}-g_{31}).
\end{gather*}

\subsection{The subgroups $Z$ and $H$}
Denote by $Z \subset O(n,1)$ and $H \subset  U(n,1)$ the following subgroups of
block matrices in model $b$:
$$
z=\begin{pmatrix}{}1&0&0\\
-{\gamma}^*&e&0\\ -\frac{{\gamma} {\gamma}^*}{2}&{\gamma}&1\end{pmatrix},
 \,\, {\gamma} \in\RR^{n-1},\quad \text{and} \quad
h=\begin{pmatrix}{}1&0&0\\
-z^*&e&0\\ it-\frac{z z^*}{2}&z&1\end{pmatrix},\,\,
t\in\RR, \,\,  z \in\CC^{n-1},
$$
where $e$ is the unit matrix of order $n{-}1$. They
are the maximal unipotent subgroups of the corresponding groups. The first of them
is isomorphic to the additive group $\RR^{n-1}$; the second one is the
Heisenberg group of dimension ${2n{-}1}$.

Elements of $H$ will be written as pairs  $h=(t,z)$ or $h=(\zeta ,z)$,
where $\zeta =it-\frac{z z^*}{2}$. In this notation, the product of
elements of $H$ is given by the formulas 
$$
(t,z)\,(t',z')=(t+t'-\opn{Im}( zz^{\prime*}),\, z+z')
$$
or
$$
(\zeta,z)\,(\zeta',z')=(\zeta+\zeta'- zz^{\prime*},\, z+z').
$$
Here and in what follows,
$zz^{\prime*}=z_1\OVER{z'_1}+ \ldots +z_{n-1}\OVER{z'_{n-1}}$.

According to the formulas for the transition from one model to 
the other, in model $a$ the subgroup $H$ consists of block matrices 
$h=\begin{pmatrix}\alpha  &\beta  \\ \gamma   &\delta \end{pmatrix}$,
where
\begin{equation}{}\label{1-3}
\begin{gathered}{}
\alpha =\begin{pmatrix}1+\frac{\zeta }{2}&\frac{z}{\sqrt 2}\\
-\frac{z^*}{\sqrt 2}&e\end{pmatrix},\quad
\beta =\begin{pmatrix}-\frac{\zeta }{2} \\  \frac{-z^*}{\sqrt 2} \end{pmatrix},
\\
\gamma =(\frac{\zeta }{2},\,\frac{z}{\sqrt 2}),\quad
\delta =1-\frac{\zeta }{2}.
\end{gathered}
\end{equation}

\subsection{The subgroups $D$, $D_0$, and $D_1$}
Denote by $D$ the subgroup of block-diagonal matrices in model $b$.
It consists of all block-diagonal matrices of the form
$$
d=\diag ( \epsilon ^{-1},u, \epsilon ), \quad
\epsilon \in\RR^*,\quad u\in O(n{-}1)\quad \text{in the case of} \quad O(n,1),
$$
$$
d=\diag (\bar \epsilon ^{-1},u, \epsilon ), \quad
\epsilon \in\CC^*,\quad u\in U(n{-}1) \quad\text{in the case of} \quad U(n,1).
$$

The group $D$ decomposes into the direct product
$$
D=D_0 \times  D_1,
$$
where, in the case of $U(n,1)$, $D_0 \cong U(1) \times U(n{-}1)$ 
is the subgroup of block-diagonal matrices of the form
$$
d=\diag (\epsilon ,u, \epsilon ), \quad \text{where} \quad
|\epsilon|=1,\quad u\in U(n{-}1),
$$
and $D_1 \cong  \RR^*_+ $ is the subgroup of block-diagonal matrices of the form
$$
d=\diag (r^{-1},e, r), \quad \text{where} \quad r>0.
$$
In the case of $O(n,1)$, the subgroups $D_0$ and $D_1$ are defined in a similar way.

\subsection{The subgroups ${{P}}$, ${{P}}_0$, and ${{P}}_1$}

Denote by ${{P}}$ the normalizer of the maximal unipotent subgroup, 
i.e., the maximal parabolic subgroup. Obviously, in model $b$ the subgroup
${{P}}$ coincides with the subgroup of all lower block triangular
matrices and thus can be written as the semidirect product
$$
{{P}}=H\leftthreetimes D.
$$
Further, denote
$$
{{P}}_0=H\leftthreetimes D_0, \quad  {{P}}_1=H\leftthreetimes D_1.
$$

\begin{PROP}{}\label{PROP:1-1} In model $b$, the groups $O(n,1)$ and $U(n,1)$
are algebraically generated by the elements of the subgroup
${{P}}$ and the element
$s=\left(\begin{array}{ccc}0&0&1\\ 0&e&0\\ 1&0&0\end{array}\right)$.
\end{PROP}

\begin{REM*}{} In view of this proposition, in order to describe a representation 
of the group $O(n,1)$ or $U(n,1)$, it suffices to determine, apart from
its restriction to the subgroup 
${{P}}$, only the involution operator $I$ corresponding 
to the element $s$.
\end{REM*}

\part{Representations of the group $O(n,1)$}

\section{Models of irreducible unitary complementary series representations
 of the group $O(n,1)$ and of the special representation of this group}\label{sect:2}

Every complementary series representation is determined by a real number
$\lambda $ from the interval $0< \lambda <n-1$. Let us describe
three models of these representations.

\subsection{Model $A$: realization in a space of functions on the unit sphere
$S \subset \RR^{n}$}
The unit sphere
$$
S = \{\omega \in \RR^n \mid  \quad  |\omega  |^2\equiv {\omega _1}^2+ \ldots + {\omega _n}^2=1 \}
$$
is a homogeneous space of the group $O(n,1)$. Namely, if this group
is realized in matrix model $a$, then its action on $S$ is given by the
formula
$$
\omega g=(\omega \beta + \delta )^{-1}\,(\omega \alpha + \gamma ) \quad
\text{for}\quad
g=\begin{pmatrix}\alpha  &\beta  \\ \gamma   &\delta \end{pmatrix}.
$$
The complementary series  representation with parameter $\lambda $ acts in
the complex Hilbert space $L^{\lambda }$ of functions $f(\omega )$ on $S$ 
with the norm
\begin{equation}{}\label{22-1-0}
\|f\|^2=\int_{S \times S} (1-  \langle \omega , \omega' \rangle)^{- {\lambda } }\,
f(\omega )\,\OVER{f(\omega ')}\,d \omega \,d \omega ',
\end{equation}
where $\langle \omega , \omega ' \rangle = \omega _1\omega '_1+ \ldots + \omega _n\omega '_n$
and $d \omega$ is the invariant measure on $S$.

\begin{REM*}{} Expanding $(1-  \langle \omega , \omega' \rangle)^{- \lambda}$
into a power series, we can also write the norm in the following form:
$$
\|f\|^2 = \frac{1}{\Gamma (\lambda) }
\sum_{k\in\ZZ^n_+} \frac{\Gamma  (\lambda +|k|)}{k!}
        \Bigl| \int_S \omega^k\,f(\omega )\,d \omega \Bigr|^2.
$$
Here $|k|=\sum k_i$, $k!=\prod k_i!$,
$\omega ^k=\prod \omega _i^{k_i}.$
\end{REM*}

The operators $T^{\lambda}_g$,
$g=\begin{pmatrix}\alpha  &\beta  \\ \gamma   &\delta \end{pmatrix}$,
of this representation are given by the formula
\begin{equation}{}\label{22-2}
T^{\lambda}_g f(\omega )=f(\omega g)\,|\omega \beta + \delta |^
{1-n+ \lambda }.
\end{equation}

The group property of these operators is a consequence of the 
following relation for the function
$b(\omega ,g)=|\omega \beta + \delta |$:
$$
b(\omega ,g_1g_2)=b(\omega ,g_1)\,b(\omega g_1 ,g_2) \quad \text{for any}
\quad g_1,g_2\in O(n,1)
$$
(the $1$-cocycle property). The unitarity follows from the relations
$$
1-  \langle \omega , \omega ' \rangle\,=
(1-  \langle \omega g, \omega 'g \rangle)\,
b(\omega ,g)\,b(\omega' ,g)
$$
and
$$
\frac{d(\omega g)}{d \omega }=b^{1-n}(\omega ,g)
\quad \text{for every $g \in O(n,1)$.}
$$

\begin{PROP}{}\label{PROP:2957}
Complementary series representations $T^{\lambda}$ and $T^{\mu}$ 
of the group $O(n,1)$ are equivalent if and only if
$\lambda +\mu=n{-}1$ or $\lambda = \mu $.

An intertwining operator
$R^{\lambda }: L^{\lambda } \to L^{n-1- \lambda }$
is given by the following formula:
$$
(R^{\lambda }f) (\omega ) = \int _{S} (1-\langle \omega ,\omega '\rangle)^{- \lambda }
\,f(\omega ) \,d \omega '.
$$
\end{PROP}

\subsection{Model $B$: realization in a space of functions on the maximal
unipotent subgroup $Z \simeq \RR^{n-1}$}

We construct this model using
matrix realization $b$ of the group
$O(n,1)$. Let us define an action $z\to z\OVER g$ of the group $O(n,1)$ on
$Z$ by the formula
$$
z\OVER g=z',
$$
where $z'\in Z$ is determined from the decomposition $zg=b^+z'$ with $b^+$ 
being an upper block triangular matrix. In the coordinates
$\gamma \in\RR^{n-1}$ on $Z$ and the block coordinates
$g_{ij}$ on $O(n,1)$, this action is given by the following formula:
\begin{equation}{}\label{veron6:1-2}
\gamma \bar g=(-\frac{|\gamma|^2}{2}g_{13}+ \gamma g_{23}+g_{33})^{-1}\,
(-\frac{|\gamma|^2}{2}g_{12}+ \gamma g_{22}+g_{32}),
\end{equation}
where $g_{ij}$ are the entries of a block matrix $g\in O(n,1)$. In particular,
\begin{quote}{}
$
\gamma\OVER g= \gamma + \gamma_0   \quad \text{for}\quad g=z(\gamma _0)\in Z;
$
\\
$
\gamma\OVER g= \epsilon^{-1} \gamma u \quad \text{for}\quad
g=\diag(\epsilon ^{-1},u, \epsilon );
$
\\
$
\gamma\OVER s= -\dfrac{2 \gamma }{|\gamma |^2}.
$
\end{quote}

\begin{REM*}{} The transformation $\gamma \to \gamma \OVER g$ 
has the following geometric interpretation. Let $Y$ be the manifold 
of one-dimensional subspaces in $\RR^{n+1}$ lying in the light cone
$$
2x_1x_{n+1}+x_2^2+ \ldots +x_n^2=0.
$$
The group $O(n,1)$, as a group of linear transformations in
$\RR^{n+1}$, acts transitively on $Y$. We use
the right notation for this action: $y\to y\bar g$.
Note that in another interpretation
$Y$ is the absolute of the $n$-dimensional Lobachevsky space realized
as the collection of one-dimensional subspaces in
$\RR^{n+1}$ lying inside the light cone.

Let us realize $Y\setminus y_0$, where $y_0=(1,0, \dots ,0)$,
as the section of the cone by the hyperplane $x_{n+1}=1$,
i.e., as the set of points in $\RR^{n+1}$ of the form
$$
(-\frac{|\gamma|^2}{2}, \gamma _i, \dots , \gamma _{n-1},1),
$$
where $\gamma =(\gamma _i, \dots , \gamma _{n-1})\in\RR^{n-1}$ and
$|\gamma | = (\sum {\gamma _i}^2)^{1/2}$.  According to this realization,
there is a natural bijection $Y\setminus y_0\to\RR^{n-1}$, so that
the action of the group $G$ on $Y$ induces its action
$\gamma \to \gamma \bar g$ on $\RR^{n-1}$.
We emphasize that this action is not linear.
\end{REM*}

Further, we define a function $\beta (\gamma,g)$ by the formula
\begin{equation}{}\label{veron6:1-3}
\beta (\gamma ,g)=|-\frac{|\gamma|^2}{2}g_{13}+ \gamma g_{23}+g_{33}|, \quad
\gamma \in\RR^{n-1}, \quad g\in G.
\end{equation}

In particular,
\begin{quote}{}
$\beta (\gamma ,g)=1$ for $g\in Z;$
\\
$\beta (\gamma ,g)= |\epsilon\,|$ for $g=\opn{diag}(\epsilon ^{-1},u, \epsilon );$
\\
$\beta (\gamma ,s)=\dfrac{|\gamma |^2}{2}.$
\end{quote}

It follows from the definition that $\beta (\gamma ,g)$ is a $1$-cocycle of the group
$O(n,1)$ with values in $\RR^*$, i.e.,
\begin{equation}{}\label{veron6:1-4}
\beta (\gamma ,g_1g_2)=\beta (\gamma ,g_1)\,\beta (\gamma\bar g_1,g_2) \quad
\text{for any} \quad \gamma \in\RR^{n-1}\quad \text{and}\quad  g_1,g_2\in G.
\end{equation}

In model $B$, the complementary series representation with parameter $\lambda$ acts in the complex 
Hilbert space of functions on $Z\cong\RR^{n-1}$ with the scalar product
\begin{equation}{}\label{veron6:1-70}
\langle f_1,f_2\rangle\,\,\,=\!\!\!\!\!\int\limits_{\RR^{n-1} \times \RR^{n-1}}\!\!\!\!\!\!\!%
| \gamma '- \gamma ''|^{- 2\lambda}
\,f_1(\gamma ')\,\OVER{f_2(\gamma '')}\,d  \gamma '\,d \gamma '',
\end{equation}
where $d \gamma =d \gamma _1 \ldots d \gamma _{n-1}$ is the Lebesgue measure on
$\RR^{n-1}$. The operators of this  representation have the form
\begin{equation}{}\label{veron6:1-8}
T^{\lambda}_g f(\gamma )=f(\gamma \bar g)\,
\beta ^{1-n+\lambda } (\gamma ,g),
\end{equation}
where $\gamma \bar g$ and $\beta (\gamma ,g)$ are given by 
\eqref{veron6:1-2} and \eqref{veron6:1-3}, respectively. In particular,
\begin{align}{}
\label{veron6:273}
&
T^{\lambda}_{z}f(\gamma )=f(\gamma + \gamma _0)\quad \text{for}\quad
z = z(\gamma _0) \in Z;
\\ &
\label{veron6:274}
T^{\lambda}_d f(\gamma )= |\epsilon|^{1-n+ \lambda }\,
f(\epsilon ^{-1} \gamma u)\quad \text{for}\quad
d=\diag(\epsilon  ^{-1},u, \epsilon );
\\ &
\label{veron6:275}
T^{\lambda}_s f(\gamma ) =  f(-\frac{2 \gamma }{|\gamma|^2})\,
\Bigl( \frac{|\gamma|^{2}}{2} \Bigr)^{1-n+ \lambda }
\quad \text{for}\quad
s=\begin{pmatrix}{}\,0&0&1\,\\\,0&e&0\,\\\,1&0&0\,\end{pmatrix}.
\end{align}

The group property of these operators immediately follows from property 
\eqref{veron6:1-4} of the function $\beta (\gamma ,g)$;
the unitarity follows from the relations
\begin{equation}{}\label{veron6:1-5}
d(\gamma \bar g)=\beta^{1-n} (\gamma,g)\,d \gamma \quad
\text{for every} \quad g\in O(n,1),
\end{equation}
where $d \gamma =d \gamma _1 \ldots d \gamma _{n-1}$, and
\begin{equation}{}\label{veron6:1-6}
|{\gamma}-{\gamma'}|^2=|{\gamma}\bar g-{\gamma'}\bar g|^2\, \beta ({\gamma},g)\,\beta ({\gamma'},g)
\end{equation}
for any ${\gamma},{\gamma'}\in\RR^{n-1}$ and $g\in O(n,1)$.

One can easily check \eqref{veron6:1-5} and \eqref{veron6:1-6} 
for elements from $Z$, $D$ and for the element $s$. Then it follows from 
the $1$-cocycle property of $\beta (\gamma ,g)$ that these relations hold for every element
$g\in O(n,1)$.

\begin{PROP}{}\label{PROP:65-34} The functions $F(\omega )$ on $S$ and $f(\gamma )$
on $\RR^{n-1}$ in models $A$ and $B$ are related by
\begin{equation}{}\label{65-14}
f(\omega )=(1+\frac{|\gamma |^2}{2})^{\lambda+1-n}\,F(\tau\gamma ),
\end{equation}
where $\tau \gamma = \omega $ is given by the formulas
$$
\omega_1=\frac{1-\frac{|\gamma |^2}{2}}{1+\frac{|\gamma |^2}{2}},
\quad \omega _i=\frac{2^{1/2}\, \gamma }{1+\frac{|\gamma |^2}{2}},
\quad i=2, \dots ,n.
$$
\end{PROP}

\subsection{Commutative models of the complementary series representations}
Let us describe a new model of a complementary series representation $T^{\lambda}$ 
of the group $O(n,1)$ in which the operators of
the subgroup $Z$ act by multiplication by functions.  We call it the
commutative model with respect to the subgroup $Z$.

This model is obtained by passing from functions $f(\gamma )$ in model $B$ 
to their Fourier transform
$$
\phi(\xi)=\int_{\RR^{n-1}} e^{\,i\langle \xi, \gamma \rangle} f( \gamma )\,d \gamma.
$$

\begin{THM}{}\label{veron6:THM:1} The commutative model of a 
complementary series representation $T^{\lambda}$ is realized
in the Hilbert space  $L^{\lambda}$ of complex-valued functions on
$\RR^{n-1}$  with the norm
\begin{equation}{}\label{veron6:1-13}
\|\phi\|^2=
\tfrac{2^{\,-2\lambda}\, \Gamma (\frac{n-1}{2} - \lambda) }{\Gamma (\lambda)}\,
\int_{\RR^{n-1}} |\xi|^{1-n+ 2\lambda}\,|\phi(\xi)|^2\,d\xi,
\quad
|\xi|=\,\langle \xi,\xi\rangle\,^{1/2}.
\end{equation}
The operators of this representation are given by the formula
\begin{equation}{}\label{veron6:1-14-1}
T^{\lambda}_g \phi(\xi)=\int_{\RR^{n-1}} A^{\lambda}(\xi,\xi',g)\,
\phi(\xi')\,d\xi',
\end{equation}
where
\begin{equation}{}\label{veron6:1-15-1}
A^{\lambda}(\xi,\xi',g)=\int_{\RR^{n-1}}
e^{\,i\,(\langle \xi,\gamma\rangle\,-\,\langle \xi',\gamma\OVER g\rangle)}\,
\beta ^{1-n+\lambda }(\gamma,g)\,d\gamma.
\end{equation}
In particular,
\begin{equation}{}\label{23-23-23}
\begin{gathered}{}T^{\lambda}_{z}\phi(\xi)=
e^{-i\langle \xi,\gamma_0\rangle\,}\,\phi(\xi)\quad \text{for}\quad
z = z(\gamma _0) \in Z;
\\
T^{\lambda}_d \phi(\xi)= |\epsilon|^{\lambda }\,
\phi(\epsilon\, \xi u)\quad \text{for}\quad
d=\diag(\epsilon  ^{-1},u, \epsilon ).
\end{gathered}
\end{equation}
\end{THM}

\begin{proof}{} In the new model, the scalar square is given by the formula
$$
\|\phi\|^2=\!\!\!\!\!\!\!\!%
\int\limits_{\RR^{n-1} \times \RR^{n-1}}\!\!\!\!\!\!\!\!%
R(\xi,\xi')\,\phi(\xi)\,\OVER{\phi(\xi')}\,d\xi\,d\xi',
$$
where
$$
R(\xi,\xi')=\!\!\!\!\!\!\!\!%
\int\limits_{\RR^{n-1} \times \RR^{n-1}}\!\!\!\!\!\!\!%
|\gamma - \gamma '|^{- 2\lambda}
e^{\,i(\langle \xi, \gamma \rangle\,-\,\langle \xi', \gamma' \rangle)}\,d \gamma \,d\xi'=
\delta (\xi-\xi')\!\!\!\int\limits_{\RR^{n-1}} |\gamma |^{- 2\lambda}\,
e^{\,i\langle \xi, \gamma \rangle\,}\,d\xi.
$$
Equation \eqref{veron6:1-13} is a consequence of the following relation:
\begin{equation}{}\label{veron6:17-4}
\int_{\RR^{n-1}} |\gamma|^{-2\lambda}\,e^{\,i\,\langle \xi,\gamma\rangle\,}\,d\gamma =
c_n\frac{2^{\,-2\lambda}\, \Gamma (\frac{n-1}{2} - \lambda) }{\Gamma (\lambda)}.
\,|\xi|^{2\lambda -n+1},
\end{equation}
where the coefficient $c_n$ depends only on $n$. For a proof of
\eqref{veron6:17-4}, see Appendix~2. 

The formulas for the operators of $T^\lambda$ in the new model can be
immediately obtained from the formulas for these operators in the 
original model by passing from functions
$f(\gamma )$ to their Fourier transforms.
\end{proof}

\begin{PROP}{}\label{veron6:PROP:3} In the commutative model of the representation
$T^{\lambda}$, the kernel
$A(\xi,\xi')=A^{\lambda}(\xi,\xi',s)$
of the operator $T^{\lambda}_s$ corresponding to the element
$s=\begin{pmatrix}{}0&0&1\\0&e&0\\1&0&0\end{pmatrix}$ has the following form:
\begin{equation}{}\label{veron6:73-1}
A(\xi,\xi')=2^{1-\lambda}\int_0^{\infty }
\cos(\xi x+\frac{2\xi'}{x})\,x^{2\lambda-2 }\,dx \quad \text{for}\quad n=2,
\end{equation}
\begin{equation}{}\label{veron6:63-1}
A(\xi,\xi')=c_n 2^{- \lambda }\,
\int_0^\infty r^{2\lambda -n}\,|r\xi+2r^{-1}\xi'|^{-\frac{n-3}{2}}\,
J_{\frac{n-3}{2}}(|r\xi+2r^{-1}\xi|)\,dr
\end{equation}
for $n>2$; here $J_{\frac{n-3}{2}}$ is a Bessel function of the second kind.
\end{PROP}

Indeed, since $\gamma \OVER s=\dfrac{-2\gamma \,\,}{|\gamma |^2}$ and
$\beta (\gamma ,s)=\dfrac{|\gamma|^2}{2}$,
it follows from \eqref{veron6:1-15-1} that
$$
A(\xi,\xi')
=  2^{n-1- \lambda }\,\int_{\RR^{n-1}} e^{\,i\,(\langle \xi, \gamma \rangle\,
                +\langle \xi',\frac{2 \gamma }{ |\gamma|^2}\rangle)
}\,
|\gamma|^{2-2n+2\lambda}\,d \gamma .
$$
For $n=2$, this immediately implies \eqref{veron6:73-1}. For $n>2$,
passing to spherical coordinates, we obtain
$$
A(\xi,\xi')=c_n 2^{- \lambda }\,
\int_0^\infty \int_0^{\pi}
e^{\,i\,|r\xi+2r^{-1}\xi'|\cos\phi}\,r^{2\lambda -n}\,
\sin^{n-3}\phi \,d\phi\,dr.
$$
Integrating with respect to $\phi$ yields \eqref{veron6:63-1}.

\begin{REM*}{} For $n=2$, the kernel $A(\xi,\xi')$ 
can be expressed in terms of Bessel functions:
\begin{gather*}{}
A(\xi,\xi')=c\,(\cos(\pi \lambda))^{-1}\,
|\xi'\xi^{-1}|^{1/2}\,[J_{2\lambda -1}(2^{3/2}|\xi\xi'|)-
J_{1-2\lambda}(2^{3/2}|\xi\xi'|)] \quad \text{for} \quad \xi\xi'<0,
\\
A(\xi,\xi')=c\,(\cos(\pi \lambda))^{-1}\,
|\xi'\xi^{-1}|^{1/2}\,[I_{2\lambda -1}(2^{3/2}|\xi\xi'|)-
I_{1-2\lambda}(2^{3/2}|\xi\xi'|)] \quad \text{for} \quad \xi\xi'>0.
\end{gather*}
\end{REM*}

\subsection{Vacuum vectors and spherical functions}
In the space $L^{\lambda}$ of an arbitrary complementary series representation
there exists a vector
$f_{\lambda}$ invariant with respect to the maximal compact subgroup
$K$ of the group $O(n,1)$, and this vector is unique up to 
a factor. Let us call it
the vacuum vector of this representation.

\begin{PROP}{}\label{PROP:75-45} Up to a factor, the vacuum vector
$f_{\lambda} \in L^{\lambda}$ has the form
\begin{align}{}
f_{\lambda}(\omega )&= const  && \text{in model} \quad A,
\\
\label{75-11}
f_{\lambda}(\gamma ) &
= (1+\frac{|\gamma |^2}{2})^{\lambda+1-n}
&& \text{in model} \quad B,
\\
\label{75-17}
f_{\lambda}( \xi )  &
= |\xi|^{-\lambda - \frac{n-1}{2} }\,
 K_{\lambda - \frac{n-1}{2}}(\sqrt2\,|\xi|))
&& \text{in the commutative model};
\end{align}
here $K_{\nu} (x)$ is a Bessel function, see \cite[Vol.~2]{BE}.
\end{PROP}

\begin{proof}{} For model $A$, the assertion follows from the explicit formulas
for the operators of the representation $T^\lambda$. By 
Proposition~\ref{PROP:65-34}, this implies the required assertion for model $B$.

The expression for $f_{\lambda}$ in the commutative model can be obtained from
\eqref{75-11} by the Fourier transform. To this end, it suffices to
apply the following formula:
\begin{equation}{}\label{veron6:17-3}
\int_{\RR^{n-1}}
(1+\,\frac{|\gamma|^2}{a^2})^{-\lambda}\,
e^{\,i\,\langle \xi,\gamma\rangle\,}\,d\gamma
= c_n\frac{2\,a^{\frac{n+1}{2}}}{\Gamma (\lambda)}\,|\xi|^{\lambda - (n-1)/2}\,
K_{\frac{n-1}{2} - \lambda}(a|\xi|).
\end{equation}
A proof of \eqref{veron6:17-3} is given in Appendix~2.
\end{proof}

By definition, the spherical function of a complementary series representation
$T^{\lambda}$ is the function on the group
$O(n,1)$ given by the formula
$$
\psi_{\lambda}(g)=\langle T^{\lambda}_g\vac,\,\vac \rangle\,,
$$
where $\vac$ is the vacuum vector with norm $1$. This function uniquely determines
the representation up to equivalence.

\begin{PROP}{}\label{PROP:81-36} In matrix model $a$ of the group $O(n,1)$,
the spherical function is given by the formula
\begin{equation}{}\label{81-35}
\psi_{\lambda}(g)=c_{\lambda}\,|\beta |^{1-\frac{n}{2}}\,
P^{1-\frac{n}{2}}_{\lambda-\frac{n}{2}}(|\delta |)
\quad \text{for} \quad
g=\begin{pmatrix}\alpha  &\beta  \\ \gamma   &\delta \end{pmatrix},
\end{equation}
where $P^{\mu}_{\nu}$ is a Legendre function and
$|\beta |=(\delta ^2-1)^{1/2}$ is the norm of a vector $\beta \in\RR^n$.
\end{PROP}

\begin{proof}{} Since $\vac = const$ in model $a$, it follows
from the description of the representation  $ T^{\lambda} $ that
$$
\psi_{\lambda}(g)=c_{\lambda}\,
\int_{S \times S} |\omega \beta + \delta |^{\lambda+1-n}\,
(1-\langle \omega ,\omega ' \rangle)^{- \lambda }
\,d \omega\,d \omega'
\quad \text{for} \quad
g = \begin{pmatrix}\alpha &\beta  \\ \gamma &\delta \end{pmatrix}.
$$
Integrating with respect to $\omega '$ yields
$$
\psi_{\lambda}(g)=c_{\lambda}\,
\int_S |\omega \beta + \delta |^{\lambda+1-n}\,d \omega.
$$
In view of the equation $|\beta |^2= \delta ^2-1$, this integral
is a function of $|\delta |$ only. Hence it suffices to calculate it for 
the matrix $g \in O(n,1) $ with the blocks
$$
\alpha =\begin{pmatrix}\ch\tau&0\\0&e \end{pmatrix},\quad
\beta=\begin{pmatrix}\sh\tau\\0 \end{pmatrix},\quad
\gamma =(\sh\tau,0),\quad \delta =\ch\tau.
$$
Passing to the spherical coordinates on $S$, we obtain
$$
\psi_{\lambda}(g)=c_{\lambda}\,
\int_0^\pi (\sin \alpha \sh\tau+\ch\tau)^{\lambda+1-n}\,
\sin^{n-2}\alpha \,d \alpha .
$$
The obtained integral can be expressed in terms of the Legendre function
$P^\mu_\nu(x)$ according to the following formula, see
\cite[Vol.~1, p.~156(7)]{BE}:
$$
P^\mu_\nu(\ch\tau)
        =\frac{
        2^\mu(\sh\tau)^{- \mu }
        }{\pi^{1/2}\,\Gamma (\frac12-\mu)
        }
\,\int_0^\pi (\sin \alpha \sh\tau+\ch\tau)^{\mu+\nu}\,
(\sin \alpha )^{-2\mu}\,\,d \alpha .
\qed
$$
\let\qed\null
\end{proof}

\setcounter{REM}{0}

\begin{REM}{} The fact that the spherical functions of two equivalent
representations, i.e., representations with parameters
$\lambda $ and $n{-}1{-}\lambda $, coincide is a consequence of the following relation for
the Legendre function:
$P^{\mu }_{-\nu -1} (\ch \tau)  = P^{\mu }_{\nu} (\ch \tau)$.
\end{REM}

\begin{REM}{} The expression for $\psi_{\lambda}(g)$ in matrix model
$b$ of the group $O(n,1)$ is obtained by replacing
$|\delta |$ in \eqref{81-35} with
$\frac12 \, |g_{11} +g_{33}  - g_{13} -  g_{31}|$.
\end{REM}

\setcounter{REM}{0}

\subsection{Embedding into the tensor product}

\begin{PROP}{}\label{veron6:PROP:5} For any real positive numbers
$\lambda _1, \dots ,\lambda _m$ with
$\sum \lambda _i<n-1$ there exists an isometric embedding 
$$
\tau:  L^{\lambda}\to \bigotimes _{i=1}^m L^{\lambda _i},\quad
\lambda =\sum \lambda _i,
$$
commuting with the action of the group $O(n,1)$. In model $B$,
it is given by the formula
\begin{equation}{}\label{veron6:31-21-1}
\tau f(\gamma^1, \dots ,\gamma^m)=f(\gamma^1)
\prod_{i=2}^m \delta (\gamma ^1- \gamma^i);
\end{equation}
in the commutative model, it is given by the formula
\begin{equation}{}\label{veron6:31-21-2}
\tau \phi(\xi _1, \dots ,\xi _m)=\phi(\xi _1+ \ldots +\xi_m).
\end{equation}
\end{PROP}

\begin{proof}{}
Let $\langle \,,\rangle\,$ and $\langle \,,\rangle_m$ be the scalar products
in the spaces $L^{\lambda}$ and $\bigotimes _{i=1}^m L^{\lambda _i}$,
respectively, in model $B$. Then, obviously,
$$
\langle \tau f,\tau f\rangle_m\,\,=\!\!\!\!\!\!\int\limits_{\RR^{n-1} \times \RR^{n-1}}
\prod_{i=1}^m| \gamma '- \gamma ''|^{- 2\lambda_i }
\,f(\gamma ')\,\OVER{f(\gamma '')}\,d  \gamma '\,d \gamma ''=\,\langle f,f\rangle.
$$
Thus the map $\tau$ is isometric. It is obvious that
it commutes with the action of $O(n,1)$.

The commutative model is obtained by passing from functions $f$ in model $B$
to their Fourier transforms. Let $\phi (\xi ) $ be the image of 
$f(\gamma )$ under  this transform. Let us check that the image of
$\tau f(\gamma ^1, \dots ,\gamma ^m)$
is $\phi (\xi _1+ \ldots +\xi _m)$. Indeed, according to
\eqref{veron6:31-21-1}, this image is equal to
$$
\int\limits_{(\RR^{n-1})^m}!\!\!\!\! f(\gamma _1)\,\prod_{i=2}^m
\delta  (\gamma _1- \gamma _i)\,\prod_{i=1}^m e^{\,i\,\langle \xi_i, \gamma_i\rangle }\,d\gamma _i
\,\, =\,\,
\phi(\xi _1+ \ldots +\xi_m). \qed
$$
\let\qed\null
\end{proof}

\begin{REM*}{} This embedding leads to one of the models
of an irreducible  nonlocal unitary representation of the current group 
$O(n,1)^X$, i.e., the group of measurable maps
$X \mapsto O(n,1)$, where $X$ is a measurable space.
\end{REM*}

\subsection{Models of the special representation of the group $O(n,1)$, $n>2$}
The special irreducible unitary representation of the group  
$O(n,1)$, $n>2$, is obtained from the complementary series representations 
by passing to the limit as $\lambda \to 0$ or $\lambda \to n{-}1$.

As $\lambda \to 0$, the norm $\|f\|_{\lambda }$ in the space of a
complementary series representation degenerates on a subspace $L^0$ of codimension 
$1$. The special representation is realized in the subspace $L^0$ with the norm
$$
\|f\|^2=\lim_{\lambda \,\to\, 0}\frac{d\,\|f\|^2_{\lambda }}{d \lambda }.
$$
For $\lambda = n{-}1$, the norm $\| f \|_{\lambda }$ degenerates on 
a one-dimensional invariant subspace $L_0$ of the group  $O(n,1)$.
The special representation is realized in the quotient 
$L^{n-1}/L_0$ with the norm
$$
\|f\|^2=\lim_{\lambda \,\to\, n-1}\|f\|^2_{\lambda }.
$$
The operators of the special representation can be obtained from the
operators of complementary series representations by passing to the limit as
$\lambda \to 0$ and $\lambda \to n{-}1$.
The unitary representations of  $O(n,1)$
in the spaces $L^0$ and $L^{n-1}/L_0$ defined in this way are irreducible and equivalent.

Since each of the spaces $L^0$ and $L^{n-1}$ can be defined in three
different models, we thus obtain six different models of the special representation.

Let us explicitly describe the spaces
$L^0$ and $L^{n-1}$ in models $A$ and $B$.

\subsubsection{Model $A$} 
Here $L^0$ is the subspace of functions $f(\omega )$
on the unit sphere $S=S^{n-1}$ in $\RR^n$ satisfying the condition
$$
\int_S f(\omega ) \,d \omega=0.
$$
The norm in the space $L^0$:
\begin{equation}{}\label{22-11}
\|f\|^2=-\int_{S \times S}
\log (1-  \langle \omega , \omega' \rangle)\,
f(\omega )\,\OVER{f(\omega ')}\,d \omega \,d \omega '.
\end{equation}
The norm in the space $L^{n-1}$ of functions on $S$:
\begin{equation}{}\label{22-11-n-1}
\|f\|^2=\int_{S \times S}
(1-  \langle \omega , \omega' \rangle)^{-n+1}\,
f(\omega )\,\OVER{f(\omega ')}\,d \omega \,d \omega '.
\end{equation}
The intertwining operator $R:L^0\to L^{n-1}$:
\begin{equation}{}\label{22-11-R}
(Rf)(\omega ) = \int_{S}
\log (1-  \langle \omega , \omega' \rangle)\,
f(\omega' ) \,d \omega '.
\end{equation}

\begin{REM*}{} In view of \eqref{22-11-R}, the norm in $L^0$ 
can be written as an integral over $S$:
$$
\|f\|^2 =\int_S (Rf)(\omega )\,\OVER{f(\omega )}\,d \omega .
$$
\end{REM*}

\subsubsection{Model $B$} Here $L^0$ is the subspace of functions
$f(\gamma) $ on $\RR^{n-1}$ satisfying the condition
$$
\int_{\RR^{n-1}} f(\gamma  ) \,d \gamma=0.
$$
The norm in the space  $L^0$:
\begin{equation}{}\label{22-121}
\|f\|^2=-\int_{\RR^{n-1} \times \RR^{n-1}}
\log |\gamma - \gamma '|\,
f(\gamma )\,\OVER{f(\gamma ')}\,d \gamma \,d \gamma '.
\end{equation}
The norm in the space $L^{n-1}$ of functions on $\RR^{n-1}$:
\begin{equation}{}\label{22-121-n-1}
\|f\|^2=\int_{\RR^{n-1} \times \RR^{n-1}}
|\gamma - \gamma '|^{-2n+2}\,
f(\gamma )\,\OVER{f(\gamma ')}\,d \gamma \,d \gamma '.
\end{equation}
The intertwining operator $R:L^0\to L^{n-1}$:
\begin{equation}{}\label{22-121-R}
\|f\|^2=-\int_{\RR^{n-1}} \log |\gamma - \gamma '|\,f(\gamma' )\,d \gamma '.
\end{equation}

\subsubsection{Commutative model}

\begin{THM}{}\label{THM:5-5-5} In the commutative models corresponding to
$\lambda =0$ and $\lambda=n{-}1$, the special representation of
the group $O(n,1)$, $n>2$, is realized in the Hilbert spaces 
of functions $f(\xi)$ on $\RR^{n-1}$ with the norms
\begin{equation}{}\label{52-13}
\|f\|^2=\int_{\RR^{n-1} } |f(\xi)|^2\,|\xi|^{1-n}\,d\xi
\quad \text{and}\quad
\|f\|^2=\int_{\RR^{n-1} } |f(\xi)|^2\,|\xi|^{n-1}\,d\xi,
\end{equation}
respectively.
\end{THM}

\begin{REM*}{} In the exceptional case of the group
$O(2,1)$, each of the spaces
$L^0$ and $L^{n-1}/L_0$ splits into two irreducible subspaces, and we obtain
not one but two special representations, 1) and 2), of this group.
In model $A$, at each endpoint of the interval
$[0,1]$, they are realized in the subspaces of functions on the unit circle 
that are the boundary values of 1) analytic and 2) antianalytic functions in the disk
$|z|<1$. This case reduces to the case of the group
$U(1,1)$ considered in~\S\,\ref{sect:6}.
\end{REM*}

\subsection{The nontrivial $1$-cocycle associated with the special representation}
\label{sect:b(g):O(n,1)}
The space $L$ of the special representation has a nontrivial $1$-cocycle,
i.e., a map $b: O(n,1)\to L$ that satisfies the condition
$$
b(g_1g_2) = b(g_1) + T_{g_2} b(g_1) \quad \text{for any $g_1,g_2 \in O(n,1)$}
$$
and cannot be written in the form $b(g) = T_g\,\xi - \xi $, $\xi \in L$.

Let us explicitly describe the $1$-cocycle $b$ in different realizations
of this representation.

\begin{THM}{}\label{THM:5-5-10}\
\begin{itemize}
\item[$\circ$]
The $1$-cocycle $b:\, O(n,1)\to L^0$ is given by the formula
$$
        b(g)=T_g \xi - \xi ,
$$
where $\xi$ is a vector that does not belong to $L^0$ and is invariant under
the maximal compact subgroup $U$ of $O(n,1)$.
The cocycle $b$ is nontrivial.

\item[$\circ$]
The $1$-cocycle $b':\, O(n,1)\to L^{n-1}$ is given by the formula
$$
        b'(g)=R\,b(g),
$$
where $R:L^0\to L^{n-1}$ is an intertwining operator.
\end{itemize}
\end{THM}

The cocycle $b$ is nontrivial, because in  $L^0$
there are no $U$-invariant vectors.

Explicit expressions for $\xi $ in different models can be obtained from
the formulas for the vacuum vector in the spaces of complementary series
representations by passing to the limit as $\lambda \to 0$. Namely,
\begin{align*}{}
\xi&=const && \text{in model} \quad A,
\\
\xi &=\Bigl(1+\frac{|\gamma |^2}{2}\Bigr)^{1-n}
&& \text{in model} \quad B,
\\
\xi&= |\xi|^{-\frac{1-n}{2} }\,
 K_{\frac{1-n}{2}}(\sqrt2\,|\xi|))
&&\text{in the commutative model};
\end{align*}
here $K_{\rho}(x)$ is a Bessel function, see \cite[Vol.~2]{BE}.
In particular, $\xi = e^{\,- \sqrt2\,|\xi |}$ for $n=2$.

\begin{REM*}{} For $\lambda = n{-}1$, a $K$-invariant vector is invariant
with respect to the whole group $O(n,1)$.
\end{REM*}

Further, we have
\begin{align*}{}
b'(\omega ,g) & = \int_S \log (1-  \langle \omega , \omega' \rangle)\,
| \omega ' \beta  + \delta |^{1-n} \,d \omega '
&& \text{in model} \quad A,
\\
b'(\omega ,g) & = \int_{\RR^{n-1}} \log (| \gamma - \gamma '|)\cdot
|-\tfrac{|\gamma' |^2}{2} \,g_{13} + \gamma' \,g_{23} + g_{33}|^{1-n}\,
d \gamma '
&& \text{in model} \quad B.
\end{align*}

\subsection{Restriction to the maximal parabolic subgroup
$P$}\label{sect:2-8(New)}

\begin{PROP}{}\label{PROP:28-28-28} The restrictions of the 
complementary series representations
$T^{\lambda }$ of the group $O(n,1)$, $n>2$,
and of the special representation $T^0\simeq T^{n-1}$ of this group
to the maximal parabolic subgroup $P$ are irreducible and
pairwise equivalent. 
In the commutative model, an intertwining operator $\tau : L^{\lambda }\to L^0$
is given by $(\tau \phi )(\xi ) = |\xi |^{\lambda }\phi (\xi )$.

The unitary representation of  $P$ arising in this way is a
special representation of this group.
\end{PROP}

\begin{proof}{} In the commutative model, the representations
$T^{\lambda }$, including the special representation
$T^{0}$, are realized in the Hilbert spaces
$L^{\lambda }$ with the norm
$$
\| \phi \|^2_{\lambda } = \int_{\RR^{n-1}} |\phi (\xi )|^2\,
|\xi |^{\,1-n+2 \lambda } \,d \xi .
$$
It follows from the explicit formulas \eqref{23-23-23}
for the operators of the subgroup ${{P}}$ that $\tau $ is an intertwining operator
$L^{\lambda } \to L^0$.
\end{proof}

\begin{REM*}{} In order to describe the extension of this special representation
of  ${{P}}$ to an arbitrary complementary series representation
$T^{\lambda }$, it suffices to determine
only the operator $T_s^{\lambda }$.
\end{REM*}

\subsection{The study of the separability
of the special representations of the groups
$O(n,1)$ and $P$ from the identity representations in the Fell topology}

By definition, a unitary representation $T$ of a locally compact group $G$
in a Hilbert space $H$ contains
almost invariant vectors if
for every compact set $K \subset G$ and every $\epsilon > 0$ there exists a unit vector
$h \subset H$ such that
$\| T_k h - h \| < \epsilon $ for all $k \in K$.

The set of all representations of $G$ that contain almost invariant vectors coincides
with the intersection of all neighborhoods of the identity representation 
of $G$ in the Fell topology on the set of unitary representations of $G$.

\begin{PROP}{}\label{PROP:27-28-29-30-31}
\
\begin{enumerate}
\item
The special representation $T$ of the group
$O(n,1)$ does not contain almost invariant vectors.
\item
The restriction $T'$ of the representation $T$ to the maximal parabolic subgroup
$P$ contains almost invariant vectors.
\end{enumerate}
\end{PROP}

\begin{proof}{} 1) The representation $T$ of the group $O(n,1)$.
Let $K$ be the maximal compact subgroup of
$G = O(n,1)$. Since $T$ is a special representation, the space 
$H$ does not contain $K$-invariant vectors, so that
$$
\int_K \langle T_kh,h\rangle\,dk = 0
$$
for every unit vector $h \in H$. Integrating the equation
$$
\| T_k h - h \|^2 = 2 - \langle T_kh,h\rangle - \langle h,T_kh\rangle
$$
over $K$, we obtain
$\int _K \| T_k h - h \|^2\,dk = 2$. Hence 
$T$ does not contain almost invariant vectors.

2) The representation $T'$ of the subgroup $P$.
Let us use the commutative model of the special representation of
$O(n,1)$. Set
$$
h(\xi ) = \chi_{\eta } (|\xi |) \cdot |\xi |^{\,\lambda }
\qquad (\lambda > 0),
$$
where $\chi _{\eta }$ is the characteristic function of the interval
$( 0 , \eta )$. Obviously, 
$h(\xi ) \in H$ for any $\eta >0 $ and $\lambda >0 $.
Let us check that for every compact set $K \subset P$ and every $\epsilon >0$,
and for all sufficiently small $\eta $ and $\lambda $,
\begin{center}{}
$\| T'_kh - h \|^2 < \epsilon \cdot \|h \|^2$ for $k \in K$.
\end{center}
Let us write elements $k \in K$ in the form $k = z(\gamma )\cdot \diag(r^{-1},u,r)$.
Then
$$
(T'_kh)(\xi )
= e^{\, - i\,\langle \xi ,\gamma \rangle\,}
\cdot \chi _{\eta } (|r \xi |)\,|r\xi |^{\lambda }.
$$
Without loss of generality we may assume that
$| r | < 1$ on $K$. Under this assumption, we have
$\chi _{\eta }(|r \xi |)\,\chi _{\eta }(|\xi |) = \chi^2 _{\eta }(|\xi|)$
for $|\xi |< \eta $, whence
\begin{align*}{}
\| T'_kh - h \|^2 &
= 2\int _{|\xi | < \eta }
(h(\xi )  -  (T'_k h)(\xi ))\,h(\xi )\,|\xi |^{\,1-n}\,d \xi
\\  &
= 2\int_{|\xi | < \eta }
(1-|r|^{\lambda } \,e^{\, - i\,\langle \xi ,\gamma \rangle\,})\,
h^2(\xi )\,|\xi |^{\,1-n}\,d \xi.
\end{align*}
This immediately implies the required assertion.
\end{proof}

\begin{REM*}{} The proposition remains valid for the special representations
$T$ of the group $U(n,1)$ considered in~\S\,\ref{sect:6}.
\end{REM*}

\subsection{Remark on representations in real spaces}
In models $A$ and $B$, the operators of the representation $T^\lambda$ 
preserve the spaces
of real-valued functions. Therefore, the spaces of complementary series 
representations and the space of the special representation, regarded as 
linear spaces over $\RR$, split into direct sums of two invariant
real Hilbert subspaces: the spaces of functions with real and purely 
imaginary values, respectively. The representations of the group
$O(n,1)$  in these real Hilbert spaces are irreducible and equivalent.

In the commutative model, the representation space splits into a direct
sum of invariant subspaces over $\RR$: the space of functions
$f(\xi) +\OVER{f(-\xi)}$ and the space of functions
$f(\xi) -\OVER{f(-\xi)}$.

\subsection{Application of the commutative model of the special representation
of the group $O(n,1)$ to constructing an irreducible unitary representation
of the current group
$O(n,1)^X$}

The role of the commutative model of the special representation of 
$O(n,1)$ manifests itself in that it leads to a natural construction of an irreducible 
unitary representation of the current group $O(n,1)^X$, i.e., the group
of measurable maps $X\to O(n,1)$, where
$X$ is a space with a finite measure $m(x)$,  with pointwise multiplication.

We start from the space $L^2(|\xi|^{1-n}\,d\xi)$, $\xi\in\RR^{n-1}$,
in the  commutative model of the special representation
of $O(n,1)$, and construct a 
$\sigma $-finite measure $\nu$ in the space $\Phi^*=\Phi(X)$ of vector 
distributions $\xi(x)=(\xi_1(x), \dots ,\xi_{n-1}(x))$
satisfying the following invariance properties:
\begin{enumerate}
\item The measure $\nu$ is invariant under the action of the group
$O(n-1)^X$, i.e.,
$$
d\nu(\xi u)=d\nu(\xi )\quad \text{for every} \quad u(x)\in O(n-1)^X.
$$
\item The measure $\nu$ is projectively invariant 
under the multiplication
by any bounded Borel function $\epsilon (x)\in (\RR^*)^X$ for which the 
integral $\int_X |\epsilon (x)|\,d m(x)$ converges, namely,
$$
d(\epsilon \xi)=e^{\int_X \log |\epsilon (x)|\,d m(x)}\,d\nu(\xi).
$$
\end{enumerate}
In particular, $d\nu(\epsilon \xi)=d\nu(\xi)$ if
$\int_X  \log |\epsilon (x)|\,d m(x)=0.$

In view of these properties, we call
$\nu$ the Lebesgue measure on $\Phi^*$.

We construct a unitary representation of the group ${{P}}^X$ in the Hilbert space
$L^2(\nu)$ of functions $f(\xi)$ on $\Phi^*$. By analogy 
with the commutative model
of the special representation of
$O(n,1)$, elements $\gamma (x)\in Z^X$ act as multiplicators:
$$
(T_{\gamma }\,f) (\xi)=e^{i\langle \xi, \gamma  \rangle\,}\,f(\xi),
$$
and the action of elements $d=(\epsilon ^{-1},u, \epsilon )\in D^X$ is given  by the formula
$$
(T_d\,f) (\xi)=e^{\frac12\, \int_X \log |\epsilon (x)|\,d m(x)}\,f(\epsilon \xi u).
$$

It is easy to check that these operators generate a representation of
the whole group ${{P}}^X$. The unitarity of
$T_{\gamma }$ is obvious, and the unitarity of 
$T_d$ follows from the Lebesgue properties of the measure $\nu$.

The unitary representation of the group ${{P}}^X$ constructed in this
way is irreducible. As proved in \cite{V-G-2005,V-G-indag}, it can be
extended to a unitary representation of the whole group
$O(n,1)^X$.

\part{Representations of subgroups of the group $U(n,1)$}

\section{Representations of the Heisenberg group $H$ and their
extensions to the groups
${{P}}_0=H\leftthreetimes D_0$ and ${{P}}=H\leftthreetimes D$}\label{sect:3}

This and the next sections are devoted to the description of 
the complementary series representations and special representations of the group
$U(n,1)$. We begin with a construction of representations
of the maximal unipotent subgroup, i.e., the Heisenberg group
$H \subset U(n,1)$, and its normalizer, the parabolic subgroup
${{P}}$.

\subsection{Description of the irreducible unitary representations of the Heisenberg
group}
It is known (see, e.g., \cite{Perelomov}) that 
the irreducible unitary representations of the Heisenberg group $H$
break into two classes: one-dimensional representations, depending on
$2(n{-}1)$ real parameters, and infinite-dimensional ones, depending
on one real parameter $\rho \ne 0$ (the Planck constant). Let us
describe the Bargmann model
of the infinite-dimensional irreducible representations.

The irreducible representation with parameter $\rho >0$
is realized in the Hilbert space
$\Cal H(\rho )$ of entire analytic functions $f(z)$
in $z=(z_1, \dots ,z_{n-1})\in \CC^{n-1}$ with the norm
\begin{equation}{}\label{444-444}
\| f \|^{2} = |\rho |^{n-1}
\int_{\CC^{n-1}} |f(z)|^{2}\,e^{\, - |\rho| \,|z|^2}\,d\mu(z),
\end{equation}
where $|z|^2=zz^*=|z_1|^2+ \ldots +|z_{n-1}|^2$ and $d\mu(z)$ is the Lebesgue measure
on $\CC^{n-1}$ normalized by the condition 
$$
\int_{\CC^{n-1}} e^{\, - |z|^2}\,d\mu(z)\,\, =\,\, 1.
$$
The operator $T_h$ 
corresponding to an element $h = (t_0,z_0)\in H$ has the form
\begin{equation}{}\label{746-24}
T_{h} f(z) = e^{\,\rho  \,(it_0 - \frac12\,|z_0|^2 - zz_0^*)}\,f(z+z_0).
\end{equation}
It is easy to check that these operators are unitary and satisfy the group property.

The irreducible representation with parameter $\rho <0$ is realized in the Hilbert 
space $\Cal H(\rho )$  of entire antianalytic functions
$f(z)$ in $z=(z_1, \dots ,z_{n-1})\in \CC^{n-1}$
with the same norm \eqref{444-444}.

The operator $T_h$ 
corresponding to an element $h=(t_0,z_0)\in H$ has the form
\begin{equation}{}\label{746-25}
T_{h} f(z) = e^{\,i\rho  t_0  - |\rho|\,(\frac12\,|z_0|^2 + z_0 z^*)}\,f(z+z_0).
\end{equation}

Note that the transition from the space $\Cal H(\rho )$ to the space
$\Cal H(-\rho )$ reduces to replacing
$f(z)$ by $\OVER{f(z)}$.

\begin{PROP}{}\label{PROP:254-47} 
The representations of the Heisenberg group $H$ in the spaces
$\Cal H(\rho )$, $\rho \ne 0$, are irreducible and pairwise nonequivalent;
they exhaust all infinite-dimensional irreducible unitary representations
of $H$.
\end{PROP}

\begin{PROP}{}\label{PROP:254-48} 
The monomials $z^k=z_1^{k_1} \ldots z_{n-1}^{k_{n-1}}$ and
$\OVER z^k=\OVER z_1^{k_1} \ldots \OVER z_{n-1}^{k_{n-1}}$,  $k\in\ZZ^{n-1}_+$,
form orthogonal bases in the spaces
$\Cal H (\rho )$ for $\rho >0$ and $\rho <0$, respectively; their norms equal
$$
\|z^k\|^2=\|\OVER z^k\|^2=k!\, |\rho |^{-|k|},
$$
where
$$
|k|=k_1+ \ldots + k_{n-1}, \quad k!=k_1! \ldots k_{n-1}!.
$$
In this basis, the operators $T_h$, $h=(t_0,z_0)$,  for $\rho >0$ have the form
\begin{equation}{}\label{254-50}
T_h z^k=   e^{\,\rho  \,(it_0 - \frac12\,|z_0|^2 )}\,
\sum_{m\in\ZZ^{n-1}_+} a_{km}(\rho ,z_0)\,z^m,
\end{equation}
where
\begin{equation}{}\label{254-51}
a_{km}(\rho ,z_0)=k!\,\sum_s
\frac{(- \rho )^{|m-s|}\,z_0^{k-s}\,\OVER z_0^{m-s}}{s!\,(k-s)!\,(m-s)!}.
\end{equation}
For $\rho <0$, the formula for $T_h$  is analogous.
\end{PROP}

\begin{REM*}{} The change of variables $z\to z'= |\rho| ^{1/2}z$
determines isomorphisms of Hilbert spaces
$\Cal H(\rho )\to\Cal H(1)$ and
$\Cal H(\rho )\to\Cal H(-1)$ for $\rho >0$ and $\rho <0$, respectively.
Thus every irreducible infinite-dimensional representation of
$H$ with parameter $\rho $ can be realized in the Hilbert space of entire
analytic (antianalytic) functions with the norm
$$
\| f \|^{2} =
\int_{\CC^{n-1}} |f(z)|^{2}\,e^{\, - |z|^2}\,d\mu(z).
$$
In this realization, the operators of this representation have the form
\begin{gather*}{}
T_{h} f(z) = e^{\,\rho  \,(it_0 - \frac12\,|z_0|^2) - \rho ^{1/2} zz_0^*}\,
f(z+ \rho ^{1/2} z_0) \quad \text{for}\quad \rho >0;
\\
T_{h} f(z) = e^{\,i\rho  t_0  - \,\frac12\,|\rho||z_0|^2
- |\rho|^{1/2} z_0 z^*}\,f(z+|\rho|^{1/2} z_0)
        \quad \text{for}\quad \rho <0.
\end{gather*}
\end{REM*}

\subsection{Another expression for the norm in the spaces $\Cal H(\rho )$}

\begin{PROP}{}\label{PROP:796-45} 
For every entire analytic function
$f(z)$ on $\CC^{n-1}$, for every $\rho >0$,
\begin{equation}{}\label{333-337}
\int\limits_{\CC^{n-1}}|f(z)|^2\,e^{\,- \rho |z|^2}\,d\mu(z)
= \rho^{n-1}\!\!\!\!\!\!\!\!%
\int\limits_{\CC^{n-1} \times \, \CC^{n-1}} \!\!\!\!\!\!%
f(z)\,\OVER{f(z')}\, e^{\, \rho\,(z'z^*-|z|^2-|z'|^2)}\,d\mu(z)\,d\mu(z').
\end{equation}
Analogously, for every entire antianalytic function
$f(z)$ on $\CC^{n-1}$,
\begin{equation}{}\label{433-337}
\int\limits_{\CC^{n-1}}|f(z)|^2\,e^{\,- |\rho| |z|^2}\,d\mu(z)
= |\rho|^{n-1}\!\!\!\!\!\!\!\!\!\!%
\int\limits_{\CC^{n-1} \times \, \CC^{n-1}} \!\!\!\!\!\!%
f(z)\,\OVER{f(z')}\, e^{\, |\rho|\,(z(z')^*-|z|^2-|z'|^2)}\,d\mu(z)\,d\mu(z').
\end{equation}
\end{PROP}

\begin{proof}{} Let us prove \eqref{333-337}. It suffices to prove this equation
for all monomials $f(z)=z^k$. Let  $I_1$ and $I_2$ be the 
left- and right-hand sides of  \eqref{333-337}, respectively, for $f(z)=z^k$.
Obviously, $I_1=k!\, \rho ^{-|k|-n+1}$. On the other hand, let us
substitute the power series expansion 
$e^{\rho z_0z^*}=\sum_l \frac{\rho ^{|l|}}{l!}\,(z')^l\OVER z^l$
of $e^{\rho z'z^*}$
into the expression for $I_2$.
Observing that the terms of the resulting sum vanish for
$l\ne k$, we have
$$
I_2=\frac{\rho ^{|k|+n-1}}{k!}\,\Bigl|\int_{\CC^{n-1}}z^k\OVER z^k\,
e^{- \rho \,|z|^2}\,d\mu(z)\Bigr|^2=\frac{\rho ^{|k|+n-1}}{k!}\,I_1^2=
k!\, \rho ^{-|k|-m+1}.
$$
Thus $I_1=I_2$. The proof of \eqref{433-337} is similar.
\end{proof}

\begin{COR*}{} The norms \eqref{444-444} in the spaces $\Cal H(\rho )$
for $\rho >0$ and $\rho <0$ coincide with the norms
\begin{equation}{}\label{454-444}
\|f\|^2=\rho^{2n-2}\int\limits_{\CC^{n-1} \times \, \CC^{n-1}} \!\!\!\!\!\!%
f(z)\,\OVER{f(z')}\, e^{\, \rho\,(z'z^*-|z|^2-|z'|^2)}\,d\mu(z)\,d\mu(z')
\end{equation}
and
\begin{equation}{}\label{464-444}
\|f\|^2=|\rho|^{2n-2}\int\limits_{\CC^{n-1} \times \, \CC^{n-1}} \!\!\!\!\!\!%
f(z)\,\OVER{f(z')}\, e^{\, |\rho|\,(z(z')^*-|z|^2-|z'|^2)}\,d\mu(z)\,d\mu(z'),
\end{equation}
respectively.
\end{COR*}

\begin{REM*}{} On the space of all functions
$f(z)$, the norm \eqref{444-444} is strictly positive, while the norms
\eqref{454-444} and \eqref{464-444} are degenerate.
\end{REM*}

\subsection{Extension of a representation of the Heisenberg group $H$ to the 
group ${{P}}_0=H\leftthreetimes D_0$}

\begin{PROP}{}\label{PROP:713-27} 
The infinite-dimensional representations of 
the Heisenberg group $H$ can be extended to the group ${{P}}_0=H\leftthreetimes D_0$.
Namely, in the Bargmann model $\Cal H(\rho )$,
the operators $T_d$ corresponding to elements
$d=\diag(\epsilon ,u, \epsilon )\in D_0$ are given by the formula
\begin{equation}{}\label{285-14}
T_d f(z)=f(\epsilon^{-1}  z u).
\end{equation}
\end{PROP}

\begin{COR*}{} The restriction of the representation of the group ${{P}}_0$ 
in the space $\Cal H(\rho )$ to the subgroup $U(n{-}1)\subset D_0$ 
splits into a direct sum of pairwise nonequivalent subspaces:
subspaces of homogeneous polynomials in
$z\in\CC^{n-1} $ for $\rho >0$ and subspaces of homogeneous polynomials in
$\OVER z$ for $\rho <0$.

In particular, the representation space contains
a $U(n{-}1)$-invariant vector, namely,
$f(z) = const$; this vector is unique up to a factor.
\end{COR*}

The explicit expression \eqref{444-444} for the norm and the expressions
for the operators of the group
${{P}}_0$ in $\Cal H(\rho )$ imply the following proposition.

\begin{PROP}{}\label{PROP:222-333} For every $r>0$, the map
$f(z)\to f(r^{-1}\,z)$ is an isomorphism of Hilbert spaces
$\Cal H(\rho )\to\Cal H(r^{-2}\,\rho)$, which send the operators 
$T_{b}$, $b \in {{P}}_0$, on $\Cal H(\rho )$ to the operators
$T_{d_1\,b\,d_1^{-1}}$ on $\Cal H(r^{-2}\,\rho)$,
where $d_1=\diag(r^{-1},e,r)$.
\end{PROP}

\subsection{Direct integrals of representations of the group
${{P}}_0$ and their extensions to representations of the parabolic
group ${{P}}={{P}}_0\leftthreetimes D_1=H\leftthreetimes D$}
Let us associate with each real number $\lambda \ge 0$ the following direct integrals
of the Hilbert spaces $\Cal H(\rho )$:
$$
\Cal H^{\lambda }_-=\int_{- \infty} ^0  \Cal H(\rho )\,| \rho |^{\lambda -1}\,
d \rho ; \quad
\Cal H^{\lambda }_+=\int^{+\infty} _0  \Cal H(\rho )\, \rho ^{\lambda -1}\,
d \rho .
$$
In more detail,
$\Cal H^{\lambda }_+$ is the Hilbert space of functions
$f(\rho ,z)$ on $\RR_+ \times \CC^{n-1}$ that are entire analytic in
$z$ with the norm
\begin{equation}{}\label{212-13}
\|f\|^2=\int_0^\infty \Bigl(
\int_{\CC^{n-1}} |f(\rho ,z)|^{2}\,e^{\, - |\rho| \,|z|^2}\,d\mu(z)\Bigr)\,
\rho ^{n+\lambda -2}\,d \rho.
\end{equation}
The space $\Cal H^{\lambda }_-$ is defined in a similar way.

The unitary representations of the group ${{P}}_0$ in $\Cal H(\rho )$
induce unitary representations of this group in 
$\Cal H^{\lambda }_{\pm}$.

\begin{PROP}{}\label{PROP:211-05} 
The representations of the group ${{P}}_0$ in the spaces
$\Cal H^{\lambda }_{\pm}$ can be extended to unitary representations of the group 
${{P}}={{P}}_0\leftthreetimes D_1$. Namely, the operators $T_d$ corresponding
to elements $d\in D_1$ are given by the formula 
\begin{equation}{}\label{273-66}
T_d f(\rho ,z)=f(r^2 \rho ,r^{-1}z)\,r^{ \lambda }\quad \text{for}\quad
d=\diag(r^{-1},e,r),\quad r>0.
\end{equation}
\end{PROP}

\begin{proof}{} The unitarity of  $T_d$ follows from the definition
of the norm in $\Cal H^{\lambda }_\pm$.
Further, Proposition~\ref{PROP:222-333} implies that
$T_d\,T_b\,T_{d^{-1}} = T_{d\,b\,d^{-1}}$ for any 
$d\in D_1$ and $b \in {{P}}_0$. Hence the operators $T_d$, $d\in D_1$,
and $T_b$, $b \in {{P}}_0$, generate a representation of the group ${{P}}$.
\end{proof}

\begin{PROP}{}\label{PROP:211-86} 
The representations of the group ${{P}}$ in the spaces
$\Cal H^{\lambda }_+$ (respectively, 
$\Cal H^{\lambda }_-$) are irreducible and pairwise equivalent. The representations
in the spaces $\Cal H^{\lambda }_+$ and $\Cal H^{\mu}_-$
are not equivalent.
\end{PROP}

\begin{proof}{} The irreducibility immediately follows from the definition
of the operators of the subgroup $D_1$ and the fact that the irreducible
representations of ${{P}}_0$ in the spaces $\Cal H(\rho )$
are pairwise nonequivalent. An intertwining operator
$\tau:\,\Cal H^{\lambda }_{\pm}\to\Cal H^{\mu}_{\pm}$ has the form
$$
\tau:\,f(\rho ,z)= |\rho |^{\frac 12(\mu-\lambda )}\,f(\rho ,z).
$$
\end{proof}

\subsection{The action of the operators $T_g$, $g \in P$,
on elements of the form $f_k (\rho )\,z^k$ and $f_k (\rho )\,\OVER z^k$}
\label{sect:3-5(L)}

Denote by $(L_k^+)^{\lambda } \subset \Cal H_+^{\lambda }$
and $(L_k^-)^{\lambda } \subset \Cal H_-^{\lambda }$, $k \in \ZZ^{n-1}_+$,
the subspaces of functions of the form $f_k (\rho )\,z^k$ and 
$f_k (\rho )\,\OVER z^k$, respectively. It follows from the definition that
$$
\Cal H^{\lambda }_{\pm}
= \bigoplus _{k \in \ZZ^{n-1}_+} (L_k^\pm)^{\lambda }.
$$
Let us describe the action of the operators of the group
${{P}}$ on elements of 
$L_k^{\lambda } = (L_k^\pm)^{\lambda }$. For definiteness, we restrict
ourselves to the case of the subspaces $\Cal H^{\lambda }_{+}$.

\begin{PROP}{}\label{PROP:744-55}
The action of the operators  $T^+_g$ corresponding to elements
$g=(t_0,z_0)\in H$ and $d=\diag(\bar \epsilon ^{-1},u, \epsilon)  \in D$
on elements $f_k(\rho)z^k \in L^+_k$ is given by
\begin{equation}{}\label{738-24}
T^+_{(t_0,z_0)}( f_k(\rho )z^k)=\sum_m a_{km}(t_0,z_0;\rho )\,z^m,
\end{equation}
where
$$
a_{km}(t_0,z_0;\rho ) = k!\,e^{\,\rho (it_0-\frac12\,|z_0|^2)}\,f_k(\rho )\,
\sum_s \frac{
        (- \rho )^{ | m-s | } z_0^{k-s}\,\bar z_0^{m-s}
        }{
        s!(k-s)!(m-s)!},
$$
\begin{equation}{}\label{744-56}
T^+_d (f_k(\rho )\,z^k) = k!\,
|\epsilon |^{\lambda }\, \bar\epsilon^{-|k|}
f_k(|\epsilon |^2 \rho ) \sum_{|l|=|k|} L_{kl}(u)\,z^l.
\end{equation}
Here
\begin{equation}{}\label{744-57}
L_{kl}(u)=\sum_{M(k,l)}\,\frac{u_{ij}^{m_{ij}}}{m_{ij}!}
\end{equation}
and the sum in \eqref{744-57} ranges over the set $M(k,l)$ of integer 
$(n{-}1)$-matrices $m=\|m_{ij}\|$, $m_{ij}\ge0$, satisfying the conditions
$$
\sum_j m_{ij}=l_i, \quad \sum_i m_{ij}=k_j.
$$
The operators $T^-_g$ act on elements from $L^-_k$ in a similar way.
\end{PROP}

\begin{proof}{} It follows from the description of the operators $T_g$ 
in the space $\Cal H_+$ that
\begin{gather*}{}
T^+_{(t_0,z_0)}( f_k(\rho )\,z^k)
= e^{\,\rho (it_0-\frac12\,|z_0|^2)}\,e^{\,- \rho\,zz^*_0}\,f_k(\rho )\,(z+z_0)^k
\\
= e^{\,\rho (it_0-\frac12\,|z_0|^2)}\,
(\sum_l\frac{( - \rho )^{ |l| } }{  l!  } \bar z_0^l z^l )\,
(k!\,\sum_{s\le k}
           \frac{z^s \bar z_0^{k-s}
                }{
                s!\,(k-s)!
                }
                )\,f_k(\rho )
\\
= k!\, e^{\,\rho (it_0-\frac12\,|z_0|^2)}\,f_k(\rho )\,\sum_m
(\sum_s \,
        \frac{(- \rho )^{  |  m-s | }
                z_0^{k-s}\,\bar z_0^{m-s}
                }{
                s!(k-s)!(m-s)!
                }
                )\,z^m.
\end{gather*}
This implies \eqref{738-24}. Further, we have
$$
T^+_d(f_k(\rho )z^k)=|\epsilon |^{\lambda }\,f_k (|\epsilon |^2 \rho )
(\bar\epsilon^{-1} z u)^k.
$$
Obviously,
$$
(\bar\epsilon^{-1} z u)^k
= k!\,\bar\epsilon^{-|k|}\,\sum_{l|=|k|} L_{kl}(u)\,z^l,
$$
where $L_{kl}(u)$ is given by \eqref{744-57}.
This implies \eqref{744-56}.
\end{proof}

\begin{REM*}{} In~\S\,\ref{sect:5} we will show that the representations 
of the group ${{P}}$ in the spaces $\Cal H_\pm^0$ can be extended to unitary representations
of the whole group  $U(n,1)$. These representations of $U(n,1)$, 
called special representations, have nontrivial $1$-cocycles.
\end{REM*}

\subsection{Embeddings of the representations of the groups
${{P}}_0$ and ${{P}}$ in
$\Cal H(\rho )$ and $\Cal H^{\lambda }_\pm$ into their tensor products}

\begin{THM}{}\label{THM:333-340} For any
$\rho _1>0, \dots , \rho_m >0$ there exists a unique embedding of Hilbert spaces
\begin{equation}{}\label{811-21}
\tau : \Cal H(\rho_1 + \ldots + \rho _m ) \to \bigotimes _{k=1}^m \Cal H(\rho _k)
\end{equation}
commuting with the action of the subgroup ${{P}}_0$. Namely,
\begin{equation}{}\label{333-341}
\tau f(z) = F(z_1, \dots ,z_n)
= f (\rho ^{\,- 1}\, {\textstyle\sum\,} \rho _i\,z_i),
\quad\text{where}\quad  \rho = {\textstyle\sum\,} \rho _i.
\end{equation}
In a similar way we define the unique embedding \eqref{811-21}
commuting with the action of the subgroup 
${{P}}_0$ for any $\rho _1<0, \dots , \rho_m <0$.
\end{THM}

\begin{proof}{}
It follows from Proposition~\ref{PROP:713-27} (corollary) that the tensor
product $\bigotimes _{k=1}^m \Cal H(\rho _k)$ contains a unique vector
invariant under the subgroup $D_0$. Hence there exists at most
one embedding  \eqref{811-21} commuting with the action of
${{P}}_0$. Let us show that \eqref{333-341} is the required embedding.

By definition,
\begin{gather*}{}
\| \,\tau f \,\|^{2} = \prod \rho _i^{n-1}\,
\int\limits_{(\CC^{n-1})^m}
| f(\rho ^{\,- 1}\, {\textstyle\sum\,} \rho _i\,z_i)|^{2} \,
e^{\,-\sum_{i=1}^m \rho _i |z_i|^2} \,\prod d\mu(z_i)
\\
= \int\limits_{(\CC^{n-1})^m}
| f(\rho ^{\,- 1}\, {\textstyle\sum\,} \rho _i^{1/2}\,z_i)|^{2} \,
e^{\,-\sum_{i=1}^m  |z_i|^2} \,\prod d\mu(z_i).
\end{gather*}
Let us replace the variables $z_1, \dots ,z_m$,
$z_i\in\CC^{n-1}$, in this integral with the new variables
$$
z= \rho ^{-1/2}\,\sum_{i=1}^m \rho _i^{1/2}\,z_i \quad \text{and}\quad
w_k=\sum_{i=1}^m a_{ki} z_i, \quad k=1, \dots ,m{-}1,
$$
where the vectors $a_k=(a_{k1}, \dots ,a_{km})$, $k=1, \dots ,m{-}1$,
together with the vector
$a=( \rho ^{-1/2}\rho_1 ^{1/2}, \dots ,\rho ^{-1/2}\rho_m ^{1/2})$
form an orthonormalized basis in $\RR^m$. We obtain
\begin{align*}{}
\|\,\tau f \,\|^{2}
&
= \int\limits_{\CC^{n-1}}|f(\rho ^{-1/2} z)|^2\,
e^{- |z|^2}\,d\mu(z)\,\int\limits_{(\CC^{n-1})^{m-1}}
e^{\,-\sum_{i-1}^{m-1}|w_i|^2} \,\prod d\mu(w_i)
\\ &
=  \rho ^{n-1}\,
\int\limits_{\CC^{n-1}}|f( z)|^2\,
e^{- \rho |z|^2}\,d\mu(z)=\|f\|^2.
\end{align*}
It follows that $\tau$ is an embedding of Hilbert spaces. One can directly
check that $\tau $ commutes with the operators of the subgroup 
${{P}}_0$.
\end{proof}

\begin{REM*}{} It is not difficult to see that the representation of
the Heisenberg group $H$ in the tensor product
$\bigotimes _{i=1}^{m} {\mathcal H}(\rho _i)$
is a multiple of the irreducible representation with parameter
$\rho =\sum \rho _i$.
\end{REM*}

\begin{THM}{}\label{THM:765-10} 
For any $\lambda _1>0, \dots ,\lambda _m>0$ there exists an embedding of Hilbert spaces
$$
\tau : \Cal H_+^\lambda \to \bigotimes _{i=1}^m \Cal H_+^{\lambda _i},
\qquad \lambda = {\textstyle\sum\,}\lambda _i,
$$
commuting with the action of the subgroup ${{P}}$. Namely,
\begin{equation}{}\label{765-10}
\tau f(\rho ,z) =
\Bigl(\frac{\prod  \Gamma (\lambda _i) }{ \Gamma (\lambda ) }
\Bigr)^{-\frac12}
f({\textstyle\sum\,}\rho _i,\,\rho ^{-1}{\textstyle\sum\,}\rho _i\,z_i),
\end{equation}
where $\rho = \sum \rho _i$. Similarly, for any
$\lambda _1<0, \dots ,\lambda _m<0$ there exists an embedding of Hilbert spaces
$$
\tau : \Cal H_-^\lambda \to \bigotimes _{i=1}^m \Cal H_-^{\lambda _i},
\qquad \lambda = {\textstyle\sum\,}\lambda _i,
$$
commuting with the action of the subgroup ${{P}}$. 
\end{THM}

\begin{proof}{} By definition,
\begin{multline*}
\|\tau f \|^2 =
\Bigl(\frac{ \prod \Gamma (\lambda _i) }{ \Gamma (\lambda ) }
\Bigr)^{-1}\,
\times \\ \times
\int_{\RR^m_+}
\Bigl( \int _{(\CC^{n-1})^m} |
f({\textstyle\sum\,}\rho _i,\,\rho ^{-1}{\textstyle\sum\,}\rho _i\,z_i)
|^2\,e^{\,- \sum \rho _i \,|z_i|^2}\,\prod d \mu (z_i)\,
\Bigr)\,
\prod \rho _i^{n+\lambda _i - 2}\,d \rho _i.
\end{multline*}
In view of Theorem~\ref{THM:333-340}, it follows that
$$
\|\tau f \|^2 =
\Bigl(\frac{\prod  \Gamma (\lambda _i) }{ \Gamma (\lambda ) }
\Bigr)^{-1}\,
\int_{\RR^m_+}
\Bigl( \int _{\CC^{n-1}} | f(\rho ,z) |^2\,
e^{\,- \rho \,|z|^2}\,d \mu (z)\,
\Bigr)\,
\rho ^{n-1}\,
\prod \rho _i^{\lambda _i - 1}\,d \rho _i.
$$
The obtained expression can be transformed into
\begin{multline*}
\|\tau f \|^2 =
\\
=\Bigl(\frac{ \prod  \Gamma (\lambda _i) }{ \Gamma (\lambda ) }
\Bigr)^{-1}\,
\int_{0}^{\infty }
\Bigl( \int _{\CC^{n-1}} | f(\rho ,z) |^2\,
e^{\,- \rho \,|z|^2}\,d \mu (z)\,
\Bigr)\,
\rho^{ n + \lambda - 2 }\,d \rho
\int _{\sum \rho _i =1\,}
\prod \rho _i^{\lambda _i - 1}\,d \rho _i.
\end{multline*}
Note that the second integral equals
$\frac{ \prod \Gamma (\lambda _i) }{ \Gamma (\lambda ) }$.
Hence $\| \tau f\| = \|f\|$, i.e.,  $\tau $ is an embedding. The fact that
it commutes with the action of ${{P}}$ follows from the explicit
formulas for the corresponding operators.
\end{proof}

\section{Representations of the groups ${{P}}_0$ and ${{P}}$ in the spaces 
$H(\rho )$ and $H^{\lambda} _{\pm}$}\label{sect:4}

\subsection{Representations of the group ${{P}}_0$ in the space $H(\rho )$}
\label{sect:3-1}
Denote by $H(\rho )$, $\rho \ne 0$, the Hilbert space of
{\it all} functions $f(z)$ on $\CC^{n-1}$ with the norm \eqref{444-444}:
$$
\| f \|^{2} = |\rho |^{n-1}
\int_{\CC^{n-1}} |f(z)|^{2}\,e^{\, - |\rho| \,|z|^2}\,d\mu(z).
$$
Note that $H(\rho )$ is spanned by the monomials $z^k\OVER z^l$,
$k,l\in\ZZ^{n-1}_+$.

Let us define a unitary representation of the group
${{P}}_0$ by the same formulas as in the case of 
$\Cal H(\rho )$, i.e., by \eqref{746-24},
\eqref{285-14} for $\rho >0$, and by  \eqref{746-25}, \eqref{285-14} for
$\rho <0$.

In particular, $\Cal H(\rho )$ is an invariant subspace of  $H(\rho )$.

\begin{PROP}{}\label{PROP:27-676} Theorem~\ref{THM:333-340},
concerning embeddings of the subspaces
$\Cal H(\rho )$ into their tensor products, translates without any change
to the spaces $H(\rho )$.
\end{PROP}

\subsection{Decomposition of $H(\rho )$ into irreducible subspaces}
For definiteness, consider the case $\rho>0$. Similar results for
$\rho <0$ can be obtained by replacing functions with their complex conjugates.

Denote by $K_m$, $m=0,1, \ldots $, the subspace of 
$H(\rho )$ consisting of all homogeneous polynomials in
$\OVER z_1, \dots , \OVER z_{n-1}$ of homogeneous degree $m$.

Obviously, the spaces $K_m$ are invariant, irreducible, and pairwise 
nonequivalent with respect to the subgroup $U(n{-}1) \subset {{P}}_0$.
It is also clear that the spaces $K_m$ are mutually orthogonal. 
Note that the lowest weight vector in $K_m$ is the vector
$\bar z_{n-1}^m$ of weight $(0, \dots ,-m)$.

\begin{DEF*}{} Denote by $H_m(\rho)$,
$m=0,1, \ldots$, the invariant subspace in $H(\rho)$ that is cyclically
generated by the vectors $\psi\in K_m$.

In particular, $H_0(\rho)$ coincides with the space $\Cal H(\rho )$
of entire analytic functions in $z_i$ introduced above.
\end{DEF*}

\begin{REM*}{}
Let us give another equivalent definition of the spaces
$H_m(\rho  )$. It follows immediately from the description of the representation
that the operators
$\partial _i=\frac{\partial }{\partial \, \bar z _i}$,
$i=1, \dots ,n-1$, satisfy the following commutation relations:
\begin{equation}{}\label{0-1-3-55}
 \partial _i\,T_h=T_h\,\partial _i \quad \text{for all}
\quad h\in H,
\end{equation}
\begin{equation}{}\label{0-1-3-56-1}
\partial _i\,T_{\epsilon ,u }=  \epsilon^{-1}\,
(\sum_{j=1}^{n-1}\bar u_{ij}\,\partial _j)\,  T_{\epsilon ,u}
\quad \text{for all} \quad (\epsilon ,u)\in D_0.
\end{equation}

Denote by $L_m(\rho )$, $m=0 ,1, \ldots  $, the Hilbert subspace of functions
$\psi\in H(\rho )$ such that
$$
p_m(\partial _1, \dots ,\partial _{n-1} )\,\psi=0
$$
for all homogeneous polynomials $p_m$ of degree $m$. Obviously,
$L_m(\rho ) \subset L_{m+1}(\rho ) $. In this notation, we have
$$
H_m(\rho )=L_{m+1}(\rho )   \ominus  L_{m}( \rho ),\quad
m=0,1, \ldots .
$$
\end{REM*}

\begin{PROP}{}\label{PROP:3-1} The representations of the subgroup ${{P}}_0$
in the spaces $H_m(\rho )$ are irreducible and pairwise nonequivalent.
\end{PROP}

\begin{proof}{} Assume that $H_m(m)$ is reducible and hence there exists a 
decomposition $H_m(\rho ) = H_1 \oplus H_2 $ into invariant subspaces.
Since $H_m(\rho )$ contains a unique lowest weight vector, this vector
belongs to one of these subspaces, say, $H_1$. But then
$K_m \subset H_1$, whence $H_1=H_m(\rho )$.

Since the lowest weight vectors in the spaces $H_m(\rho )$ have different
weights, these spaces are pairwise nonequivalent.
\end{proof}

\begin{COR*}{}
The subspaces $H_m(\rho )$ are mutually orthogonal, and the space
$H(\rho )$ is their direct sum:
$$
H(\rho )= \bigoplus _{m=0}^\infty H_m (\rho ).
$$
\end{COR*}

\begin{PROP}{}\label{PROP:27-348} The restriction of the representation
of the group ${{P}}_0$ in the space $H_m(\rho )$ to the subgroup $H$ 
is a multiple of the irreducible representation of $H$ with parameter 
$\rho $, the multiplicity being equal to $\dim K_m$. The decomposition of $H_m(\rho )$ 
into $H$-irreducible subspaces can be written in the form
$$
H_m(\rho )= \bigoplus _{|k|=m} H_{mk}(\rho ),
$$
where $ H_{mk}(\rho )$ is the $H$-irreducible
subspace generated by the vector $\OVER z^k\in K_m$.
\end{PROP}

\subsection{The orthogonal basis in $H_m(\rho )$ associated with the action of 
the Heisenberg group $H$}
Set
\begin{equation}{}\label{55-22-1}
f_p(z) =\OVER z^p=\OVER z _1^{p_1} \dots \OVER z _{n-1}^{p_{n-1}},
\quad p\in\ZZ_+^{n-1}.
\end{equation}
It follows from the definition of the norm in $H(\rho )$ that
\begin{equation}{}\label{65-22-1}
\|f_p\|^2=p!\, \rho ^{-|p|}.
\end{equation}

Let $L_p(\rho ) \subset H(\rho )  $ be the $H$-invariant subspace
cyclically generated by the vector $f_p(\rho )$. By
Proposition~\ref{PROP:27-348}, the subspace $L_p(\rho )$ is irreducible, and
$$
H_m(\rho )= \bigoplus _{|p|=m}L_p(\rho ), \quad
H(\rho )= \bigoplus _{p\in\ZZ^{n-1}_+}L_p(\rho ).
$$

Let us define an orthogonal basis in each subspace $L_p(\rho )$. By definition,
the action of the operators $T_h$, $h=(t,z)\in H$, on the vectors
$f_p$ is given by
$$
T_hf_p(z)=(\OVER z +\OVER a)^p\,
e^{\,\rho (it-\frac{aa^*}{2}-  za^*)},
$$
where $\OVER z ^p=
\OVER z _1^{p_1} \ldots \OVER z_{n-1}^{p_{n-1}}.$

\begin{DEF*}{} Let us define functions $f_{p,q}( z)$, $p,q\in\ZZ_+^{n-1}$,
as the coefficients of the following expansion:
\begin{equation}{}\label{010-11}
(\OVER z +\OVER a)^p\,e^{\,-\rho\,z a^*}=\sum_{q\in\ZZ_+^{n-1}}
f_{p,q}(z )\,\OVER a^q, \quad
\OVER a^q= \OVER a_1^{q_1} \ldots \OVER a_{n-1}^{q_{n-1}},
\end{equation}
or, which is equivalent, as the coefficients of the expansion
\begin{equation}{}\label{910-92}
e^{b(z +a)^* - \rho\,z a^*}=\sum_{p,q\in\ZZ_+^{n-1}} (p!)^{-1}\,
f_{p,q}(z )\,b^p\OVER a^q.
\end{equation}
\end{DEF*}
In particular, $f_{p,0}=f_p$. It follows from this definition that
\begin{enumerate}
\item $f_{p,q}$ are polynomials in $z_i$ and $\OVER z _i$;

\item for every fixed $p$, the functions $f_{p,q}$ linearly span
$L_p(\rho)$.
\end{enumerate}

\begin{PROP}{}\label{PROP:3-3} The polynomials $f_{p,q}$ 
are given by the following explicit expression:
\begin{equation}{}\label{939-22}
f_{p,q}( z )=p!\,\sum_s (- \rho) ^{|q-s|}\,
\frac{z ^{q-s}\,\OVER z^{p-s}}{s!\,(p-s)!\,(q-s)!},
\end{equation}
where the sum ranges over $s \in \ZZ_+^{n-1}$ such that $s\le p$, $s\le q$.
\end{PROP}

Indeed, we have
\begin{gather*}{}
e^{b(z +a)^* - \rho\,z a^*} =e^{bz^*}\,
e^{(b- \rho\,z )a^*}
=\sum_r\frac{b^r\OVER z ^r}{r!}\,
\sum_q\frac{(b- \rho\,z) ^q\OVER a^q}{q!}
\\
=\sum_{r,s,q}(- \rho) ^{|q-s|}\,\frac{\OVER z ^r z ^{q-s}}{r!s!(q-s)!}
\,b^{r+s}\,\OVER a^q
=\sum_{p,q}\Bigl(\sum_s\frac{(- \rho) ^{|q-s|}\,\OVER z ^{p-s}\,
z ^{q-s}}{s!(p-s)!(q-s)!}\Bigr)\,b^p\,\OVER a^q.
\end{gather*}
This implies \eqref{939-22}.

\begin{THM}{}\label{THM:3-3} The polynomials $f_{p,q}$, $|p|=m$,
form an orthogonal basis in $H_m(\rho )$, $m=0,1,\ldots$, and
\begin{equation}{}\label{010-12}
\|f_{p,q}\|^2=\frac{\rho ^{|q|}}{q!}\,\|f_p\|^2
\quad \text{for any}
\quad p,q.
\end{equation}
\end{THM}
\begin{proof}{} The orthogonality of $f_{p,q}$ and $f_{p',q'}$ for $p\ne p'$
follows from the mutual orthogonality of the spaces $L_p(\rho )$.
Further, substituting the expression for  $T_h f_p$ into the equation
$ \langle T_hf_p,T_hf_p  \rangle \,
= \langle f_p,f_p \rangle$, we obtain
$$
\sum_{q,q'} \langle f_{p,q},f_{p,q'} \rangle \,\OVER a^q\, a^{q'}
= \langle f_p,f_p \rangle \,e^{\,\rho a a^*}
= \langle f_p,f_p \rangle\,\sum_q \frac{\rho ^{|q|}}{q!}\,\OVER a^q\,a^q.
$$
This implies \eqref{010-12}.
\end{proof}

\begin{COR*}{} For any $p,q\in\ZZ^{n-1}_+$,
\begin{equation}{}\label{713-35}
\|f_{pq}\|^2=\frac{p!}{q!}\, \rho ^{|q-p|}.
\end{equation}
\end{COR*}

\subsection{Formulas for the operators of the group 
${{P}}_0$ in terms of the basis $\{f_{p,q}\}$ on $H(\rho )$}

\begin{PROP}{}\label{PROP:3-4} The operators $T_{t,a}$, $(t,a)\in H$, have the form
\begin{equation}{}\label{777-11}
T_{t,a}\,f_{p,q}=e^{\,\rho (it-\frac{aa^*}{2})}\sum_{q'} c_{qq'}(\rho , a)\,
f_{p,q'},
\end{equation}
where
\begin{equation}{}\label{777-12}
c_{qq'}(\rho , a)=\sum_s \frac{(- \rho )^{|s|}}{s!\,(q-s)!\,(q'-q+s)!}\,
a^s \OVER a^{q'-q+s}
\end{equation}
and the sum ranges over all $s\in\ZZ^{n-1}_+$ such that
$s\le q$ and $q'-q+s\ge 0$.
\end{PROP}

\begin{proof}{} For every $b\in\CC^{n-1}$, we have
\begin{gather*}{}
\sum_q (T_{t,a}f_{p,q}( z ))\,\OVER b^q=
T_{t,a}[(\OVER z +\OVER b)^p\,e^{\,- \rho\,z b^*}]
=e^{\,\rho (it-\frac{aa^*}{2}- z a^* )}\,
(\OVER z +\OVER b+\OVER a)^p\,e^{\,- \rho (z+a) b^*}
\\
=e^{\,\rho (it-\frac{aa^*}{2})}\,
(\OVER z +\OVER b+\OVER a)^p\,e^{\,-\rho\,z(b+ a)^* }
\,e^{\,- \rho a b^*} .
\end{gather*}
Since
\begin{gather*}{}
(\OVER z +\OVER b+\OVER a)^p\,e^{\,-\rho\,z(b+ a)^* }=
\sum_{q'} f_{p,q'}\,(\OVER b+\OVER a)^{q'}
=\sum_{r, q'}\frac{q'!}{r!\,(q'-r)!} f_{p,q'}\,\OVER b^r\OVER a^{q'-r}
\end{gather*}
and
$
\displaystyle
e^{\,- \rho a b^*}=\sum_s \frac{(- \rho )^{|s|}}{s!}a^s\OVER b^s,
$
we obtain
\begin{gather*}{}
\sum_q (T_{t,a}f_{p,q}(\rho , z ))\,\OVER b^q=
e^{\,\rho (it-\frac{aa^*}{2})}\,
\sum_{r,s,q'}\frac{(- \rho )^{|s|}}{r!\,(q'-r)!\,s!} a^s\,\OVER a^{q'-r}\,
\OVER b^{r+s}.
\end{gather*}
Proposition~\ref{PROP:3-4} follows.
\end{proof}

To describe the operators $T_d$, $d=\diag( \epsilon ,u, \epsilon )$, where
$|\epsilon |=1,$ $u\in U(n{-}1)$, we use the polynomials introduced above 
in~\S\,\ref{sect:3-5(L)}:
$$
L_{pq}(x)=\sum_{M(p,q)}\bigl(\prod\frac{x_{ij}^{l_{ij}}}{l_{ij}!}\bigr),
\quad p,q\in\ZZ^{n-1}_+,
$$
where the sum ranges over the set $M(p,q)$ of integer
$(n{-}1)$-matrices $m=\|m_{ij}\|$, $ m_{i,j}\ge0$, satisfying the condition
\begin{equation}{}\label{178-28}
\sum_i l_{ij}=p_j \quad \text{for all} \quad j; \quad
\sum_j l_{ij}=q_i \quad \text{for all} \quad i
\end{equation}
(Louck polynomials).
Note that these polynomials are the coefficients of the expansion
\begin{equation}{}\label{735-45}
e^{axb}=\sum_{p,q} L_{pq}(x)\,a^q\,b^p, \quad a,b\in\CC^{n-1},
\end{equation}
where $axb=\sum_{i,j=1}^{n-1}a_ix_{ij}b_j.$

Obviously, $L_{pq}(x)=L_{qp}(x')$, where the prime denotes transpose.

Note that $L_{pq}(x)$ can be written as a sum over the subgroup
$M_0$ of all integer matrices $l$ for which all rows and all columns sum to $0$:
$$
L_{pq}(x)=\sum_{l\in M_0}
\bigl(\prod\frac{x_{ij}^{l_{ij}+l^0_{ij}}}{\Gamma(l_{ij}+l^0_{ij}+1) }\bigr),
$$
where $l^0$ is an arbitrary fixed integer matrix satisfying \eqref{178-28}.

\begin{PROP}{}\label{PROP:3-7} The action of the operators $T_d$, 
$d=\diag( \epsilon ,u, \epsilon )$, where 
$|\epsilon |=1$, $u\in U(n{-}1)$, on the basis vectors is given by the formula
\begin{equation}{}\label{721-13}
T_{\epsilon ,u}\,f_{p,q}= p\,!\epsilon^{|q-p|}\,
\sum_{p',q'}q'\,!\, L_{pp'}(\OVER u)\,L_{q q'}(u)\,f_{p', q'}.
\end{equation}
\end{PROP}

\begin{proof}{} First of all, we have
\begin{gather*}{}
\sum_q (T_{\epsilon ,e}\,f_{p,q})\,\OVER a^q =T_{\epsilon ,e}\,
(\OVER z +\OVER a)^p\,e^{\,- \rho\,z a^*}
=(\OVER {\epsilon z} +\OVER a)^p\,e^{\,- \rho \epsilon  z a^*}
\\
=\OVER \epsilon ^{|p|}\,(\OVER { z} +\epsilon\OVER a)^p\,
e^{\,- \rho \epsilon  z a^*}
=\OVER \epsilon ^{|p|}\,\sum_q f_{p,q}\,(\epsilon \OVER a)^q=
\epsilon ^{|q-p|}\sum_q f_{p,q}\,\OVER a^q.
\end{gather*}
Hence
\begin{equation}{}\label{378-12}
T_{\epsilon ,e}\,f_{p,q}= \epsilon ^{|q-p|}\,f_{p,q}.
\end{equation}
Further, for every $u\in U(n{-}1)$ we have
\begin{gather*}{}
\sum_q (T_{1,u}f_{p,q})\,\OVER a^q =(\OVER{z  u}+ \OVER a)^p\,
e^{\,- \rho\,z ua^*}
=((z +\OVER{au^*})\OVER u)^p\,e^{\,- \rho\,z (au^*)^*}
\\
=p\,!\sum_{p'} L_{pp'}(\OVER u)\,(z +\OVER{au^*})^{p'}\,
e^{\,- \rho\,z (au^*)^*}
= p\,!\sum_{p',q'} L_{pp'}(\OVER u)\,f_{p', q'}\,(\OVER{au^*})^{q'}
\\
= p\,!\sum_{p',q',q}q'\,!\, L_{pp'}(\OVER u)\,f_{p', q'}\,L_{q' q}(u')\,\OVER a^q.
\end{gather*}
Therefore
\begin{equation}{}\label{378-13}
T_{1,u}\,f_{p,q}=p\,!\sum_{p',q'}q'\,!\,
 L_{pp'}(\OVER u)\,L_{q' q}(u')\,f_{p', q'}.
\end{equation}
Combining \eqref{378-12} and \eqref{378-13} and taking into account that
$L_{q' q}(u')=L_{q q'}(u)$, we obtain  \eqref{721-13}.
\end{proof}

\begin{REM*}{}
Observe the following property of the polynomials $L_{pq}$.
The matrices
$$
D^m(x)=\|D_{pq}(x)\|_{|p|=|q|=m},
\quad\text{where}\quad D_{pq}(x)=(p!q!)^{1/2}\,L_{pq}(x),
$$
satisfy the relations
\begin{equation}{}\label{978-11}
D^m(xy)=D^m(x)\,D^m(y).
\end{equation}
Relation \eqref{978-11} remains valid if we replace
$(p!q!)^{1/2}$ in the definition of  $D_{pq}(x)$ by
$\prod_{i=1}^{n-1} p_i^{\alpha _i} q_i^{\beta _i}$, where $\alpha _i+ \beta _i=1$,
$i=1, \dots ,n{-}1$, in particular, if we replace $(p!q!)^{1/2}$ by $p!$ or $q!$.
\end{REM*}

\subsection{Representations of the group ${{P}}$ in the spaces 
$H^{\lambda }_{\pm}$}
By analogy with the spaces $\Cal H^{\lambda }_{\pm}$, we define spaces
$H^{\lambda }_{\pm}$ as the following direct integrals:
$$
H^{\lambda }_-=\int_{- \infty} ^0  H(\rho )\,| \rho |^{\lambda -1}\,
d \rho ; \quad
H^{\lambda }_+=\int_{-\infty} ^0  H(\rho )\, \rho ^{\lambda -1}\,
d \rho .
$$
In more detail, $H^{\lambda }_+$ is the Hilbert space of ALL functions
$f(\rho ,z)$ on $\RR_+ \times \CC^{n-1}$ with the norm 
\begin{equation}{}\label{222-13}
\|f\|^2=\int_0^\infty \Bigl(
\int_{\CC^{n-1}} |f(\rho ,z)|^{2}\,e^{\, - |\rho| \,|z|^2}\,d\mu(z)\Bigr)\,
\rho ^{n+\lambda -2}\,d \rho.
\end{equation}
The space $H^{\lambda }_-$ is defined in a similar way.

The unitary representations of the group ${{P}}_0$ in the spaces $\Cal H(\rho )$
induce unitary representations of this group in the spaces $H^{\lambda }_{\pm}$.

\begin{PROP}{}\label{PROP:211-055} The representations of the group ${{P}}$ in the spaces 
$\Cal H^{\lambda }_{\pm}$ can be extended to unitary representations
of the group ${{P}}={{P}}_0\leftthreetimes D_1$. Namely, the operators
$T_d$ corresponding to elements $d\in D$ are given by the formula \eqref{273-66}:
$$
T_d f(\rho ,z)=f(r^2 \rho ,r^{-1}z)\,r^{ \lambda }\quad \text{for}\quad
d=\diag(r^{-2},e,r),\quad r>0.
$$
\end{PROP}

\begin{PROP}{}\label{PROP:76-73-2} The spaces $H^{\lambda }_{\pm}$
decompose into a direct sum of pairwise nonequivalent $P$-irreducible subspaces:
$$
H^{\lambda }_{\pm}= \bigoplus _{m=0}^\infty (H^{\lambda }_{\pm})_m,
$$
where
$$
(H^{\lambda }_-)_m=\int_{- \infty }^0 H(\rho ) |\rho |^{\lambda -1}\,d \rho ,
\quad
(H^{\lambda }_+)_m=\int^{ \infty }_0 H(\rho ) \rho ^{\lambda -1}\,d \rho .
$$
\end{PROP}

\begin{COR*}{} The representations of the group ${{P}}$ in the spaces
$ ( H^{\lambda }_+ )_m  $ (respectively,
$ ( H^{\lambda} _- )_m  $) are pairwise nonequivalent.
\end{COR*}

\begin{PROP}{}\label{PROP:227-696} Theorem~\ref{THM:765-10},
concerning embeddings of the spaces $\Cal H^{\lambda  }_{\pm}$ 
into their tensor product, translates without any change to the spaces
$H^{\lambda  }_{\pm}$.
\end{PROP}

\subsection{Another realization of representations of the groups ${{P}}_0$ and ${{P}}$}

\begin{PROP}{}\label{PROP:713-26} Every Hermitian functional in the space of 
functions $f(z)$ on $\CC^{n-1}$ that is invariant with respect to the operators
\eqref{746-24}, \eqref{285-14} (for $\rho >0$)
of the group
${{P}}_0$ can be written in the form
\begin{equation}{}\label{412-32}
\Phi(f)=\int\limits_{\CC^{n-1} \times \CC^{n-1}} f(z)\,\OVER{f(z')}\,
e^{\rho (z'z^*-|z|^2-|z'|^2)}\,Q(\rho |z-z'|)\,d\mu(z)\,d\mu(z'),
\end{equation}
where the kernel $Q$ is an arbitrary (generalized) function.
\end{PROP}

In particular, if $Q$ is the delta function, then this functional coincides,
up to a factor, with the squared norm in the space $H_+(\rho )$.

Assume that $\Phi$ is a positive definite functional and denote by
$H(Q, \rho )$ the Hilbert space of functions $f(z)$ on $\CC^{n-1}$
with the norm $\|f\|^2=\Phi(f)$. Formulas \eqref{746-24}, \eqref{285-14}
define a unitary representation of the group ${{P}}_0$ in this space.

In a similar way we define a representation of $P_0$ in $H(Q, \rho )$ for $\rho <0$.

\begin{PROP}{}\label{PROP:970-31}  If $H(Q, \rho )$ contains all monomials
$z^k\OVER z^l$, $k,l\in\ZZ^{n-1}_+$, and is generated by them, then the
representation of the group ${{P}}_0$ in $H(Q, \rho )$ is equivalent to its
representation in $H(\rho )$.
\end{PROP}

Further, let us introduce the following direct integrals of the Hilbert spaces $H(Q, \rho )$:
$$
H^{\lambda }_+(Q)=\int_0^\infty H(Q, \rho )\, \rho ^{2n+ \lambda -3}\,d \rho.
$$
The representations of ${{P}}_0$ in the spaces $H^{\lambda }_+(Q)$
can be extended to unitary representations of the group
${{P}}$ by the same formulas as for the space $H^{\lambda }_+$.

\begin{PROP}{}\label{PROP:970-331}  For all $Q$ satisfying the
positivity condition and the assumptions of Proposition~\ref{PROP:970-31}, 
and for all $\lambda $, the representation of the group ${{P}}$ in the space
$H^{\lambda }_+$ is equivalent to its representation in the space
$\Cal H^{\lambda }_+$.
\end{PROP}

\part{Representations of the group $U(n,1)$}

\section{Models of the complementary series of irreducible unitary
representations of the group $U(n,1)$}\label{sect:5}

Every complementary series representation is determined by a real number
$\lambda  $ from the interval $0< \lambda <2n$. Let us describe
several models of these representations.

\subsection{Model $A$: realization in a space of functions on the 
unitary unit sphere  $S=S^{2n-1}$ in $\CC^{n}$}
The unitary unit sphere in $\CC^n$,
$$
S: \quad  |\omega  |\equiv |\omega _1|^2+ \ldots + |\omega _n|^2=1,
$$
is a homogeneous space of the group $U(n,1)$. Namely, if this group
is realized in matrix model  $a$ (see~\S\,\ref{sect:1}),
then its action on  $S$ is given by the formula
$$
\omega g=(\omega \beta + \delta )^{-1}\,(\omega \alpha + \gamma ) \quad
\text{for}\quad
g=\begin{pmatrix}\alpha  &\beta  \\ \gamma   &\delta \end{pmatrix}.
$$

The complementary series representation with parameter $\lambda $ 
acts in the Hilbert space $L^{\lambda }$ of functions $f(\omega )$ on 
$S$ with the norm
\begin{equation}{}\label{22-1}
\|f\|^2=\int_{S \times S} |1-  \langle \omega , \omega' \rangle\,|^{- \lambda }\,
f(\omega )\,\OVER{f(\omega ')}\,d \omega \,d \omega ',
\end{equation}
where $\langle \omega , \omega ' \rangle
= \omega _1\OVER{\omega '_1}+ \ldots + \omega _n\OVER{\omega '_n}$
and $d \omega $ is the invariant measure on $S$.

\begin{REM*}{} In the interval $n< \lambda <2n$ the integral should be
understood in the sense  of the regularized value,
see~\cite{V-G-6}.
\end{REM*}

The operators $T^{\lambda }_g$,
$g=\begin{pmatrix}\alpha  &\beta  \\ \gamma   &\delta \end{pmatrix}$,
of this representation are given by the formula
\begin{equation}{}\label{22-2-222}
T^{\lambda }_g f(\omega )=f(\omega g)\,|\omega \beta + \delta |^{\lambda - 2n}.
\end{equation}

The group property of these operators is a consequence of the following
relation for the function
$b(\omega ,g)=|\omega \beta + \delta |$:
$$
b(\omega ,g_1g_2)=b(\omega ,g_1)\,b(\omega g_1 ,g_2) \quad \text{for any}
\quad g_1,g_2\in U(n,1)
$$
(the $1$-cocycle property). The unitarity follows from the relations
$$
|1-  \langle \omega , \omega ' \rangle\,|=
|1-  \langle \omega g, \omega 'g \rangle\,|\,
b(\omega ,g)\,b(\omega' ,g)
$$
and
$$
\frac{d(\omega g)}{d \omega }=b^{-2n}(\omega ,g)
\quad \text{for every $g \in U(n,1)$}.
$$

\begin{PROP}{}\label{PROP:2957-2}
Complementary series representations  $T^{\lambda}$ and $T^{\mu}$
of the group $U(n,1)$ are equivalent if and only if
$\lambda +\mu=2n$ or $\lambda = \mu $.

An intertwining operator
$R^{\lambda }: L^{\lambda } \to L^{2n- \lambda }$ has the form
$$
(R^{\lambda }f) (\omega ) = \int _{S}
|1- \langle \omega ,\omega \rangle\, '|^{- \lambda }
\,f(\omega ) \,d \omega '.
$$
\end{PROP}

\subsection{Model $B$: realization in a space of functions on the 
Heisenberg group $H$}
In model $B$ we use matrix model $b$ of the group $U(n,1)$.

The complementary series representation with parameter $\lambda $  is realized in 
the Hilbert space $L^{\lambda  }$ of functions $f(h)=f(t,z)$
on the Heisenberg group $H$ with the norm
\begin{equation}{}\label{2-1}
\|f\|^2 =\int_{H \times H}R^{- \lambda }(h,h')\,f(h)\,\OVER{f(h')}\,dh\,dh',
\end{equation}
where $R(h_1,h_2)=|(h_1h_2^{-1})_{31}|$ and  $dh=dt\,d\mu(z)$
is the invariant measure on $H$.

In the coordinate form:
\begin{equation}{}\label{2-2}
R(h,h')=|\zeta +\bar \zeta'+zz^{\prime*}|=
|i(t-t'+ \opn{Im}(zz^{\prime*})-\frac{1}{2}|z-z'|^2|,
\end{equation}
where $\zeta =it-\frac{zz^*}{2}$.

To describe the operators of this representation, let us define an action 
of the group $U(n,1)$ on $H$ by the formula
$h\bar g=h'$, where $h'\in H$ is determined by the relation
$hg=b^*h'$, $b^*\in s{{P}}s$ (the subgroup of upper block triangular matrices).
In the coordinate form:
$$
(\zeta , z)\bar g= (\zeta' , z' ),
$$
where
\begin{equation}{}\label{01-2}
\zeta '=(\zeta g_{13}+ z g_{23}+g_{33})^{-1}\,
(\zeta g_{11}+ z g_{21}+g_{31}),
\end{equation}
\begin{equation}{}\label{001-2}
z'=(\zeta g_{13}+ z g_{23}+g_{33})^{-1}\,
(\zeta g_{12}+ z  g_{22}+g_{32}).
\end{equation}
Further, set
\begin{equation}{}\label{3-2-0}
\beta (h, g)=\beta (\zeta, z ,g )=|\zeta g_{13}+ z g_{23}+g_{33}|
\quad \text{for any}\quad h\in H,\quad g\in U(n,1).
\end{equation}

The operators $T^{\lambda }_g$ of the group $U(n,1)$ 
in the space $L^{\lambda }$ are given by the formula
\begin{equation}{}\label{2-3}
T^{\lambda }_g f(h)=f(h\bar g)\, \beta ^{\lambda -2n}(h,g).
\end{equation}

In particular,
\begin{gather*}{}
T^{\lambda }_{t_0,z_0}f(\zeta ,z)=f(\zeta + \zeta _0- zz_0^*, z+z_0), \quad
(\zeta _0,z_0)\in H,
\\
T^{\lambda }_{\epsilon ,u}f(\zeta ,z)=f(|\epsilon |^{-2} \zeta ,
\bar \epsilon ^{-1}zu)\,|\epsilon |^{\lambda -2n},\quad (\epsilon ,u)\in D,
\\
T^{\lambda}_s f(\zeta ,z)=f(1/ \zeta , z/ \zeta )\,|\zeta |^{\lambda -2n}.
\end{gather*}
The group property and the unitarity follow from the relations
$$
R(h,h')=R(h\bar g, h'\bar g)\,\beta (h,g)\,\beta (h',g)
$$
and
$$
\frac{d(h\bar g)}{dh} =\beta^{-2n} (h,g).
$$

\begin{PROP}{}\label{PROP:111-1}
An intertwining operator $A^{\lambda }:\, L^{\lambda }\to L^{2n- \lambda }$
for equivalent representations has the form
\begin{equation}{}\label{111-1}
A^{\lambda }\,f(h)=\int_H R^{- \lambda }(h,h')\,f(h')\,dh.
\end{equation}
Thus the following relations hold:
\begin{equation}{}\label{111-2}
T^{2n - \lambda  }_g\, A^{\lambda }=A^{\lambda }\,T^{\lambda  }_g \quad
\text{for every}\qquad g\in U(n,1),
\end{equation}
\begin{equation}{}\label{111-3}
\|f\|_{\lambda }= c(\lambda )\,\|A^{\lambda }f\|_{2n- \lambda },
\end{equation}
that is,
\begin{equation}{}\label{111-4}
R^{- \lambda }(h_1,h_2)= c(\lambda)\,
\int_{Z \times Z} R^{- \lambda }(h_1,h_1')\,
R^{- \lambda }(h_2,h_2')\,R^{ \lambda-2n }(h_1',h_2')\,dh_1'\,dh_2'.
\end{equation}
(The integral is understood in the sense of the regularized value.)
\end{PROP}

\subsection{Relation between models $A$ and $B$}

In matrix model $b$, the stationary subgroup of the point
$\omega _0=(1,0, \dots ,0)\in S$ coincides with the subgroup
${{P}}^*=s{{P}}s$ of upper block triangular matrices, i.e.,
$S= {{P}}^* \backslash  U(n,1)$. It follows that the group $H$ is a section of 
the fiber bundle 
$\tau:U(n,1)\to S$ over the punctured sphere, 
$$
\tau(h)= \omega _0 h, \quad  h \in H.
$$
Set
$$
\phi(h) = f(\tau h).
$$

\begin{PROP}{}\label{PROP:1}
The norm in the space of functions on the sphere $S$, written 
in the coordinates
$h=(\zeta ,z)$,
has, up to a constant factor, the following form:
\begin{equation}{}\label{2-4}
\|\phi\|^2=\int_{H \times  H}|\zeta +\bar \zeta '+ zz^{\prime*}|^{- \lambda }\,
|1- \zeta |^{\lambda -2n}\,|1- \bar\zeta' |^{\lambda -2n}\,\phi(h)
\OVER{\phi(h')}\,dh\,dh'.
\end{equation}
\end{PROP}

\begin{proof}{} The relation $\omega = \omega _0 h$ implies that for
$h=(\zeta ,z)$,
$$
\omega _1=\frac{1+\zeta }{1- \zeta },\quad \omega '=\frac{\sqrt 2 z}{1-\zeta},
$$
where $\omega '= (\omega _2, \dots , \omega _{n})$. Therefore
$$
|1-  \langle \omega , \omega ' \rangle\,|=2\,|\zeta +\bar \zeta '+zz^{\prime*}|
\,|1- \zeta |\,|1- \zeta '|.
$$
Further, it is easy to check that the measure $d \omega $ on $S$ 
and the measure $dh$ on $H$ are related by the formula
$$
d \omega =|1- \zeta |^{-2n}\,dh.
$$
Substituting these expressions into  \eqref{22-1}, we obtain \eqref{2-4}.
\end{proof}

Comparing \eqref{2-4} with the expression \eqref{2-1} for the norm, we
deduce the following corollary.

\begin{COR*}{} The function $F(\omega )$, $\omega \in S$, in model $A$
and the function $f(h)$, $h\in R$, in model $B$ are related by the formula
$$
f(h)=|1- \zeta |^{\lambda -2n}\,F(\tau h).
$$
\end{COR*}

\subsection{Model $C$: realization in a space of functions
$\psi(\rho ,z)$ on $\RR \times \CC^{n-1}$}
Model $C$ of the complementary series representations of the group
$U(n,1)$ is obtained from the previous one by passing from functions
$f(h)=f(t,z)$ in model $B$ to functions
$\psi(\rho ,z)$ on $\RR \times \CC^{n-1}$, where
$$
\psi(\rho, z ) = |\rho |^{-n+1}\,e^{\frac{|\rho |}{2}|z|^2}\,
\int_{- \infty }^{+ \infty }e^{\,-i\rho t}\,
f(t, z )\,d t.
$$
According to the inversion formula for the Fourier transform, we have
$$
f(t, z )=|\rho |^{n-1}\int_{- \infty }^{+ \infty } \psi(\rho, z )\,|\rho |^{n-1}\,
e^{\,i \rho t- \frac{|\rho |}{2}\, |z|^2 }\,d \rho.
$$

To describe the space $L^{\lambda }$ and the operators of the representation in this 
model, set
\begin{equation}{}\label{33-8}
Q_{\lambda }(t)=t^{\frac{1- \lambda }{2}}\,e^{\frac{t}{2}}\\
K_{\frac{1- \lambda }{2}}(\frac{|t|}{2}),
\quad \text{where $K_{\nu}$ is a Bessel function},
\end{equation}
\begin{equation}{}\label{33-81}
a_{\rho }(z_1,z_2)=
\begin{cases}
 z_1z_2^* & \text{for $\rho >0$,}\\[1ex]
 z_2z_1^* & \text{for $\rho <0$},
\end{cases}
\end{equation}
\begin{equation}{}\label{33-82}
R_{\lambda }(\rho ,z_1,z_2)=e^{|\rho |(a_{\rho }(z_1,z_2)-|z_1|^2-|z_2|^2)}\,
Q_{\lambda }(|\rho |\,|z_1-z_2|^2).
\end{equation}

\begin{PROP}{}\label{PROP:2}
In model $C$, the complementary series representation $T^{\lambda }$ 
of the group $U(n,1)$ is realized in the Hilbert space
$L^{\lambda }$ of functions
$\psi(\rho, z )$, $\rho\in\RR$, $z \in \CC^{n-1}$, with the norm
\begin{equation}{}\label{493-78}
\|\psi\|^2=\frac{2\pi^{1/2}}{\Gamma (\lambda /2)}\,
\int_{- \infty }^{+\infty } \|\psi\|^2_{\rho }\,
|\rho |^{ \lambda -1}\,d \rho ,
\end{equation}
where
\begin{equation}{}\label{33-33}
\|\psi\|^2_{\rho }=|\rho |^{2n-2}\int_{\CC^{n-1} \times \CC^{n-1}}
R_{\lambda }(\rho ,z_1,z_2)\,\psi(\rho ,z)\,
\OVER{\psi(\rho ,z)}\,d\mu(z_1)\,d\mu(z_2).
\end{equation}

The operators of this representation are given by the formulas
\begin{align}{}\label{78-24-1}
T^{\lambda }_{t_0, z _0}\psi(\rho, z ) &=
e^{\,i\rho\ t_0-|\rho |( \frac{|z_0|^2}{2}+a_{\rho }(z_0,z))}\,
\psi(\rho ,z+z_0 ) \quad  \text{for}\quad (t_0,z_0)\in H,
\\
\label{78-25-1}
T^{\lambda }_d \psi(\rho, z )
&
= \begin{cases}
        \psi(\rho,\,\epsilon  z u)&\text{for}\quad
        d=\diag(\epsilon ,u, \epsilon  )\in D_0\,,
\\[1ex]
        r^{\lambda }\, \psi(r^2\rho,\,r^{-1} z)&\text{for}\quad
        d=\diag(r^{-1},e,r)\in D_1,   \end{cases}
\end{align}
\begin{equation}{}\label{378-25}
T^{\lambda }_s \psi(\rho, z )=|\rho |^{-n+1}\,e^{\frac{|\rho | }{2}|z|^2}\,
\int_{- \infty }^{+ \infty }\int_{- \infty }^{+ \infty }
e^{-it(\rho + \rho '|\zeta |^{-2}}\,|\rho '|^{n-1}\,|\zeta |^{\lambda-2n }\,
e^{-\frac{|\rho ||z|^2 }{2|\zeta |^2 }}\,
\psi(\rho ',\frac{z}{\zeta })\,dt\,d \rho ',
\end{equation}
where $\zeta =it-\frac{|z|^2}{2}$.
\end{PROP}

\begin{proof}{}
The formulas for the operators of $T^\lambda$ immediately follow from
the corresponding formulas in the previous model.

Let us derive formula \eqref{493-78} for the norm in model $C$.

It follows from the definition of the norm in model $B$ that
$$
\|\psi\|^2 =
\int F(\rho _1,z_1; \rho _2,z_2)\,\psi(\rho_1,z_1)\,
\OVER{\psi(\rho_2,z_2)}\,d\mu(z_1)\,d\mu(z_2)\,
|\rho_1\rho_2|^{n-1}\,d \rho _1\,d \rho _2,
$$
where
\begin{gather*}{}
F(\rho _1,z_1; \rho _2,z_2)=
e^{\,-\frac12\,(|\rho _1||z_1|^2+|\rho _2||z_2|^2)} \times
\\
\times  \int_{- \infty }^{+ \infty }\int_{- \infty }^{+ \infty }
e^{i(\rho _1t_1- \rho _2t_2)}\,
|(t_1-t_2+ \opn{Im}(z_1z_2^*))^2+\frac14\,|z_1-z_2|^4|^{- \lambda /2}\,dt_1\,dt_2
\\
=\delta (\rho _1- \rho _2)\,2\,
e^{\,-i\rho_1\opn{Im}(z_1z_2^*)-\frac{|\rho _1|}{2}\,(|z_1|^2+|z_2|^2)}\,
\int_0^{+ \infty }
|t^2+\frac14\,|z_1-z_2|^4|^{- \lambda /2}\,\cos {\rho_1 t}\,dt.
\end{gather*}
According to~\cite[Vol.~1, p.~11, formula~(7)]{Bateman-Tables},
$$
\int_0^{+ \infty } (t^2+a^2)^{-\nu-\frac{1}{2}}\,\cos {\rho t}\,dt=
\frac{\pi^{1/2}}{\Gamma (\nu+\frac{1}{2})}\,(\frac{|\rho| }{2a})^{\nu}\,
K_{\nu}(a |\rho|).
$$
Therefore
\begin{gather*}{}
F(\rho _1,z_1; \rho _2,z_2)
\\
=\delta (\rho _1- \rho _2)\,\frac{2\pi^{1/2}}{\Gamma (\lambda /2)}\,
e^{\,-i\rho_1\,\opn{Im}(z_1z_2^*)-\frac{|\rho _1|}{2}\,(|z_1|^2+|z_2|^2)}\,
(\frac{|\rho| }{|z_1-z_2|^2})^{\frac{\lambda-1}{2}}\,
K_{\frac{\lambda -1}{2}}(\frac{|\rho| }{2}\,|z_1-z_2|^2)
\\
=\delta (\rho _1- \rho _2)\,\frac{2\pi^{1/2}}{\Gamma (\lambda /2)}\,
e^{s(\rho ,z_1,z_2)}\,
Q_{\lambda } (|\rho|\,| z _1- z_2|^2)\,|\rho |^{\lambda -1},
\end{gather*}
where
$$
s(\rho ,z_1,z_2)=
-i\rho_1\opn{Im}(z_1z_2^*)-\frac{|\rho _1|}{2}(|z_1|^2+|z_2|^2)
-\frac{|\rho |}{2}\, |z_1-z_2|^2.
$$
To prove \eqref{493-78}, it suffices to observe that
\begin{equation}{}\label{333-8}
s(\rho, z_1,z_2)=
\begin{cases}
\rho ( z_2z_1^*-|z_1|^2-|z_2|^2)& \text{for $\rho >0$},\\[1ex]
|\rho| ( z_1z_2^*-|z_1|^2-|z_2|^2)& \text{for $\rho <0$}.
\end{cases}
\end{equation}
The proposition is proved.
\end{proof}

\subsection{The spaces $L^{\lambda }(\rho )$ and their decomposition into
${{P}}_0$-invariant subspaces}

In what follows, we will use model $C$ of the complementary series representations.

\begin{DEF*}{} Denote by $L^{\lambda }_+$ and $L^{\lambda }_-$
the subspaces of functions $\psi(\rho ,z)\in L^{\lambda }$
concentrated, with respect to $\rho $, on the half-lines
$\rho >0$ and $\rho <0$, respectively.
\end{DEF*}

By definition, $L^{\lambda }_+$ and $L^{\lambda }_-$ are Hilbert spaces 
with the norms
$$
\|\psi\|^2_+=\int_0^{+ \infty }\|\psi\|^2_{\rho }\,
\rho ^{ \lambda -1}\,d \rho \quad \text{and} \quad
\|\psi\|^2_-=\int_{- \infty }^0 \|\psi\|^2_{\rho }\,
|\rho| ^{ \lambda -1}\,d \rho,
$$
and the space $L^{\lambda } $  is their direct sum:
$$
L^{\lambda }=L^{\lambda }_+ \bigoplus L^{\lambda }_-.
$$
It also follows that $L^{\lambda  }_+$ and $L^{\lambda  }_-$
are direct integrals of the Hilbert spaces  $L^{\lambda  }(\rho )$
of functions $\psi(z)$ on $\CC^{n-1}$ with the norms $\|\psi\|_{\rho }$:
$$
L^{\lambda }_+=\int_0^{+ \infty } L^{\lambda  }(\rho )\,
\rho ^{ \lambda -1}\,d \rho \quad \text{and} \quad
L^{\lambda }_-=\int_{- \infty }^0 L^{\lambda  }(\rho )\,
|\rho| ^{ \lambda -1}\,d \rho .
$$
The explicit formulas for the operators of $T^\lambda$ in model $C$ imply
the following proposition.

\begin{PROP}{}\label{PROP:811-55}\

\begin{enumerate}
\item
The subspaces $L^{\lambda }_+$ and $L^{\lambda }_-$ are invariant under the subgroup
${{P}}$.

\item
The subspaces $L^{\lambda }(\rho )$
 are invariant under the subgroup ${{P}}_0.$
\end{enumerate}
\end{PROP}

According to \eqref{78-25-1}, \eqref{78-24-1}, the formulas for the operators
$T_d$, $d=\diag(\epsilon, u, \epsilon ) \in D$,
and $T_h$, $h = (t_0,z_0) \in H$, do not depend on $\lambda $ and have the form
$$
T_d \,\psi (z)  = \psi (\epsilon \,z\,u ),
\qquad
T_h \psi (z)  =
\begin{cases}{}
e^{\,\rho (it_0 - \frac{|z_0|^2}{2} - zz_0^*)} \,\psi (z+z_0)
&\text{for $\rho>0$,}
\\[1ex]
e^{\,i\rho t_0 - |\rho | (\frac{|z_0|^2}{2} - zz_0^*)} \,\psi (z+z_0)
&\text{for $\rho<0$.}
\end{cases}
$$

Let us compare the representations of the group ${{P}}_0$ 
in $L^{\lambda  } (\rho )$
and $H(\rho )$, where $H(\rho )$ is the space introduced in Section~\ref{sect:3-1}.

According to Proposition~\ref{PROP:3-1}, the space $H(\rho )$
decomposes into a direct sum of  pairwise nonequivalent
$P_0$-irreducible subspaces $H_m(\rho )$:
$$
H(\rho )= \bigoplus _{m=0}^\infty H_m (\rho ).
$$

\begin{PROP}{}\label{PROP:48-96} The representations of ${{P}}_0$ in the spaces
$L^{\lambda }(\rho )$ and $H(\rho )$ are equivalent.
\end{PROP}

\begin{proof}{} Note that the space $L^{\lambda }(\rho )$
is linearly generated by the same polynomials in
$z$ and $\OVER z$ as $H(\rho )$, and the formulas for the
operators of the subgroup ${{P}}_0$ coincide. It follows
that $L^{\lambda }(\rho )$ is equivalent to a subspace in
$H(\rho )$. The isomorphism of $L^{\lambda }(\rho )$ and
$H(\rho ) $ follows from the fact that the norm on $L^{\lambda }(\rho )$
is nondegenerate.
\end{proof}

\begin{COR}{}
The space $L^{\lambda }(\rho )$ decomposes into a direct sum
of pairwise nonequivalent
$P_0$-irreducible subspaces $L^{\lambda }_m(\rho )$:
$$
L^{\lambda }(\rho )= \bigoplus _{m=0}^\infty L^{\lambda }_m (\rho ),
$$
where $L^{\lambda }_m(\rho ) $ is the preimage of $H_m(\rho ) $
under the isomorphism $J:L^{\lambda }(\rho ) \to H(\rho ) $.
\end{COR}

\begin{COR}{}
The intertwining operator $J:L^{\lambda }(\rho )\to H(\rho )$ 
is a multiple of the identity operator on each subspace
$L^{\lambda }_m( \rho  )$, i.e.,
$$
Jf = c_m \,f \quad \text{for every}\quad f \in L_m(\rho ),
$$
where $c_m$ does not depend on $f$.
\end{COR}

In view of the isomorphism $L^{\lambda }(\rho ) \to H(\rho )$,
as an orthogonal basis in $L^{\lambda }( \rho  )$ we can take the same vectors
$f_{p,q}(\rho ,z)$, $p,q \in \ZZ^{n-1}_+$, as in $H(\rho )$
(Theorem~\ref{THM:3-3}). It was mentioned there that
$$
\|f_{p,q}\|^2=\frac{\rho ^{|q|}}{q!}\,\|f_p\|^2
\quad \text{for any $p,q$, where $f_p = f_{p,0}$}.
$$

\begin{PROP}{}\label{PROP:49-98} The norm of the vector $f_p$ in the space
$L^{\lambda }(\rho )$ equals
\begin{equation}{}\label{44-44-1}
\|\, f_p \, \|^{\,2}_{L} = const\,\rho ^{-|p|}\,p!\,
        \frac{
        \Gamma (n - \lambda )\,\Gamma (|p| + \frac{\lambda }{2})
        }{
        \Gamma (|p| + n - \frac{\lambda }{2})\,\Gamma^{\,2} (\frac{\lambda }{2})
        }.
\end{equation}
\end{PROP}

For a proof of formula \eqref{44-44-1}, see Appendix~2. This
formula implies the following corollary.

\begin{COR*}{}
The coefficient $c_m$ in the definition of the intertwining operator
$J:L^{\lambda }(\rho )\to H(\rho )$ is given by
\begin{equation}{}\label{312-12}
c_m^2 = c\,
        \frac{
        \Gamma (m + n - \frac{\lambda }{2})\,\Gamma^{\,2} (\frac{\lambda }{2})
        }{
        \Gamma (n - \lambda )\,\Gamma (m + \frac{\lambda }{2})
        },
\end{equation}
where $c$ does not depend on $m$.
\end{COR*}

\begin{proof}[Proof of the corollary]{} By definition,
$
c_m^2 = \frac{ \|f\|^{\,2}_{H} }{ \|f\|^{\,2}_{L} }
$
for every $f \in L_m^{\lambda }(\rho )$,
where $\|f\|^{\,2}_{H}$ and $\|f\|^{\,2}_{L}$ are the norms of
$f$ in the subspaces $H_m(\rho )$ and $L^{\lambda }_m(\rho )$, respectively.
Let $f = f_p$, $|p| = m$, i.e.,
$$
f_p =
\begin{cases}{}
\OVER{z}^p&\text{for $\rho > 0$,}\\[1ex]
z^p&\text{for $\rho < 0$.}\end{cases}
$$
In this case, \eqref{312-12} immediately follows from the expression
\eqref{44-44-1} for $\|\, f_p \, \|^{\,2}_{L}$ and the equality
$\|\, f_p \, \|^{\,2}_{H} = p!\,\rho ^{-|p|}$.
\end{proof}

\subsection{Restriction of the complementary series representations to the 
subgroup $P$}
The isomorphisms $L^{\lambda }(\rho ) \to H(\rho )$ induce an isomorphism
of Hilbert spaces
$$
J:
L^{\lambda }=\int_{- \infty }^{+ \infty } L^{\lambda  }(\rho )\,
|\rho| ^{ \lambda -1}\,d \rho
\quad  \longrightarrow \quad
H^\lambda =
\int_{- \infty }^{+ \infty } H(\rho )\,|\rho| ^{ \lambda -1}\,d \rho
$$
commuting with the representation of the subgroup ${{P}}_0$. It follows
from the explicit formulas for the operators of the subgroup
${{P}}={{P}}_0\leftthreetimes D_1$ that this isomorphism commutes with
the representations of the whole group ${{P}}$ in these spaces.
Recall that $H^{\lambda}$ decomposes into a direct sum of 
pairwise nonequivalent $P$-irreducible subspaces:
$H^{\lambda } = \bigoplus _{m = 0}^{\infty } (H^{\lambda }_\pm)_m$.
This implies the following proposition.

\begin{PROP}{}\label{PROP:4591}
The subspaces $L^{\lambda }_+$ and  $L^{\lambda }_-$ decompose into direct
sums of pairwise nonequivalent $P$-irreducible subspaces
$$
L^{\lambda }_\pm = \bigoplus _{m = 0}^{\infty } (L^{\lambda }_\pm)_m\,,
$$
where
$$
(L^{\lambda }_+)_m =
\int_0^{\infty } L^\lambda _m (\rho ) \, \rho ^{ \lambda -1}\,d \rho,
\qquad
(L^{\lambda }_-)_m =
\int_{-\infty }^0 L^\lambda _m (\rho ) \, |\rho| ^{ \lambda -1}\,d \rho.
$$
\end{PROP}

In view of the isomorphism $L^{\lambda } \to H^\lambda $, we can transfer
the action of the whole group  $U(n,1)$ from
$L^{\lambda }$ to $H^\lambda $.

\begin{DEF*}{}
The realization of the complementary series representations in
the spaces $H^\lambda $ will be called the Bargmann model of these representations.
\end{DEF*}

We should mention two particular features of the Bargmann model:
\begin{enumerate}
\item
In the Bargmann model, it is essential that the irreducible
${{P}}_0$-invariant subspaces $H(\rho )$ occurring in the decomposition of 
$H^\lambda $ do not depend on $\lambda $.
\item
In the Bargmann model, the action of the operators of the subgroup
${{P}}$ is defined from the beginning; hence in order to describe a
representation of the whole group $U(n,1)$, it suffices to define only 
the action of the operator $T_s$.
\end{enumerate}

To describe the operators $T_s$ in the Bargmann model, we must use the decomposition
$L^\lambda = \bigoplus_{m=0}^\infty
( (L^\lambda_+) _m  \oplus  (L^\lambda_-) _m )$ of the space $L^\lambda $
into   pairwise nonequivalent $P$-irreducible subspaces.
Let $f^\pm _m \in  (L^\lambda_\pm) _m $
be the components of a vector $f \in L^\lambda $
in the subspaces $(L^\lambda_\pm) _m$. Let us identify the elements
from $L^\lambda $ and $H^\lambda $ and denote by
$T^L_s$ and $T^H_s$ the corresponding operators $T_s$.
Then the action of $T^H_s$ in $H^\lambda $ is given by the formula
$$
T^H_s \,f = \sum _m c_m \,((T^L_s\,f)^+ _m + (T^L_s\,f)^- _m  )\,,
$$
where the coefficients $c_m$ are given by \eqref{312-12}.

\subsection{Vacuum vectors and spherical functions}

\begin{PROP}{}\label{PROP:35-45} Up to a factor, the vacuum vector
$f_\lambda$ in the space $L^{\lambda }$ 
has the form
\begin{align}{}
f_{\lambda }(\omega ) &= const  && \text{in model} \quad A,
\\
\label{75-111}
f_{\lambda }(t,z)  &
=(t^2+(1+\tfrac{|z|^2}{2})^2)^{\frac{\lambda }{2}-n}
&&\text{in model} \quad B,
\\
\label{75-171}
f_{\lambda }(\rho,z) &
=|\rho |^{-\frac{\lambda-1}{2}}\,\, e^{\frac{|\rho |}{2}|z|^2}\,
(1+\tfrac{|z|^2}{2})^{\frac{\lambda+1}{2}-n}\,
K_{n-\frac{\lambda+1}{2}}(|\rho |\,(1+\tfrac{|z|^2}{2}))
&&\text{in model} \quad C.
\end{align}
\end{PROP}

\begin{proof}{} For model $A$, the assertion follows from the explicit formulas
for the operators of the representation. When we pass to model $B$, the vector
$const$ goes to the vector \eqref{75-111}. Finally, when we pass to model $C$,
the vector \eqref{75-111} goes to the vector
$$
 f_{\lambda }(\rho,z)=|\rho |^{-n+1}\,\,e^{\frac{|\rho |}{2}|z|^2}\,
\int_0^\infty e^{-i \rho t}\,
(t^2 + ( 1+\frac{|z|^2}{2})^2\,)^{\frac{\lambda}{2}-n}\,dt.
$$
According to \cite[p.~11(7)]{Bateman-Tables}, the obtained expression
coincides, up to a factor, with \eqref{75-171}.
\end{proof}

\begin{PROP}{}\label{PROP:734-9} In model $A$ of the representation
$T^{\lambda }$, the spherical function $\psi_{\lambda }(g)$ is given by the formula
\begin{equation}{}\label{751-43}
\psi_{\lambda }(g)=|\delta |^{- \lambda }\,
F(\frac{\lambda }{2},\frac{\lambda }{2},n,1-|\delta |^{-2} ) \quad
\text{for} \quad
g=\begin{pmatrix}\alpha  &\beta  \\ \gamma   &\delta \end{pmatrix},
\end{equation}
where $F(a,b,c,x)$ is the Gauss hypergeometric function, see \cite{BE}.
\end{PROP}

\begin{proof}{} By definition,
$$
\psi_{\lambda }(g)=c_{\lambda }\,
\int_{S \times S} |\omega \beta + \delta |^{\lambda -2n }\,
|1-\langle \omega ,\omega ' \rangle\,^{1 - \lambda }
\,d \omega\,d \omega'
\quad \text{for} \quad
g = \begin{pmatrix}\alpha &\beta  \\ \gamma &\delta \end{pmatrix}.
$$
Similarly to the case of $O(n,1)$, this expression can be reduced to the form
$$
\psi_{\lambda }(g)=c_{\lambda }\,
\int_S |\omega_1\sh\tau+\ch\tau|^{\lambda-2n}\,d \omega=c_{\lambda }\,
(\ch\tau)^{\lambda-2n}\,
\int_S |1+a\omega_1|^{\lambda-2n}\,d \omega,
$$
where $\ch\tau=|\delta |,$ $a=\frac{\sh\tau}{\ch\tau}$.

Let us write this integral as a power series in $a$. We have
\begin{gather*}{}
|1+a\omega_1|^{\lambda-2n}=(1+a\omega_1)^{\frac{\lambda }{2}-n}\,
(1+a\OVER\omega_1)^{\frac{\lambda }{2}-n}
\\
=\sum_{k=0}^\infty \frac{\Gamma (n-\frac{\lambda }{2}+k)}
{\Gamma (n-\frac{\lambda }{2})\,k!}\,(-a \omega _1)^k\,
\sum_{l=0}^\infty \frac{\Gamma (n-\frac{\lambda }{2}+l)}
{\Gamma (\frac{n-\lambda }{2})\,l!}\,(-a \OVER\omega _1)^l.
\end{gather*}
Substituting this expression into the integral and taking into account 
that the integrals of the monomials
$\omega _1\OVER \omega  _1^l$ with $k\ne l$ vanish, we obtain
$$
\psi_{\lambda }(g)=c_{\lambda }\,(\ch\tau)^{\lambda-2n}\,
\sum_{k=0}^\infty \frac{\Gamma^2 (n-\frac{\lambda }{2}+k)}{(k!)^2}\,
\Bigl(  \int_S |\omega _1^k|^2\,d \omega \Bigr) \,a^{2k}.
$$
Since $\int_S |\omega _1^k|^2\,d \omega= \frac{k!}{(k+n-1)!}$
(see Appendix~2), it follows that
$$
\psi_{\lambda }(g)=c_{\lambda }\,(\ch\tau)^{\lambda-2n}\,
\sum_{k=0}^\infty \frac{\Gamma^2 (n-\frac{\lambda }{2}+k)}{k!(n+k-1)!}\,
a^{2k}=c_{\lambda }\,(\ch\tau)^{\lambda-2n}\,
F(n-\tfrac{\lambda }{2},n-\tfrac{\lambda }{2}n,a^2),
$$
where $F(a,b,c,x)$ is the Gauss hypergeometric function. Since
$\psi_{\lambda }(e)=1$, we have $c_{\lambda }=1$. Further, in view
of the Kummer relations for the Gauss function (see \cite{BE}), we have
$$
F(n-\tfrac{\lambda }{2},n-\tfrac{\lambda }{2},n,a^2)
= (1-a^2)^{\lambda -n}\,F(\tfrac{\lambda }{2},\tfrac{\lambda }{2},n,a^2).
$$
This implies \eqref{751-43}.
\end{proof}

\begin{REM*}{} The above relation for the Gauss function is equivalent 
to the assertion that the spherical functions of the representations with
parameters
$\lambda $ and $2n - \lambda $ coincide (since these representations are equivalent).
\end{REM*}

\section{Models of the special representations of the group $U(n,1)$}\label{sect:6}

\subsection{Definition of the special representations}
Similarly to the case of $O(n,1)$, the special representations of the group
$U(n,1)$ arise as the  
$\lambda \to 0$ or $\lambda \to 2n$ limits
of the complementary series representations $T^{\lambda }$.

Let us describe the cases $\lambda =0$ and $\lambda =2n$ separately.

For $\lambda =0$, the norm $\|f\|_{\lambda }$ degenerates on a subspace
$L^0$ of codimension $1$, and on this space we can introduce the norm
$$
\|f\|^2=\lim_{\lambda \to 0}\frac{d\|f\|^2_{\lambda }}{d \lambda } .
$$
In contrast to the case of $O(n,1)$, this norm for ${{n>1}}$ degenerates on some
infinite-dimensional subspace $L_0 \subset L^0$.

\begin{THM}{}\label{THM:354-2}
The completion of the quotient space $L^0/L_0$ with respect to the norm 
$\|f\|$ splits into a direct sum of two irreducible nonequivalent subspaces.
The representations of the group $U(n,1)$ in these spaces are special 
representations.
\end{THM}

Note that on the subspace $L_0$ there is a nondegenerate norm, given by
$$
\|f\|^2_0=\lim_{\lambda \to 0}\frac{d^{\,2}\|f\|^2_{\lambda }}{d \lambda ^2},
$$
and the representation of $U(n,1)$ in the Hilbert space with
this norm is irreducible.

The space $L$, which is the $\lambda \to 2n$
limit of the spaces $L^{\lambda }$ of complementary 
series representations, has a dual
structure. Namely, $\| f\|_{\lambda }$ degenerates for $\lambda =2n$ and
${{{n>1}}}$ on a subspace $L_{2n}$ of infinite codimension. In the space
$L_{2n}$ there is a nondegenerate norm
$$
\|f\|^2=\lim_{\lambda \to 0}\frac{d\|f\|^2_{\lambda }}{d \lambda } ,
$$
which is degenerate on a one-dimensional subspace $L^{2n}$.

\begin{THM}{}\label{THM:354-3}
The completion of the quotient space $L_{ 2n }/L^{ 2n }$ with respect to this norm
splits into a direct sum of two irreducible nonequivalent subspaces.
The representations of the group $U(n,1)$ in these spaces are special
representations equivalent to those defined above for $\lambda =0$.
\end{THM}

Let us describe several realizations of the special representations.

\subsection{Model $A$: realization in a space of functions on the unitary sphere 
$S \subset \CC^n$ as $\lambda \to 0$ }\label{sect:62}

We assume that the group $U(n,1)$ is realized in matrix model  $a$.

Denote by $L$ the subspace of functions $f(\omega )$ on the unitary unit sphere
$S \subset \CC^n$ satisfying the condition 
$$
\int_S f(\omega )\,d \omega =0,
$$
where $d \omega $ is the invariant measure on $S$. Let us define an action 
of the group $U(n,1)$ on $L$ by the formula
$$
T_g f(\omega )=f(\omega g)\,|b(\omega,g)|^{-2n},
$$
where
\begin{equation}{}\label{822-3}
\omega g=(\omega \beta + \delta )^{-1}\,(\omega \alpha + \gamma ),\quad
b(\omega,g)= \omega \beta + \delta \quad \text{for}
\quad g=\begin{pmatrix}\alpha  &\beta  \\ \gamma   &\delta \end{pmatrix}.
\end{equation}
The fact that $L$ is invariant under the operators
$T_g$ follows from the equality
$d (\omega  g)={|\omega \beta + \delta |^{-2n}\,d \omega }$.

We introduce on $L$ the following scalar products:
\begin{equation}{}\label{54-1}
\langle f_1,f_2 \rangle\,_-=-\int_{S \times S}
\log(1- \langle \omega , \omega ' \rangle)\,
f_1(\omega )\,\OVER{f_2(\omega ')}\,
d \omega \, d \omega ',
\end{equation}
\begin{equation}{}\label{54-2}
\langle f_1,f_2 \rangle\,_+=-\int_{S \times S}
\log(1- \OVER{\langle \omega , \omega ' \rangle\,})\,
f_1(\omega )\,\OVER{f_2(\omega ')}\,
d \omega \, d \omega ',
\end{equation}
\begin{equation}{}\label{54-3}
\langle f_1,f_2 \rangle\,=\langle f_1,f_2 \rangle\,_+
+\langle f_1,f_2 \rangle\,_-
= -2\,\int_{S \times S}
\log|1- \langle \omega , \omega ' \rangle\,|\,
f_1(\omega )\,\OVER{f_2(\omega ')}\,
d \omega \, d \omega ',
\end{equation}
where $\langle \omega , \omega ' \rangle\,
= \omega _1 \bar\omega '_1 + \ldots +\omega _n \bar\omega '_n$.

Note that $\langle f_1,f_2 \rangle\,$ is obtained by passing to the limit
from the scalar products in the spaces of complementary series representations.

\begin{PROP}{}\label{PROP:4-3} The norms
$\|f\|_{\pm}= \langle f,f \rangle\,_{\pm}^{1/2} $
can be written in the form
\begin{equation}{}
\begin{gathered}{}\label{4-4}
\|f\|^2_-=\sum_{k>0}
\frac{(|k|-1)!}{k!}\Bigl|\int_S \omega ^k\,f(\omega )\,d \omega \Bigr|^2,
\\
\|f\|^2_+=\sum_{k>0}
\frac{(|k|-1)!}{k!}\Bigl|\int_S \OVER{\omega} ^k\,f(\omega )\,d \omega \Bigr|^2.
\end{gathered}
\end{equation}
Here we use the notation $k!=k_1! \ldots k_n!$,\, $|k|=k_1+ \ldots +k_n,$\,
$\omega ^k= \omega _1^{k_1} \ldots \omega _n^{k_n}$.
\end{PROP}

Indeed, in order to prove, for example, the first equality, it suffices
to use the series expansion
$$
\log (1-  \langle \omega , \omega ' \rangle)=
-\sum_{k>0}\frac{1}{|k|}\,\frac{|k|!}{k!}\, \omega ^k \OVER{\omega '}^k.
$$

\begin{COR*}{}\ 

1) The scalar products  $\langle f,f \rangle\,_{\pm}$
are sign definite on $L$.

2) The conditions $\|f\|_+=0$ and $\|f\|_-=0$ are equivalent, respectively,
to the conditions
$$
\int_S \omega ^k\,f(\omega )\,d \omega =0 \quad \text{and}\quad
\int_S \bar\omega ^k\,f(\omega )\,d \omega =0 \quad \text{for all} \quad k>0.
$$
\end{COR*}

\begin{PROP}{}\label{PROP:4-2} The operators $T_g$ preserve the norms 
$\|f\|_{\pm}$, i.e., $\|T_g f\|_{\pm}=\|f\|_{\pm}$ for all $g\in U(n,1)$.
\end{PROP}

\begin{proof}{}  Since
$|\omega \beta + \delta|^{-2n}\,d \omega =d(\omega  g)$, it follows that 
\begin{gather*}{}
\|T_g f\|_-^{\,2} = -\int_{S \times S}\log (1-  \langle \omega , \omega ' \rangle)\,
f(\omega g)\,\OVER{f(\omega 'g)}\,d( \omega g)\,d (\omega 'g)=
\\
-\int_{S \times S}\log (1-  \langle \omega g^{-1}, \omega 'g^{-1} \rangle)
\,f(\omega )\,\OVER{f(\omega ')}\,d \omega \,d \omega '.
\end{gather*}
Note that
$$
1-  \langle \omega g^{-1}, \omega 'g^{-1} \rangle\,=
(1-  \langle \omega , \omega ' \rangle)\,b^{-1}(\omega ,g^{-1})\,
\OVER{b^{-1}(\omega' ,g^{-1})}.
$$
Hence
$$
\|T_g f\|_-^2 = \| f\|_-^2 - \int_{S \times S}
(
\log b(\omega g^{-1}) + \log \OVER{b(\omega' g^{-1})}\,
)
f(\omega )\,\OVER{f(\omega ')}\,d \omega \,d \omega ' .
$$
Since $f\in L$, the second term in this sum vanishes.

For the norm $\|f\|_+$, the proof is the same.
\end{proof}

\begin{COR*}{} The following subspaces in $L$ are invariant under the action
of $U(n,1)$:
$$
L_+= \{f\in L \mid \|f\|_-=0 \};\quad
L_-= \{f\in L \mid \|f\|_+=0 \};\quad  
L_0=L_+ \cap  L_-.
$$
\end{COR*}

Let $L_0 \subset L$ be the subspace of functions $f$ with $\|f\|=0$.
By definition, the norm $\|f\|$ is nondegenerate on the quotient space $L/L_0$.

\begin{DEF*}{} Denote by $\Cal L$, $\Cal L_+$, and $\Cal L_-$ 
the complex Hilbert spaces obtained by completing
$L/L_0$, $L_+/L_0$, and $L_-/L_0$ with respect to the norm $\|f\|$. \end{DEF*}

Obviously,
$$
\Cal L=\Cal L_+ \oplus \Cal L_-
$$
and
$$
\|f\|=\|f\|_+ \quad \text{on} \quad \Cal L_+;\quad
\|f\|=\|f\|_- \quad \text{on} \quad \Cal L_-.
$$

As a result, we obtain the following theorem.

\begin{THM}{}\label{THM:683-5}
In model $A$ with $\lambda =0$ of the two special representations
$T^+$ and $T^-$ of the group $U(n,1)$,
\begin{enumerate}
\item
The direct sum of $T^+$ and $T^-$ is realized
in the Hilbert space $\Cal L$, the completion with respect to the norm
$
\| f\|^2 = -2\,\int_{S \times S}
\log|1- \langle \omega , \omega ' \rangle\,|\,
f(\omega )\,\OVER{f(\omega ')}\,
d \omega \, d \omega '
$
of the quotient space $L/L_0$, where $L$ is the space of functions
$f(\omega ) $ on the sphere $S=S^{2n-1}$ with 
$\int _S f(\omega ) \,d \omega = 0$ and
$L_0$ is the subspace of functions with norm $0$.

\item
The subspaces $\Cal L_+$ and $\Cal L_-$ of the special representations
are determined by the conditions
$\Cal L_{\pm} = \{f \in \Cal L\,\, \mid\,\,\, \| f \|_{\mp} =0\}$,
where $\| f \|_{+}$ and $\| f \|_{-}$ are given by Proposition~\ref{PROP:4-3}.

\item
The operators of the group $U(n,1)$ act in $\Cal L_{\pm}$
as follows:
$$
T^{\pm}_g f(\omega )=f(\omega g)\,|b(\omega,g)|^{-2n} \pmod {L_0},
$$
where
\enskip
$
\omega g=(\omega \beta + \delta )^{-1}\,(\omega \alpha + \gamma ),\quad
b(\omega,g)= \omega \beta + \delta \quad \text{for}
\quad g=\begin{pmatrix}\alpha  &\beta  \\ \gamma   &\delta \end{pmatrix}.
$
\end{enumerate}
\end{THM}

\subsection{The structure of the special representations in model
$A$ for $\lambda =0$}

Let us describe the subspaces $L$, $L_+$, $L_-$, and $L_0$
defined in Section~\ref{sect:62} in terms of generating elements.
Note that the original space $\TIL L$ is linearly generated by the monomials 
$\omega ^k \bar \omega ^l=\prod \omega _i^{k_i} \bar \omega _i^{l_i}$.
Let us write $k>l$ if $k_i\ge l_i$ and $k\ne l$. Further, set
$$
c_k=\int_S |\omega ^k|^2\,d \omega , \quad k\in\ZZ_+^n;
$$
an explicit expression for $c_k$ will be given below.

The description of the norms $\|f\|_{\pm}$
implies the following proposition.

\begin{PROP}{}\label{PROP:55-1} \
\begin{enumerate}
\item If $k>l$, then $\|\omega ^k\bar \omega ^l\|_+
\ne 0$ and $\|\omega ^k\bar \omega ^l\|_-=0$.

\item If $l>k$, then
$\|\omega ^k\bar \omega ^l\|_-\ne 0$ and $\|\omega ^k\bar \omega ^l\|_+=0.$

\item If $k\ngtr l$, $l\ngtr k$, $k\ne l$, then
$\|c_{k+l} \omega ^k- c_k \omega ^{k+l}\bar \omega ^l\|_{\pm}=
\|c_{k+l}\bar\omega ^k- c_k \omega ^l\bar\omega ^{k+l} \|_{\pm}=0$.
\end{enumerate}
(Here $k,l\in\ZZ_+^n$.)
\end{PROP}

\begin{COR*}{}\
\begin{enumerate}\item The subspace $L$ is linearly generated
by the monomials $\omega ^k\bar \omega ^l$, $k\ne l$, and the binomials
$c_k- c_0 |\omega ^k|^2$.
\item The subspace $L_0$  is linearly generated
by the monomials $\omega ^k\bar \omega ^l$, $k\ne l$, where 
$k\ngtr l$, $l\ngtr k$, $k\ne l$, and the binomials
$c_{k+l} \omega ^k- c_k \omega ^{k+l}\bar \omega ^l$,
$c_{k+l}\bar\omega ^k- c_k \omega ^l\bar\omega ^{k+l}$,
$k,l\in\ZZ_+^n$.
\item The quotient spaces $L_+/L_0$ and $L_-/L_0$ are linearly generated
by the monomials $\omega ^k$, $k>0$, and 
$\bar \omega ^k$, $k>0$, respectively.
\end{enumerate}
\end{COR*}

Note that $L_0=0$ for $n=1$.

\begin{PROP}{}\label{PROP:4-44} The monomials $\omega ^k$ and
$\bar \omega ^k$, $k>0$, form orthogonal bases in $\Cal L_+$ and $\Cal L_-$,
respectively; their norms equal
\begin{equation}{}\label{4-44}
\|\omega ^k\|^2_+=\|\bar \omega ^k\|^2_-=
\frac{4\pi^{2n}\,k!\,(|k|-1)!}{[(|k|+n-1)!]^2}.
\end{equation}
\end{PROP}

\begin{proof}{} The fact that these monomials really generate the corresponding
spaces follows from Proposition~\ref{PROP:55-1} (corollary).
The orthogonality follows from the explicit formulas
\eqref{4-4} for the norm. The same explicit formulas and
\eqref{953--7} imply that
$$
\|\omega ^k\|^2_+=\|\bar \omega ^k\|^2_-=
\frac{(|k|-1)!}{k!}\Bigl|\int_S |\omega ^k|^2\,d \omega\Bigr|^2.
$$
Since
$$
\int_S |\omega ^k|^2\,d \omega =\frac{2\,\pi^n\,k!}{(|k|+n-1)!},\quad
k\in\ZZ_+^n,
$$
this implies \eqref{953--7}. (For a proof of the latter formula, see
Appendix~2.)
\end{proof}

\begin{COR*}{} Elements of the spaces $\Cal L_+$ and $\Cal L_-$ 
can be interpreted as the boundary values of regular 
analytic (respectively, antianalytic) functions in the unit ball in
$\CC^n$. In this interpretation, if a function $f$ in the ball
is given by the series  $f(z)=\sum_{k>0} a_k z^k$, then
$\|f\|_+^2=4\pi^{2n}\,\sum_{k>0}\frac{ k!(|k|-1)!}{[(|k|+n-1)!]^2}\,|a_k|^2.$
\end{COR*}

\subsection{Formulas for the $1$-cocycles $U(n,1)\to\Cal L$ and
$U(n,1)\to\Cal L_{\pm}$ in model $A$}
Let us find explicit formulas for the nontrivial $1$-cocycles
on the group $U(n,1)$ associated with the two special representations
of $U(n,1)$ in the spaces $\Cal L_+$ and $\Cal L_-$.

\smallskip

Note that for $\lambda =0$ the nontrivial $1$-cocycle
$a:\,U(n,1)\to\Cal L = \Cal L_+ \oplus \Cal L_- $
is given by the formula
\begin{equation}{}\label{4-10}
a(\omega  ,g)=T_g f_0- f_0=
|\omega \beta + \delta |^{-2n}-1 \quad \text{for}
\quad g=\begin{pmatrix}\alpha  &\beta  \\ \gamma   &\delta \end{pmatrix}
\end{equation}
(since the vector ${f_0(\omega )\equiv 1}$ is invariant under
the maximal compact subgroup of $U(n,1)$ and $f_0\notin\Cal L$).

Obviously, the $1$-cocycles $a^{\pm}:\,U(n,1)\to\Cal L_{\pm}$ associated
with the special representations are the projections of $a$
to the corresponding subspaces:
\begin{equation}{}\label{717-43}
a(\omega  ,g)=a^+(\omega  ,g)+a^-(\omega  ,g),\quad a^{\pm}\in\Cal L_{\pm}.
\end{equation}
Let us find explicit expressions for these cocycles for $\lambda =0$ and $\lambda =2n$.

\begin{PROP}{}\label{PROP:4-10} For $\lambda =0$, the nontrivial $1$-cocycles
$a^{\pm}(\omega  ,g)$ are given (modulo $L_0$) by the following formula:
\begin{equation}{}\label{4-11}
a^+(\omega  ,g)=\OVER{a^-(\omega  ,g)} =
\delta ^n\,(\omega \beta + \delta )^{-n}-1.
\end{equation}
In particular, $a^+(\omega  ,g)=-\frac{\omega \beta }{\omega \beta+ \delta }$
for $n=1$.
\end{PROP}

\begin{proof}{} Setting $\gamma = \beta \delta ^{-1}$, we have
\begin{align*}{}
a(\omega ,g)&=|\delta |^{-2n}\, (1+ \omega \gamma )^{-n}\,
(1+ \bar\omega \bar\gamma )^{-n}-1
\\
&=|\delta |^{-2n}\,\sum_{k,l}(-1)^{|k+l|}\,
\frac{(|k|+n-1)!\,(|l|+n-1)!}{k!\,l!}\, \omega ^k \gamma ^ k
\bar \omega ^l\bar \gamma ^l -1.
\end{align*}
Since $\omega ^k\bar \omega ^l\ne 0$ on $\Cal L_+$ only for
$k>l$, it follows that
\begin{equation}{}\label{968-3}
a^+(\omega ,g)=|\delta |^{-2n}\, \sum_{k>0}\sum_l (-1)^{|k|}\,
\frac{(|k+l|+n-1)!\,(|l|+n-1)!}{(k+l)!\,l!}\, \omega ^{k+l} \bar \omega ^l\
\gamma  ^ {k+l}\bar \gamma ^l.
\end{equation}
According to Proposition~\ref{PROP:55-1}, on $\Cal L_+$ we have
$$
\omega ^{k+l}\bar \omega ^l=\frac{c_{k+l}}{c_k}\, \omega ^k=
\frac{(k+l)!}{(|k+l|+n-1)!}\frac{(|k|+n-1)!}{k!}\, \omega ^k.
$$
Substituting this expression into \eqref{968-3}, we obtain
\begin{align*}{}
a^+(\omega ,g)&=|\delta |^{-2n}\, \sum_{k>0} (-1)^{|k|}\,
\frac{(|k|+n-1)!}{k!}\, \omega ^k \gamma ^k\,
\sum_l \frac{(|l|+n-1)!}{l!}\,\gamma ^l\bar \gamma ^l
\\
&=|\delta |^{-2n}\,((1+ \omega \gamma )^{-n}-1)\,(1 - \gamma^*\gamma )^{-n}.
\end{align*}
Note that $1 - \gamma ^*\gamma =
|\delta |^{-2}(|\delta|^2 - \beta ^*\beta ) = |\delta |^{-2}$. Hence
$$
a^+(\omega ,g) = (1+ \omega \gamma )^{-n}-1 = \delta ^n\,
(\omega \beta + \delta )^{-n}-1.
$$
The proposition is proved.
\end{proof}

\begin{PROP}{}\label{PROP:4-10:2n} For $\lambda = 2n$, the nontrivial $1$-cocycle
$a^{+}(\omega  ,g)$ is given, up to a factor, by the following formula: 
\begin{equation}{}\label{4-11:2n}
a^+(\omega  ,g) = e^{\,- \delta ^{-1}\omega \beta } - 1;
\end{equation}
the cocycle $a^-$ is obtained from $a^+$ by taking the complex conjugate.
\end{PROP}

\begin{proof}{} It follows from the formulas for $a^+(\omega ,g)$ in the case $\lambda =0$
and for the intertwining operator $R:L^0\to L^{2n}$  that
$$
a^+(\omega ,g) = \int_S \log(1-\langle \omega ,\omega ' \rangle)\,
        (1+ \omega ' \gamma )^{-n}\,d \omega ',
$$
where $\gamma = \beta \delta ^{-1}$. Expanding
$\log(1-\langle \omega ,\omega ' \rangle)$ and $(1+ \omega ' \gamma )^{-n}$
into a power series and throwing off the terms with zero integrals, we obtain
\begin{align*}{}
a^+(\omega ,g) &
= \sum _{k>0}
\Bigl[
        \frac{ (-1)^{|k|}\,(n+|k| - 1)\opn{!} }{ (n-1)\opn{!}\,(k\opn{!})^2 }\,
        \omega ^k \,\gamma ^k\, \int_S  |\omega ^{\prime k}|^2\,d \omega '
\Bigr]
\\ &
= \frac{ 2\pi^n }{ (n-1)\opn{!} }
  \sum _{k>0} \frac{ (-1)^{|k|} }{ k\opn{!} }\,\omega ^k \,\gamma ^k
= \frac{ 2\pi^n }{ (n-1)\opn{!} }\,(e^{\,- \omega \,\gamma } - 1).
\end{align*}
The proposition is proved.
\end{proof}

\subsection{The ``strange'' irreducible representation of the group
$U(n,1)$ in the space $\Cal L_0$}

Besides two special representations, there exists another irreducible
unitary representation of the group $U(n,1)$ that can be obtained from
complementary series representations in the
$\lambda \to 0$ or $\lambda \to 2n$ limit (the so-called ``strange representation'',
see \cite{21-1958,19-1957,18-1956}).
Note that there is no analog of this representation for the groups $O(n,1)$.

Let us construct this representation starting from model $A$ of
the complementary series representations for $\lambda \to 0$.

We define a scalar product on the subspace $L_0=L_+\cap L_-$
by the formula
\begin{equation}{}\label{4-5}
\langle f_1,f_2 \rangle\, _0 =
\lim\limits_{\lambda \to 0}
\frac{
        d^2\langle f_1,f_2 \rangle\, _{\lambda }
        }{
        d \lambda^2
        }
=\int_{S \times S}\log^2 |1-  \langle \omega , \omega ' \rangle\,|\,
f_1(\omega )\,\OVER{f_2(\omega ')}\,d \omega \,d \omega '.
\end{equation}

\begin{PROP}{}\label{PROP:4-12} The norm
$\|f\|_0=\langle f_1,f_2 \rangle\, _0^{1/2}$ can be written in the form
\begin{equation}{}\label{4-12}
\|f\|^2_0=\sideset{}{'}\sum_{k,l>0} c_kc_l\,
|\int_S \omega^k \bar\omega  ^l\,f(\omega )\,d \omega |^2,
\end{equation}
where $c_k=\frac{(|k|-1)!}{k!}$. The sum ranges only over multi-indices
$k>0$ and $l>0$ such that neither $k>l$ nor $l>k$.
\end{PROP}

\begin{proof}{} We have
\begin{gather*}{}
\log^2 |1-\langle \omega , \omega ' \rangle\,|=
(\sum_{k>0}c_k \omega ^k \bar{\omega '}^k+
\sum_{l>0}c_l \bar\omega ^l \omega ^{\prime l}  )^2
\\
= \sum_{k,k'>0}c_kc_{k'} \omega ^{k+k'}\bar{\omega '}^{k+k'}+
\sum_{l,l'>0}c_lc_{l'} \bar\omega ^{l+l'} \omega ^{\prime l+l'}
+\sum_{k,l>0}c_k c_l \omega ^k \bar\omega^l
\,\OVER{\omega ^{\prime k} \bar\omega^{\prime l}}.
\end{gather*}
Let us substitute this expression into \eqref{4-5}. Since
$L_0=L_+\cap L_-$, the terms of the first two sums make zero contributions.
The same holds for the terms of the third sum with
$k>l$ or $l>k$. Thus \eqref{4-12} holds.
\end{proof}

\begin{COR*}{} The norm $\|f\|_0$ is positive definite on $L_0$.
\end{COR*}

Denote by $\Cal L_0$ the complex Hilbert space obtained by completing
$L_0$ with respect to the norm $\|f\|_0$.

\begin{PROP}{}\label{PROP:4-6} The norm $\|f\|_0$ is invariant under the operators
$T_g$; thus the operators $T_g$ form a unitary representation of the group
$U(n,1)$ in the space $\Cal L_0$.
\end{PROP}

\begin{PROP}{}\label{PROP:4-7} The representation of the group $U(n,1)$ 
in the space $\Cal L_0$ is irreducible.
\end{PROP}

\subsection{Model $A$ of the special representations in a space of functions 
on the unitary sphere for $\lambda \to  2n$}\label{sect:63}

In model $A$ with $\lambda =2n$, the direct sum
$\Cal L=\Cal L_+ \oplus \Cal L_-$ of the special representations is realized
in the quotient of $L^{2n}/L_{2n}$
by the subspace of constants
with the norm
\begin{equation}{}\label{norma-L/L_0}
\| f \|^2 = (-1)^n \int_{S \times S}  \delta ^{(n)}
(1-\langle \omega ,\omega ' \rangle)\,
f(\omega )\,\OVER{f(\omega ')}\,d \omega \,d \omega ' ,
\end{equation}
where $\delta ^{(n)}$ is the derivative of the delta function on the sphere
$S$, given by the formula 
${ \int_S \delta (1-\langle \omega ,\omega ' \rangle)\,
f(\omega' ) \,d \omega '} = f(\omega )$.
The operators of the representation are given by the formula
$$
(T_gf)(\omega ) = f(\omega g).
$$
An intertwining operator $R$ from model $A$ with $\lambda =0$
to model $A$ with $\lambda = 2n$ has the form
$$
(Rf)(\omega ) = \int_{S}
\log (|1-  \langle \omega , \omega' \rangle\,|)\,
f(\omega' ) \,d \omega '.
$$

Note that the norm in model $A$ with $\lambda =0$ can be written as an integral
over the sphere $S$:
$$
\|f\|^2 =\int_S (Rf)(\omega )\,\OVER{f(\omega )}\,d \omega .
$$

In the model under consideration with $\lambda =2n$, the subspaces
$\Cal L_+$ and $\Cal L_-$ coincide with the subspaces of boundary values
of analytic (respectively, antianalytic) functions in the unit ball
$D \subset \CC^n$: $|z|^2 \equiv |z_1|^2+ \dots +|z_n|^2<1$.

The special representations $T^+$ and $T^-$ can be realized directly in the 
Hilbert spaces of analytic and antianalytic functions on $D$.
For this, it suffices to use the following assertions.

\begin{enumerate}
\item The ball $D$ is a homogeneous space with respect to the action
$z\to zg$ of the group $U(n,1)$, where
$zg = (\omega \beta + \delta )^{-1} (\omega \alpha + \gamma )$.
\item On $D$ there exists a $U(n,1)$-invariant measure, namely,
$(1-|z|^2)^{-n-1}\,d \mu (z)$, where $d \mu (z)$ is the Lebesgue measure on $D$
(see \cite{CH}).
\end{enumerate}

\begin{THM}{}\label{THM:7-34-52} The special representation $T=T^+$
(respectively, $T=T^-$) of the group $U(n,1)$ is equivalent to 
the representation in the quotient of the Hilbert space of analytic 
(respectively, antianalytic) functions
on the ball $D$ by the subspace of constants with the norm
$$
\| f\|^2 = \int_D |f(z)|^2 \, \delta ^{(n)}(1-|z|^2)\,d \mu (z).
$$
The operators of this representation have the form $(T_gf)(z) = f(zg)$.
\end{THM}

Note that the vectors $z^k$, $k \in \ZZ_+^n$, form an orthogonal basis
in this space, with $\|z^k\|^2 = \dfrac{(|k|+n-1)!}{(|k|-1)!}$.

\subsection{Model $B$: realization in a space of functions on
the Heisenberg group $H$}
Model $B$ of the two special representations of the group 
$U(n,1)$ is obtained from model $A$ by passing from functions 
$ F(\omega )$ on the sphere $S$ to functions
$f(h)=f(\zeta ,z)$, $ \zeta =it-\frac{|z|^2}{2}$, on the Heisenberg
group $H$ according to the formulas
$$
f(\zeta ,z)=|1 - \zeta |^{-2n}\,F(\omega ),
$$
where
$$
\omega = \Bigl(\frac{1+\zeta }{1- \zeta }, \frac{\sqrt 2 z_1}{1-\zeta}, \dots ,
\frac{\sqrt 2 z_{n-1}}{1-\zeta}\Bigr).
$$
We will describe model $B$ only for $\lambda =0$
(for $\lambda =2n$ the description is similar). In this description,
we use matrix model $b$ of the group $U(n,1)$
(see~\S\,\ref{sect:1}).

\begin{THM}{}\label{THM:6-7-8-9-10}
In model $B$ of the two special representations $T^+$ and $T^-$ 
of $U(n,1)$ with $\lambda =0$,
\begin{enumerate}
\item
The direct sum of $T^+$ and $T^-$ is realized
in the Hilbert space $\Cal L$, the completion with respect to the norm
\begin{equation}{}\label{5-7}
\|f\|^2=
-\int_{H \times H}\log R(h,h')\,f(h)\OVER{f(h')}\,dh\,dh',
\end{equation}
where $dh$ is the invariant measure on $H$ and
\begin{equation}{}\label{5-8}
R(h,h')=R(\zeta ,z; \zeta ',z')=|\zeta +\bar{\zeta '}+zz^*|,
\end{equation}
of the quotient space $L/L_0$,
where $L$ is the space of functions $f(h ) $ on the group $H$
with $\int _H f(h ) \,d h = 0$
and $L_0$ is the subspace of functions with norm $0$.

\item In the coordinates $h=(\zeta ,z)$ on $H$,
the operators of the representation $T = T^+ \oplus T^-$ have the following form:
\begin{equation}{}\label{5-2}
T_g f(\zeta ,z )=f(\zeta ',z')\,|b(\zeta ,z\,;\,g)|^{-2n} \pmod{L_0},
\end{equation}
where
\begin{equation}{}\label{5-3}
\zeta '=(\zeta g_{13}+ z g_{23}+g_{33})^{-1}\,
(\zeta g_{11}+ z g_{21}+g_{31}),
\end{equation}
\begin{equation}{}\label{5-4}
z'=(\zeta g_{13}+ z g_{23}+g_{33})^{-1}\,
(\zeta g_{12}+ z  g_{22}+g_{32}),
\end{equation}
\begin{equation}{}\label{5-5}
b(\zeta ,z;g)=\zeta g_{13}+ z g_{23}+g_{33}.
\end{equation}

\item
The subspaces $\Cal L_+$ and $\Cal L_-$ of the special representations
are determined by the conditions
$\Cal L_{\pm} = \{f \in \Cal L\,\, \mid\,\,\, \| f \|_{\mp} =0\}$,
where $\|f\|_-$ and $\|f\|_+$ are determined by the following formulas:
\begin{equation}{}\label{5-9}
\begin{aligned}{}
\|f\|^2_-
&
=\sum_{k>0} 2^{|k|-1}
\frac{(|k|-1)!}{ k!}\,
\Bigl|\int_H B_k(h)\,f(h)\,d h \,\Bigr|^2, \quad
\\
\|f\|^2_+
&
=\sum_{k>0} 2^{|k|-1}
\frac{(|k|-1)!}{ k!}\,
\Bigl|\int_H \OVER{B_k(h)}\,f(h )\,d h \,\Bigr|^2,
\end{aligned}
\end{equation}
where $ k=(k_1, \dots ,k_n)$, $|k|=k_1 + \dots + k_n $,
\begin{equation}{}\label{5-10}
B_k(h)\,=\,B_k(\zeta ,z)
\,=\,\frac{(1+ \zeta )^{k_1}}{(1- \zeta )^{|k|}}\,
\cdot z_1^{k_2} \cdots z_{n-1}^{k_n}.
\end{equation}
\end{enumerate}
\end{THM}

\begin{REM*}{} Obviously, the conditions $\|f\|_-=0$ and $\|f\|_+=0$
are equivalent, respectively, to the equations
\begin{equation}{}\label{783-25}
\int_H B_k(h)\,f(h)\,dh=0 \quad \text{and}\quad
\int_H \OVER{B_k(h)}\,f(h)\,dh=0\quad \text{for}\quad k>0.
\end{equation}
\end{REM*}

\begin{PROP}{}\label{PROP:6-7-8-9-11}
In model $B$ with $\lambda =0$, the nontrivial $1$-cocycles
$a^{\pm}:\,U(n,1)\to\Cal L_{\pm}$ associated with the special representations 
of the group $U(n,1)$ are the projections to $\Cal L_+$ and $\Cal L_-$, 
respectively, of the following function $a:U(n,1)\to\Cal L$:
$$
a(\zeta ,z;g)=T_gf_0-f_0, \quad \text{where} \quad f_0=|1- \zeta |^{-2n}.
$$
\end{PROP}

Proposition~\ref{PROP:4-44} implies the following result.

\begin{PROP}{}\label{PROP:270-74} The following functions $f_k(\zeta ,z)$,
$k=(k_0,k_1, \dots ,k_{n-1})\in\ZZ^n_+$, $|k|>0$, form an orthogonal basis in
$\Cal L_+$:
\begin{equation}{}\label{4-44-33}
f_k(\zeta ,z)=(1-\OVER\zeta) ^{-n}\,(1-\zeta )^{-n-|k|}\,
    (1+\zeta) ^{k_0}\,z_1^{k_1}\ldots z_{n-1}^{n-1}.
\end{equation}
A similar basis in $\Cal L_-$ is obtained by replacing $\zeta $ and $z_i$
with their complex conjugates.
\end{PROP}

\subsection{Model $C$: realization in a space of functions
$\psi(\rho ,z)$ on $\RR \times \CC^{n-1}$}
Model $C$ is obtained by passing from functions $f(h)=f(\zeta ,z)$,
$\zeta =it-\frac{|z|^2}{2}$, on $H$ in model $B$ (with $\lambda =0$)
to functions $\psi(\rho ,z)$ according to the formula
$$
\psi(\rho, z ) =|\rho |^{-n+1}\,e^{\frac{|\rho |}{2}|z|^2}\,
\int_{- \infty }^{+ \infty }e^{\,i\rho t}\,
f(t, z )\,d t.
$$

\begin{PROP}{}\label{PROP:5-3} In model $C$, the subspace  $L$ of functions
$\psi(\rho ,z)$ is determined by the relation
\begin{equation}{}\label{5-111}
\int_{\CC^{n-1}}\psi(0,z)\,d\mu(z)=0.
\end{equation}
The norm $\|\psi\|$ on $L$ is given by the formula
\begin{equation}{}\label{497-78}
\|\psi\|^2=
\int_{- \infty }^{+\infty } \|\psi\|^2_{\rho }\,
|\rho |^{ -1}\,d \rho ,
\end{equation}
where
\begin{equation}{}\label{33-337}
\|\psi\|^2_{\rho }=|\rho |^{2n-2}\int_{\CC^{n-1} \times \CC^{n-1}}
e^{|\rho |(a_{\rho }(z_1,z_2)-|z_1|^2-|z_2|^2)}\,
\psi(\rho ,z_1)\,\OVER{\psi(\rho ,z_2)}\,d\mu(z_1)\,d\mu(z_2),
\end{equation}
\begin{equation}{}\label{33-818-2}
a_{\rho }(z_1,z_2)=
\begin{cases}
 z_2z_1^* & \text{for $\rho <0$},\\[1ex]
 z_1z_2^* & \text{for $\rho >0$}.
\end{cases}
\end{equation}
In another way,
\begin{equation}{}\label{53-818}
\|\psi\|^2_{\rho }=
\begin{cases}
\sum_{k\in\ZZ^{n-1}_+}  \frac{|\rho|^{2n+|k|-2} }{k!}\,\Bigl|
\int_{\CC^{n-1}}z^k\,e^{\,-|\rho| zz^*}\,\psi(\rho ,z)\,d\mu(z)\,\Bigr|^2\,
& \text{for $\rho <0$},\\[1.5ex]
\sum_{k\in\ZZ^{n-1}_+}  \frac{\,\rho^{2n+|k|-2} }{k!}\,\Bigl|
\int_{\CC^{n-1}}\OVER z^k\,e^{\,-|\rho| zz^*}\,\psi(\rho ,z)\,d\mu(z)\,\Bigr|^2\,
& \text{for $\rho >0$}.
\end{cases}
\end{equation}
\end{PROP}

\begin{proof}{} Relation \eqref{5-111} is equivalent to relation
\eqref{5-2} in the realization on $H$. To calculate
$\|\psi\|$, we will use the equation
$\|\psi\|^2=\lim\limits_{\lambda \to 0}\frac{d\|\psi\|^2_{\lambda }}{d\lambda }$,
where $\|\psi\|_{\lambda }$ is the norm in the space of the
complementary series representation. Substituting into this equation the explicit expression for
$\|\psi\|^2_{\lambda }$ and taking into account that
$Q_{\lambda }(t)= const$ for $\lambda =0$, we obtain \eqref{497-78}
and \eqref{33-337}. Substituting into \eqref{33-337} the power
series expansion of the function $e^{a_{\rho }(z_1,z_2)}$ for $\rho >0$ and
$\rho <0$,
we obtain \eqref{53-818}.
\end{proof}

\begin{COR*}{} In model $C$, the space $\Cal L$ decomposes into
a direct integral of subspaces:
\begin{equation}{}\label{719-55}
\Cal L=\int_{- \infty }^{+ \infty }\Cal L(\rho )\,|\rho |^{-1}\,d \rho ,
\end{equation}
where $\Cal L(\rho )$ is the quotient of the Hilbert space of functions
$\psi$ on $\CC^{n-1}$ with the norm \eqref{33-337} by
the subspace of functions with norm $0$.
\end{COR*}

\begin{THM}{}\label{THM:5-1} In model $C$, the subspaces $\Cal L_+$ and
$\Cal L_-$ of the special representations coincide with the subspaces of functions
$\psi(\rho ,z)$ concentrated, with respect to  $\rho $, on the half-lines
$\rho >0$ and $\rho <0$, respectively.
\end{THM}

\begin{proof}{} It suffices to prove that $\phi(\rho ,z)\in\Cal L_+$
if $\phi(\rho ,z)=0$ for $\rho <0$, and, similarly, $\phi(\rho ,z)\in\Cal L_-$
if  $\phi(\rho ,z)=0$ for $\rho >0$.

Let us use condition \eqref{783-25}, which describes when a function 
$f(h)=f(t,z)$ in model $B$ belongs to the subspace $\Cal L_+$. Passing
from functions $f(t,z)$ to functions $\psi(\rho ,z)$, we obtain a condition
for $\psi\in\Cal L_+$ in the form 
\begin{equation}{}\label{48-23}
I_k(\psi)\equiv \int_{- \infty }^0\int_{\CC^{n-1}}
A _k(\rho ,z)\,z^{k'}\,e^{\,-|\rho |\,\frac{|z|^2}{2}}\,
\psi(- \rho ,z)\,|\rho |^{n-1}\, d \rho d\mu(z)=0 \quad \text{for}\quad
|k|>0,
\end{equation}
where
\begin{equation}{}\label{48-24}
A_k(\rho ,z) = \int_{- \infty }^\infty
\frac{(1+ \zeta )^{k_0}}{(1- \zeta )^{|k|}}\,e^{\,i\rho  t}\,dt,\quad
\zeta =it-\frac{|z|^2}{2}.
\end{equation}
Let us check that this condition is satisfied if
$\psi(\rho ,z)=0$ for $\rho <0$. Indeed, the function $A_k(\rho ,z)$ 
can be written as a linear combination of the integrals
$$
\int_{- \infty }^\infty \frac{e^{\,i\rho t}\,dt}{(1- \zeta )^m}
= \int_{- \infty }^\infty \frac{e^{\,i\rho t}\,dt}{(a- it)^m},\quad
|k'|\le m\le |k|,
$$
where $a=1+\frac{|z|^2}{2}$. It is known that for $\rho >0$ these integrals vanish
for every $m>0$. Hence if $\psi(\rho ,z)=0$ on the half-line
$\rho <0$, then $I_k(\psi)=0$ for $|k'|>0$. For $|k'|=0$, the expansion
of $A_{k}(\rho ,z)$ contains the term $\int_{- \infty }^\infty
e^{\,-i \rho  t}\,dt= \delta (\rho )$. In this case, the corresponding term in 
$I_k(\psi)$ vanishes because of \eqref{5-111}.

In a similar way one can prove that if $\psi(\rho ,z)\equiv 0$ for $\rho>0$,
then $\psi\in \Cal L_-$.
\end{proof}

\begin{COR*}{} In model $C$, the spaces $\Cal L_+$ and $\Cal L_-$
decompose into direct integrals of subspaces:
$$
\Cal L_+=\int_0^{+ \infty }\Cal L(\rho )\,\rho ^{-1}\,d \rho ,\quad
\Cal L_-=\int_{- \infty }^0\Cal L(\rho )\,|\rho |^{-1}\,d \rho;
$$
the norms in the spaces $\Cal L_+$ and $\Cal L_-$, respectively, are given by
\begin{equation}{}\label{573-15}
\aligned
\|\psi\|_+^2
&
=\int_0^{+\infty } \|\psi\|^2_{\rho }\,
\rho ^{ -1}\,d \rho ,
& \qquad 
\|\psi\|_-^2
&
=\int_{- \infty }^0 \|\psi\|^2_{\rho }\,
|\rho |^{ -1}\,d \rho , 
\endaligned
\end{equation}
where $\| \psi \|_{\,\rho }$ is given by \eqref{33-337} or
\eqref{53-818}. In more detail,
\begin{equation}{}
\begin{aligned}{}\label{5-123}
\|\psi\|_-^2 &
=\sum_{k\in\ZZ^{n-1}_+}  \frac{1}{k!}\,\int_{- \infty }^0\Bigl|
\int_{\CC^{n-1}}z^k\,e^{\,-|\rho| zz^*}\,\psi(\rho ,z)\,d\mu(z)\,\Bigr|^2\,
|\rho |^{2n+2|k|-3}\,d \rho,
\\
\|\psi\|_+^2 &
=\sum_{k\in\ZZ^{n-1}_+}\frac{1}{k!}\,\int_0^{ \infty }\Bigl|
\int_{\CC^{n-1}}\bar z^k\,e^{\,-\rho\,zz^*}\,\psi(\rho ,z)\,d\mu(z)\,\Bigr|^2\,
\rho ^{2n+2|k|-3}\,d \rho.
\end{aligned}
\end{equation}
\end{COR*}

Formulas for the operators $T^+_g$ and $T^-_g$ of the representations of $U(n,1)$
in $\Cal L_+$ and $\Cal L_-$, respectively, are obtained by passing to the limit
as $\lambda \to 0$ from formulas
\eqref{78-24-1}, \eqref{78-25-1}, and \eqref{378-25}
for the corresponding operators in complementary series representations.
Thus we obtain the following theorem.

\begin{THM}{}\label{THM:731-24} In model $C$, the operators of the representations of
$U(n,1)$ in $ \Cal L_{\pm}$ have the following form.
In the space $\Cal L_+$:
\begin{align}{}
\label{78-74}
T^+_{t_0, z _0}\psi(\rho, z ) &
= e^{\,\rho\,(it_0-\frac{|z_0|^2}{2}-zz_0^*)}\,\psi(\rho ,z+z_0 ),
& (t_0, z _0) &\in H,
\\
\label{78-75}
T^+_{\epsilon  , u} \psi(\rho, z ) &
=\psi(| \epsilon|^2\rho,\,\bar\epsilon ^{-1} z u)\,,
& (\epsilon  , u) &\in D,
\end{align}
\begin{align}
\label{78-76}
T_s^+ \psi (\rho,z) &
= \rho ^{-n+1} \int\limits_0^{ \infty }
\int\limits_{- \infty }^{ \infty }
e^{\,-\rho \zeta +\frac{\,\rho '\mathstrut}{\zeta \mathstrut}}\,
|\zeta |^{-2n}\, \rho ^{\prime n-1}\, \psi (\rho ',\,\frac{z}{\zeta })\,dt\,d \rho ',
& \zeta = it &- \frac{|z|^2}{2}.
\end{align}
In the space $\Cal L_-$:
\begin{align}{}\label{78-77}
T^-_{t_0, z _0}\psi(\rho, z ) &
= e^{\,\rho\,(it_0+\frac{|z_0|^2}{2}+z_0z^*)}\,\psi(\rho ,z+z_0 ),
& (t_0, z _0)  &\in H,
\\
\label{78-78}
T^-_{\epsilon  , u} \psi(\rho, z ) &
= \psi(| \epsilon|^2\rho,\,\bar\epsilon ^{-1} z u)\,,
& (\epsilon  , u) & \in D,
\end{align}
\begin{align}
\label{78-79}
T_s^- \psi (\rho,z) &
= |\rho |^{-n+1} \int\limits_{- \infty }^0
\int\limits_{- \infty }^{ \infty }
e^{\,-\rho \OVER\zeta +\frac{\,\rho '\mathstrut}{\OVER\zeta \mathstrut}}\,
|\zeta |^{-2n}\,|\rho' |^{n-1}\psi (\rho ',\,\frac{z}{\zeta })\,dt\,d \rho ',
& \zeta = it & - \frac{|z|^2}{2}.
\end{align}
\end{THM}

\subsection{Relation to the Bargmann models of representations
of the Heisenberg group $H$}
Above, by the Bargmann model we meant the realization of the irreducible representation
of the group $H$ with parameter $\rho$ in the Hilbert space
$\Cal H(\rho )$ of entire analytic (for $\rho >0$) or antianalytic (for
$\rho<0$) functions on $\CC^{n-1}$ with the norm
$$
\| f \|^{2} = |\rho |^{n-1}
\int_{\CC^{n-1}} |f(z)|^{2}\,e^{\, - |\rho| \,|z|^2}\,d\mu(z).
$$
According to Proposition~\ref{PROP:796-45}, for $\rho >0$ 
this norm coincides with the norm
$$
\|f\|^2=\rho^{2n-2}\int\limits_{\CC^{n-1} \times \, \CC^{n-1}} \!\!\!\!\!\!%
f(z)\,\OVER{f(z')}\, e^{\, \rho\,(z'z^*-|z|^2-|z'|^2)}\,d\mu(z)\,d\mu(z');
$$
and for $\rho <0$, it coincides with the same norm in which
$\rho $ is replaced by $|\rho |$ and $z'z^*$ in the exponent
is replaced by $z(z')^*$.

Also recall that the representation of $H$ in the space $\Cal H(\rho  )$
can be extended to a unitary representation of the group ${{P}}_0$.

Comparing the description of $\Cal L(\rho )$ and $\Cal H(\rho )$
and the formulas for the operators of the group ${{P}}_0$ 
in these spaces, we conclude the following.

\begin{PROP}{}\label{PROP:554-1} 
There is a natural isomorphism of Hilbert spaces
$$
\tau:\, \Cal L(\rho )\to\Cal H(\rho )
$$
commuting with the action of the group ${{P}}_0$. Namely,
$\tau\psi\in\Cal H(\rho )$ is determined for $\rho>0$ by the equations
\begin{equation}{}\label{7-2}
\int_{\CC^{n-1}}\bar z^k\,e^{\,-\rho\,zz^*}\,\psi(z)\,d \mu (z)
= \int_{\CC^{n-1}}\bar z^k\,e^{\,-\rho\,zz^*}\,\tau\psi(z)\, d \mu (z)
\quad \text{for every} \quad k\in\ZZ_+^{n-1}.
\end{equation}
In other words,
\begin{equation}{}\label{7-3}
\tau\psi(z)=\sum_{k\in\ZZ_+^{n-1}} a_k\,z^k,
\quad \text{where} \quad
a_k=\frac{\rho ^{|k|+n-1}}{k!}\,
\int_{\CC^{n-1}}\bar z^k\,e^{\,-\rho\,zz^*}\,\psi(z)\,d \mu (z).
\end{equation}
In a similar way we define $\tau$ for $\rho <0$.
\end{PROP}

\subsection{The Bargmann model of the special representations}
Consider the following direct integrals of the Hilbert spaces $\Cal H(\rho )$:
$$
\Cal H_+=\int_0^{+ \infty }\Cal H(\rho )\,\rho ^{-1}\,d \rho ,
\quad
\Cal H_-=\int_{- \infty }^0\Cal H(\rho )\,|\rho |^{-1}\,d \rho.
$$
In other words, $\Cal H_+$ is the space of functions
$f(\rho ,z)=\sum_{k\in\ZZ^{n-1}_+} f_k(\rho )\,z^k$, where $f_k(\rho )$ 
are functions on the half-line $\rho>0$, with the norm
\begin{equation}{}\label{678-92}
\|f\|^2=\sum_{k\in\ZZ^{n-1}_+} k!\,\int_{- \infty }^0
|f_k(\rho )|^2\, |\rho |^{-|k|-1}\,d \rho .
\end{equation}
The space $\Cal H_-$ is defined in a similar way.

The isomorphisms $\tau:\, \Cal L(\rho )\to\Cal H(\rho )$
induce isomorphisms of Hilbert spaces
$$
\tau:\quad   \Cal  L_+\to\Cal H_+,\quad \Cal L_-\to\Cal H_-
$$
commuting with the action of the group ${{P}}$ in these spaces. Namely,
for $\psi\in\Cal L_+$ we have
\begin{equation}{}\label{377-70}
\tau\psi(\rho , z)=\sum_{k\in\ZZ_+^{n-1}} f_k(\rho )\,z^k,
\quad \text{where} \quad
f_k(\rho )=\frac{\rho ^{|k|+n-1}}{k!}\,
\int_{\CC^{n-1}}\bar z^k\,e^{\,-\rho\,zz^*}\,\psi(\rho , z)\,d\mu(z).
\end{equation}
In a similar way we define $\tau$ on $\Cal L_-$.

In view of this isomorphism, the spaces $\Cal H_{\pm}$ are new models
of the special representations of the group $U(n,1)$. Let us call them
the Bargmann models of the special representations.

In the Bargmann model, the action of the operators of the subgroup
${{P}}$ is given by the same formulas as in model $C$. The action of the whole
group $U(n,1)$ is determined by indicating only one operator
$T_s$. Let us study the structure of this operator. For definiteness,  
we restrict ourselves with the case of $\Cal H_+$.

By definition, the operator $T_s$ on $\Cal L_+$ is given by
$$
T_s\psi(\rho ,z)=\tau\psi_1(\rho ,z),
$$
where $\psi_1(\rho ,z)\in\Cal L_+$ stands for the right-hand side of \eqref{78-76}.

Denote by $L_k \subset \Cal H_+$, $k\in\ZZ^{n-1}_+$, the subspace of functions
of the form $f(\rho )\,z^k$. Obviously,
$$
\Cal H_+= \bigoplus _{k\in\ZZ^{n-1}_+} L_k.
$$

\begin{PROP}{}\label{PROP:L_k,T_s}
The subspaces $L_k \subset \Cal H_\pm$ are invariant under the
operator $T_s$.
\end{PROP}

\begin{proof}{}  Let $\psi\in L_k^+$, i.e., $\psi(\rho ,z)=f(\rho )\,z^k$.
Then, according to the formula \eqref{78-76} for the operator $T_s^+$ in the space
$\Cal L_+$, we have  $T_s^+ \psi (\rho,z)=\tau\,\psi_1(\rho ,z)$, where 
\begin{equation}{}\label{249-34}
\psi_1(\rho ,z)= \rho ^{-n+1}\,
\Bigl(\int\limits_0^{ \infty }
\int\limits_{- \infty }^{ \infty }
e^{\,-\rho \zeta +\frac{\,\rho '\mathstrut}{\zeta \mathstrut}}\,
|\zeta |^{-2n}\,\zeta^{-|k|}\,\rho ^{\prime\, n-1} \,f (\rho ')\,dt\,d \rho '\Bigr)\,z^k,
\quad \zeta =it-\frac{|z|^2}{2}.
\end{equation}
Obviously, the function $\psi_1$ can be written in the form
$\psi_1(\rho ,z)=(\sum_0^\infty f_n(\rho )\,|z|^{2n})\,z^k$.
It follows from the definition of $\tau$ that $\tau(f_n(\rho )\,|z|^{2n}\,z^k)\in L_k$
for all $n$. Hence $\psi_1\in L_k$.
\end{proof}

Denote $\rho _+^{\alpha }=\rho ^{\alpha }$ for $\rho >0$;
$\rho _+^{\alpha }=0$ for $\rho <0$.

\begin{PROP}{}\label{PROP:278-32} The action of the operator $T^+_s$ on functions
of the form $\rho_+ ^{\alpha }\,z^k$, $\alpha \in\CC$, regarded as distributions on
$\Cal H_+$, is given by
\begin{equation}{}\label{7-4}
T_s^+(\rho_+ ^{\alpha }\,z^k)=(-1)^{|k|}\,\frac{\Gamma (\alpha )}
{\Gamma (- \alpha +|k|)}\rho_+^{\,-\alpha+|k| }\,z^k.
\end{equation}
\end{PROP}

For a proof of \eqref{7-4}, see Appendix~2.

\begin{PROP}{} \label{PROP:278-33} The action of the operator $T_{s}$
on functions $\xi _k(a; \rho ,z)
= \rho ^{|k|} \,e^{\, - a \,\rho }\, z^{k} \in L_k $,
$\opn{Re}(a)>0$, is given by
\begin{equation}{}\label{278-33}
T_{s} \xi _k(a;\, \rho ,z) = a^{\,- |k|}\,\xi _k(a^{-1};\, \rho ,z).
\end{equation}
\end{PROP}

\begin{proof}{} Let us use the equation
$$
\int _{0}^{\infty }  e^{\, - a \,\rho }\,
\rho ^{\,\frac{|k|}{2} + i \sigma -1 } \,d \rho
= \Gamma (\tfrac{|k|}{2} + i \sigma )\,a^{\,-\tfrac{|k|}{2} - i \sigma}.
$$
By the inversion formula for the Mellin transform,
$$
\rho ^{\,|k|} \,e^{\, - a \,\rho }\,
= \int _{- \infty }^{+ \infty }
        a^{\,-\tfrac{|k|}{2} - i \sigma}
        \,\Gamma (\tfrac{|k|}{2} + i \sigma )
        \,\rho ^{\,\tfrac{|k|}{2} - i \sigma} \,d \sigma .
$$
According to \eqref{7-4}, we have
$$
T_s ( \rho ^{\,\tfrac{|k|}{2} - i \sigma} \,z^{k} )
= \frac{
        \Gamma (\frac{|k|}{2} - i \sigma )
        }{
        \Gamma (\frac{|k|}{2} + i \sigma )
        }
        \,\rho ^{\,\frac{|k|}{2} + i \sigma }.
$$
Therefore
\begin{align*}{}
T_s ( \rho ^{|k|}\,e^{\,- a \rho } \,z^{k} )
&
= \Bigl(
\int _{- \infty }^{+ \infty }
        a^{\,-\tfrac{|k|}{2} - i \sigma}
        \,\Gamma (\tfrac{|k|}{2} - i \sigma )
        \,\rho ^{\,\tfrac{|k|}{2} + i \sigma} \,d \sigma
\Bigr)\,z^{k}
\\  &
= a^{\,-|k|}
\Bigl(
\int _{- \infty }^{+ \infty }
\Bigl( \frac1a \Bigr)^{\,-\tfrac{|k|}{2} - i \sigma}
        \,\Gamma (\tfrac{|k|}{2} + i \sigma )
        \,\rho ^{\,\tfrac{|k|}{2} - i \sigma} \,d \sigma
\Bigr)\,z^{k}
\\ &
= a^{\,- |k|}\,\xi _k(a^{-1};\, \rho ,z).
\end{align*}
The proposition is proved.
\end{proof}

It is convenient to write the action of the operator $T_s$ in $L_k$ 
passing from functions $f_k(\rho ) \,z^k$ to functions $F_k(\sigma  ) \,z^k$, where
\begin{equation}{}\label{895-1}
F_k(\sigma  )=\int_0^\infty f_k(\rho )\,
\rho^{\,-\frac{|k|}{2}-i\,\sigma-1 }
\,d \rho .
\end{equation}

By the inversion formula for the Mellin transform,
\begin{equation}{}\label{895-2}
f_k(\rho ) = \rho^{\frac{|k|}{2}}\,
\int_{- \infty }^\infty F_k(\sigma )\, \rho ^{i \,\sigma }\,d \sigma.
\end{equation}

Let us introduce the Hilbert space $\TIL{\Cal H} _+$
of functions $F(\sigma ,z)$ entire in $z\in\CC^{n-1}$ with the norm
\begin{equation}{}\label{895-3}
\|F\|_+^2 =\sum_{k\in\ZZ^{n-1}_+}k!\,\int_{- \infty }^\infty
|F_k(\sigma)|^2\,d \sigma\quad \text{for}\quad
F(\sigma,z)=\sum_{k\in\ZZ_+^{n-1}}F_k(\sigma)\,z^k.
\end{equation}

\begin{PROP}{}\label{PROP:768-45}
The map
$$
\psi(\rho ,z)=\sum_{k\in\ZZ_+^{n-1}}f_k(\rho )\,z^k\to
F(\lambda ,z)=\sum_{k\in\ZZ_+^{n-1}}F_k(\lambda  )\,z^k,
$$
where $f_k(\rho )$ and $F_k(\lambda )$ are related by \eqref{895-1} and
\eqref{895-2}, is an isomorphism of Hilbert spaces
$\Cal H_+\to\TIL{\Cal H}_+$.
\end{PROP}

The assertion follows immediately from the Plancherel formula
for the Mellin transform.

Denote by $\TIL T^+_s$ the image of the operator $T^+_s$ on $\Cal H_+$
under the map $\Cal H_+\to\TIL{\Cal H}_+$.

\begin{PROP}{}\label{PROP:891-6} The action of the operator $\TIL T^+_s$ on functions
$F_k(\lambda )\,z^k\in\TIL{\Cal H}_+$ is given by
\begin{equation}\label{891-6}
\TIL T^+_s(F_k(\lambda )\,z^k)=
(-1)^{|k|}\,
\frac{\Gamma (-i \lambda +\frac{|k|}{2})}
{\Gamma (i \lambda +\frac{|k|}{2})}\,F_k(- \lambda )\,z^k.
\end{equation}
\end{PROP}

\begin{proof}{} According to Proposition~\ref{PROP:278-32}, we have
$$
T_s^+ (\rho ^{\frac{|k|}{2} + i \sigma}\,z^k)=
(-1)^{|k|}\,
        \frac{
        \Gamma (  i \sigma+\frac{|k|}{2}   )
        }{
        \Gamma ( -i \lambda +\frac{|k|}{2} )
        }\,
\rho^{\frac{|k|}{2}-i \sigma}\,z^k.
$$
Hence \eqref{895-2} implies
\begin{align*}{}
T^+_s(f_k(\rho )\,z^k)&=
(-1)^{|k|}\,
\Bigl(
\int_{- \infty }^\infty F_k(\sigma)\,
\frac{\Gamma (i \sigma+\frac{|k|}{2})}{\Gamma (-i \sigma+\frac{|k|}{2})}\,
\rho^{\frac{|k|}{2}-i \sigma}\,d \sigma
\Bigr)
\, z ^k
\\
&=(-1)^{|k|}\,
\Bigl(
\int_{- \infty }^\infty F_k(-\sigma)\,
\frac{\Gamma (-i \sigma+\frac{|k|}{2})}{\Gamma (i \sigma+\frac{|k|}{2})}\,
\rho^{\frac{|k|}{2}+i \sigma}\,d \sigma
\Bigr)
\, z ^k.
\end{align*}
On the other hand,
$$
T^+_s(f_k(\rho )\,z^k)=
\int_{- \infty }^\infty \TIL T^+_s(F_k(\sigma)\,z^k)\,
 \rho ^{ \frac{|k|}{2} +i \,\sigma} \, d \sigma.
$$
Comparing these equations, we obtain (\ref{891-6}).
\end{proof}

\subsection{Description of the nontrivial $1$-cocycles
$a^{\pm}:  U(n,1)\to\Cal H_{\pm}$}
The nontrivial $1$-cocycle $a:\,U(n,1)\to\Cal H$, where 
$\Cal H=\Cal H_+ \oplus \Cal H_-$, is given by the formula
$$
a(g)=T_g \xi - \xi ,
$$
where $\xi =f(\rho ,z)$ is invariant under the action
of the maximal compact subgroup $U \subset U(n,1)$ and $\xi \notin L$.
The $1$-cocycles $a^{\pm}(g)$ are the projections of $a(g)$ 
to the subspaces $\Cal H_{\pm}$.

By definition, $a(u)=0$ for $u\in U$, whence
$$
a(ug)=a(g)\quad \text{for} \quad u\in U.
$$
Since an element $g\in U(n,1)$ can be written as the product $g=ub$, where
$u\in U$, $b\in {{P}}_1$, the $1$-cocycle $a(g)$ is uniquely determined by
its restriction to the subgroup ${{P}}_1$.

\begin{THM}{}\label{THM:} In the Bargmann model, the function
$\xi $ in the formula for the $1$-cocycle has the form
\begin{equation}{}\label{758-45}
\xi (\rho ,z)=e^{-|\rho |}.
\end{equation}
\end{THM}

\begin{proof}{}
We use the fact that in model $B$ the function $\xi $ has the form
$$
\phi(\zeta ,z)=|1- \zeta |^{-2n},  \quad \text{where} \quad
\zeta =it-\frac{|z|^2}{2}.
$$
It follows that in model $C$ it has the form
$$
\psi(\rho ,z)= |\rho|^{-n+1}\,\,e^{|\rho |\frac{|z|^2}{2}}\,
\int_{- \infty }^{+ \infty }|1- \zeta |^{-2n}\,e^{-i \rho t }\,dt.
$$
Note that $\psi$ is an even function in $\rho $; thus below we assume that
$\rho >0$. The required vector $\xi $ in the Bargmann model for
$\rho >0$ has the form
$$
\xi (\rho ,z)=\tau\psi(\rho ,z)=\sum_{k\in\ZZ^{n-1}_+} f_k(\rho )\,z^k,
$$
where the functions $f_k(\rho )$ are given by \eqref{377-70}.
In view of this equation, it is clear that $f_k(\rho )=0$ for $k>0$. Hence
\begin{equation}{}\label{743-86}
\xi (\rho ,z)=f_0(\rho )=\int_{- \infty }^{+ \infty }
\int_{\CC^{n-1}} e^{\rho \zeta }\,|1- \zeta |^{-2n}\,
d\mu(z)\,dt.
\end{equation}
Let us calculate the integral \eqref{743-86}. We have
$$
\Bigl(1-\frac{ d }{ dx }\Bigr)^{n} f_0(\rho )
=\int_{- \infty }^{+ \infty }
\int_{\CC^{n-1}} e^{\rho \zeta }\,(1- \OVER\zeta )^{-n}\,d\mu(z)\,dt.
$$
Obviously,
$$
\gathered
\int_{- \infty }^{+ \infty }
        \frac{  e^{\rho \zeta }  }{  (1- \OVER\zeta )^{n}  }\,dt
=\int_{- \infty }^{+ \infty }
        \frac{  e^{\rho \zeta }  }{  (1+ |z|^2 + \zeta  )^{n}  }\,dt
= \opl{Res}_{\zeta = -1- |z|^2} \Bigl(
        \frac{  e^{\rho \zeta }  }{  (1+ |z|^2 + \zeta  )^{n}  }
\Bigr)
\\  \\
= const\,\rho ^{n-1} \,e^{\,- \rho\,(1+|z|^2)}.
\endgathered
$$
Hence
$$
\Bigl(1-\frac{ d }{ dx }\Bigr)^{n} f_0(\rho )
= const\,\rho ^{n-1} \,e^{\,- \rho}
\int_{\CC^{n-1}} e^{\rho\, |z|^2}\,d\mu(z)
= const\,e^{\,- \rho}.
$$
Since $f_0(\infty ) = 0$, it follows that $f_0 (\rho ) =  const\,e^{\,- \rho}$.
\end{proof}

\part*{Appendices}

\section{Appendix 1. Canonical states and canonical representations}\label{sect:7}

Canonical states on groups and the corresponding canonical
representations were first introduced for the group
$SL(2,\RR)$ in \cite{V-G-1}, and then for the groups
$O(n,1)$ and $U(n,1)$ in \cite{V-G-6}; they were used for constructing irreducible unitary 
representations of the current groups $O(n,1)^X$ and $U(n,1)^X$. 
These representations are unitary, reducible, and, similarly to 
complementary series representations, determined by one real parameter
$\lambda\in (0, \lambda _0)$, $\lambda _0>0$.

The following property shows the relation between the canonical representations 
of the groups $O(n,1)$ and $U(n,1)$ and
the complementary series representations of these groups:
\par\hangindent=\parindent

{\it Each complementary series representation $T^{\lambda }$ is contained
as a direct summand with multiplicity one in exactly one 
canonical representation, the one with the same parameter $\lambda$.}

There is a simple relation between the canonical representations of the groups
$O(n,1)^X$ and $U(n,1)^X$ and their special representations. Namely,
for $\lambda \to 0$ the norm $\|f\|_{\lambda }$ in the space $H^{\lambda }$
of the canonical representation degenerates on a space $H \subset H^0$ 
of codimension $1$. Therefore, as $\lambda \to 0$, the canonical representations
$T^{\lambda }$ tend to the identity representation on
the quotient space $H^0/H$. The space $H$, glued to this one-dimensional subspace,
has a positive invariant norm $\|f\|$. In the case of the groups
$O(n,1)$, $n>2$, this norm is nondegenerate; the completion
$\Cal H$ of the space $H$ with respect to this norm contains the space of the
special representation as a direct summand. In the case of the groups
$U(n,1)$, this norm degenerates on some infinite-dimensional subspace
$H_0 \subset H$; the completion $\Cal H$ of the space
$H/H_0$ with respect to this norm is the direct sum of the two special subspaces.

The role of canonical representations in the representation theory of
current groups is that their spherical functions  $\Psi _{\lambda } (g)$
satisfy the infinite divisibility condition:
$$
\Psi _{\lambda_1 + \lambda _2 } (g)  = \Psi _{\lambda _1} (g) \,\Psi _{\lambda _2} (g).
$$

In Secs.~1--3 of  this Appendix we give a general definition
of a canonical state and the associated family of canonical representations
$T^{\lambda }$ for an arbitrary locally compact group $G$ and discuss their
basic properties. The subsequent sections are devoted to the canonical representations 
of the groups $O(n,1)$ and $U(n,1)$.

For the classical simple Lie groups, other approaches to the notion
of canonical representation, various generalizations of this notion,
and the related harmonic analysis are presented in a series of papers.
See, e.g., \cite{GvD-H,Mol} and the bibliography in \cite{Mol}.

\subsection{General definition of canonical states and canonical representations}
Let $G$ be a locally compact, noncompact topological group, 
let $K$ be a subgroup of $G$, and consider the homogeneous space $X=K \backslash G $ of 
the group $G$, on which $G$ acts by right shifts: $x\to xg$. We assume that on $X$
there exists a $G$-invariant measure $d\mu(x)$.

Denote by $L(X)$ the subspace of compactly supported continuous functions on $X$.
We will define elements from $L(X)$ by their liftings to the group $G$, i.e.,
as functions $f(g)$ satisfying the relation
$f(kg) = f(g)$ for $k \in K$, $g \in G$.

\setcounter{DEF}{0}

\begin{DEF}{}\label{DEF:1} We say that a continuous positive function
$\phi (g)$ on $G$ is positive definite on $X$ if
\begin{enumerate}
\item $\phi (g^{-1}) =  \phi (g)$,
\item $\phi (k_1 g k_2) = \phi (g)$ for $k_1,k_2 \in K$,
\item for every nonzero function $f \in L(X)$,
\begin{equation}{}\label{7-7-7}
\| f \|^2 \sim \int_{ X \times X } \phi (g_1g_2^{-1})\,f(g_1)\,\OVER{f(g_2)}
\,d \mu (x_1)\,d \mu (x_2) > 0.
\end{equation}
\end{enumerate}
\end{DEF}

The expression in the right-hand side of \eqref{7-7-7} is well defined,
because the integrand function can be projected 
from $G \times G$ to $X \times X$.

\begin{REM*}{} By property 2), the function $\phi (g)$ is projected to $X$
by the map $G \to K \backslash G$.
\end{REM*}

\begin{DEF}{}\label{DEF:2} A canonical state on the group $G$ with respect to
the subgroup $K$ is a family of positive continuous functions
$\Psi ^{\,- \lambda }(g)= (\Psi (g))^{\,- \lambda }$ on $G$ depending
on a real positive parameter $\lambda $ (degree) and satisfying the conditions
\begin{enumerate}
\item $\Psi ^{\,- \lambda }(g) $ is positive definite on $X$ in the sense of
Definition~1 for small values of $\lambda >0$,
\item for the function $\left. \frac{ d \, \Psi ^{\,- \lambda }(g)  }{ d \lambda }
\right|_{\lambda =0} = - \log \Psi (g)$, the following integral diverges:
\begin{equation}{}\label{7-7-7-8}
- \int_{ X  } \log \Psi (g)\,d \mu (x).
\end{equation}
\end{enumerate}
\end{DEF}

\begin{PROP}{}\label{PROP:7-7-7} Let $B(x_1,x_2)$ be a continuous positive function
on $X \times X$ satisfying the following conditions $1)$--$4)$:
\begin{enumerate}
\item symmetry: $B(x_1,x_2) = B(x_2,x_1)$,
\item $G$-invariance:
$$
B(x_1g,x_2g)=B(x_1,x_2) \quad \text{for every}\quad g \in G,
$$
\item for every nonzero function $f \in L(X)$, for small $\lambda >0$,
\begin{equation}{}\label{3513}
\int_{X \times X}
B^{- \lambda }(x_1,x_2)\,f(x_1) \,\OVER{f(x_2)}\,d\mu(x_1)\,d\mu(x_2) >0,
\end{equation}
\item
the following integral diverges:
$$
-\int _X \log B(x_1,x_2) \,d\mu (x_2)
$$
(this integral does not depend on $x_1$ in view of condition $2$)).
\end{enumerate}
Then $\Psi ^{\, - \lambda }(g) = B^{\, - \lambda }(x_0g,x_0)$, where
$x_0 \in X$ is a fixed point of the subgroup $K$,
is a canonical state.
\end{PROP}

\begin{REM*}{} In the case when $G$ is the group $O(n,1)$ or $U(n,1)$ and
$K$ is its maximal compact subgroup, there is a simple relation between 
the kernel $B^{- \lambda }(x_1,x_2)$ and the Berezin kernel, see Section~4 below.
\end{REM*}

\begin{DEF}{}\label{DEF:3} 
The family of Hilbert spaces associated with a canonical state
$\Psi ^{\,- \lambda }(g)$ is
the family of Hilbert spaces
$H^{\, \lambda }$ of functions on $X$ with the norms
\begin{equation}{}\label{7-7-21}
\| f \|_{ \lambda }^2
= \int_{ X \times X } \Psi^{\, - \lambda } (g_1g_2^{-1})\,f(g_1)\,\OVER{f(g_2)}
\,d \mu (x_1)\,d \mu (x_2),
\end{equation}
or, in terms of the kernels $B^{ - \lambda }(x_1,x_2)$,
\begin{equation}{}\label{7-7-231}
\| f \|_{ \lambda }^2
= \int_{ X \times X } B^{\, - \lambda } (x_1,x_2) \,f(x_1)\,\OVER{f(x_2)}
\,d \mu (x_1)\,d \mu (x_2).
\end{equation}
\end{DEF}

\begin{DEF}{}\label{DEF:4} The representations of the group $G$
associated with a canonical state $\Psi ^{\,- \lambda }(g)$,
or simply the canonical representations, are the unitary representations of $G$
in the spaces $H^{\lambda }$ by operators of the form
$$
(T_g f)(x) = f(xg).
$$
\end{DEF}

\begin{REM*}{} If we do not require that 
$\Psi ^{\,- \lambda }(g)$ is positive definite for all $\lambda$,
then this representation can be nonunitary.
\end{REM*}

\subsection{Spherical functions}

Note that the norm \eqref{7-7-21} is defined for any 
compactly supported distribution on $X$, and elements of the space
$H^{\lambda }$ can be interpreted as such distributions. In particular,
the delta function $\delta _{x_0}(x)$ at an arbitrary point
$x_0\in X$, defined by the formula 
$$
\int_X \delta _{x_0}(x)\,f(x_0)\,d\mu(x)=f(x_0),
$$
is an element of $H^{\lambda }$ for every $\lambda\in (0, \lambda _0)$.

\begin{PROP}{}\label{PROP:71-80} The functions
$\delta _{x_0}(x)\in H^{\lambda }$ have the following properties:
\begin{enumerate}
\item any finite number of them are linearly independent, 
and they form an overdefined system in
$H^{\lambda }$,
\item $\|\delta _{x_0}(x)\|_{\lambda }=1$,
\item for every $x_0$, the vector $\delta _{x_0}(x)$ is cyclic in
$H^{\lambda }$,
\item the vector $\delta _{x_0}(x)$ is invariant under the stationary
subgroup of the point $x_0$.
\end{enumerate}
\end{PROP}

\begin{DEF*}{}  Let us fix the special cyclic vector
$\xi_{\lambda }=\delta _{x_0}(x)$, where $x_0\in X$ is a point with
stationary subgroup $K$, and define the spherical function
$\psi _{\lambda }(g)$ on  $G$ associated with the space
$H^{\lambda }$ by the formula
$$
\psi _{\lambda }(g) = {\langle T_{g}\xi _{\lambda },\, \xi _{\lambda } \rangle}_{\lambda }\,.
$$
\end{DEF*}

\begin{PROP}{}\label{PROP:7:5} The spherical function can be expressed
in terms of the canonical state by the formula
$$
\psi _{\lambda }(g) = \Psi^{- \lambda }(g);
$$
thus it is an infinitely divisible function, i.e.,
$$
\psi _{\lambda_1 + \lambda _2 } (g)  = \psi _{\lambda _1} (g) \,\psi _{\lambda _2} (g).
$$
\end{PROP}

Indeed, setting $x_1=x_0g_1$,  $x_2=x_0g_2$, we have
\begin{gather*}{}
{\langle T_{g}\xi _{\lambda },\, \xi _{\lambda } \rangle}_{\lambda }
=\int_{{X} \times {X}} \Psi^{- \lambda } g_1g_2^{})\,
\delta _{x_0}(x_1g)\,\delta _{x_0}(x_2)\,d \mu (x_1)\, d \mu (x_2)
\\
=\int_X \Psi^{- \lambda }(g_1)\, \delta _{x_0}(x_1 g)\,
d \mu (x_1)
=\int_X \Psi^{- \lambda }(g_1 g^{-1})\, \delta _{x_0}(x_1 )\,d \mu (x_1)=
\Psi^{- \lambda }( g^{-1}).
\end{gather*}

The infinite divisibility of the spherical function implies the following proposition.

\begin{PROP}{}\label{PROP:4356} For $\lambda = \lambda _1+ \ldots + \lambda _n$
there exists an isometric embedding
$$
H ^{\lambda }\to H ^{\lambda_1 } \otimes \ldots \otimes  H ^{\lambda _n}
$$
commuting with the action of the group $G$. Namely, the map
$$
\xi_{\lambda }\to \xi_{\lambda _1}\otimes \ldots \otimes \xi_{\lambda _n}
$$
can be extended to a required morphism.
\end{PROP}

Indeed, the assertion follows from the equation
$$
{\langle T^{\lambda}_{g} \xi_{\lambda }, \xi_{\lambda} \rangle}
= \prod_{i=1}^n
{\langle T^{\lambda_i}_{g} \xi_{\lambda _i},\xi_{\lambda_i} \rangle},
$$
valid for every $g\in G$ with $\lambda = \lambda _1+ \ldots + \lambda _n$.

\subsection{The special representations of the group $G$ associated with
the canonical state $\Psi ^{\,- \lambda }(g)$}

As $\lambda \to 0$, the canonical state tends to one, hence the norm
$\|f\|^2_0=\lim_{l\to 0}\|f\|^2_{\lambda }$ can be written as
$\|f\|^2_0=|\int_X f(x)\,d\mu(x)|^2$;
the operators of the representation are still given by the formula
$$
(T_g f)(x)=f(xg).
$$
Let us study the structure of the subspace on which this norm
degenerates, i.e, the invariant subspace $H$ of functions satisfying the condition
$$
\int_X f(x)\,d\mu(x)=0.
$$

In the space $H$ there is a pre-norm invariant under the operators
of the group $G$:
$$
\|f\|^2= -\int_{ X \times X } \log\Psi (g_1g_2^{-1})\,
f(g_1)\,\OVER{f(g_2)}\,d \mu (x_1)\,d \mu (x_2).
$$

Let $H_0 \subset  H$ be the subspace on which this pre-norm degenerates.
Then we can define a unitary representation
$T_g$ of $G$ in the Hilbert space $\Cal H$ obtained by completing the quotient
$H/H_0$ with respect to the norm $\|f\|$.

\begin{PROP}{}\label{PROP:132-43} The representation $T_g$ of the group $G$ in
the space $\Cal H$ has a nontrivial $1$-cocycle
$a:\,G\to\Cal H$, namely,
$$
a(g)=T_g\xi-\xi,  \quad \text{where} \quad  \xi= \delta _{x_0}(x)
\in H \setminus H_0.
$$
\end{PROP}

\begin{proof}{} First of all, note that $\xi\notin H$ and $a(g)\in H$, i.e., 
$a(g)$ is a $1$-cocycle $G\to\Cal H$. Assume that this cocycle is trivial, i.e.,
$T_g \xi - \xi = T_g \xi_0 - \xi_0$ for some
$\xi _0 \in H$. Then the vector $\xi _1 = \xi - \xi _0$ is
$G$-invariant; hence $\xi _1 = const$. On the other hand, the integral
$$
\int_{ X } \log \Psi (g)\,f(g)\,d \mu (x) 
$$
converges for $f=\xi_1$, since it converges for $f=\xi$ and for $f=\xi_0$.
Hence, since $\xi_1\ne 0$, the integral
$$
\int_{ X  } \log \Psi (g)\,d \mu (x) 
$$
also converges, which  contradicts \eqref{7-7-7-8}. 
\end{proof}

\begin{COR*}{} The space $\Cal H$ contains the special subspaces as direct summands.
\end{COR*}

\begin{PROP}{}\label{PROP:856-7} The norm of the nontrivial $1$-cocycle 
$a:\,G\to\Cal H$ can be expressed in terms of the canonical state by the formula
${\|a(g)\|^2 = 2\log\Psi(g)}$.
\end{PROP}

\begin{proof}{} Indeed, we have
\begin{gather*}{}
\|a(g)\|^2
={\langle T_g\xi-\xi,T_g\xi-\xi \rangle}
= \left. \frac{d}{d \lambda } \right|_{\lambda =0}\!\!%
{\langle T^{\lambda }_g\xi_{\lambda }-\xi_{\lambda },
T^{\lambda }_g\xi_{\lambda }-\xi_{\lambda } \rangle}_{\lambda }
\\
= \left. 2\frac{d}{d \lambda } \right|_{\lambda =0}\!\!%
(1-{\langle T^{\lambda }_g\xi_{\lambda },\xi_{\lambda } \rangle})_{\lambda }
= \left. -2\frac{d\Psi^{\,-\lambda}  (g)}{d \lambda }\right|_{\lambda =0}
= 2\log\Psi(g).
\end{gather*}
\end{proof}

\subsection{Canonical representations of the groups $O(n,1)$ and $U(n,1)$}

Let us realize $O(n,1)$ and $U(n,1)$ as the groups of linear transformations
in $\RR^{n+1}$ and $\CC^{n+1}$ preserving, respectively, the bilinear form
$x_1x'_1+ \ldots +x_{n}x'_{n}-x_{n+1}x'_{n+1}$ and the Hermitian form
$x_1\OVER{x'_1}+ \ldots +x_{n}\OVER{x'_{n}}-x_{n+1}\OVER{x'_{n+1}}$.
In this model, elements of the groups are block matrices
$g= \begin{pmatrix} \alpha & \beta \\ \gamma & \delta \end{pmatrix}$,
where $ \alpha $ and $ \delta$ are square matrices of orders
$n$ and $1$, and the subgroup $K$ of block-diagonal matrices
$k= \begin{pmatrix} \alpha & 0\\ 0& \delta \end{pmatrix}$,
where $\alpha \in U(n)$ and $\delta \in U(1)$,
is their maximal compact subgroup.

Consider the homogeneous space $K \backslash  U(n,1)$. In the case of $O(n,1)$,
it is realized as the open unit ball $D \subset \RR^n$, and in the case of
$U(n,1)$, as the open unitary unit ball $D \subset \CC^n$, i.e.,
$$
D: \quad |z|^2\equiv |z_1|^2 + \ldots + |z_n|^2 <1,
$$
where $z_i\in\RR$ in the case of $O(n,1)$ and $z_i\in\CC$ in the case of $U(n,1)$.

In both cases the action of the group is given by the formula 
$$
z\to zg = (z\,\beta + \delta )^{-1}\,(z\,\alpha + \gamma )
\quad \text{for
$g= \begin{pmatrix} \alpha & \beta \\ \gamma & \delta \end{pmatrix}$
}.
$$

\begin{PROP}{}\label{PROP:23-45} The functions
\begin{equation}{}\label{27-95}
B({z},{z'})=\frac{|1-{z}{z'}^*|}
{(1-{z}{z}^*)^{1/2}\,(1-{z'}{z'}^*)^{1/2}}
\end{equation}
on the real and complex unitary ball, respectively, satisfy the conditions
of Proposition~\ref{PROP:7-7-7}.
\end{PROP}

\begin{COR*}{} The function
$$
\Psi ^{\, -\lambda }(g) = B^{\, -\lambda }(z_0g,z_0),
\quad \text{where}\quad z_0=0,
$$
is a canonical state.
\end{COR*}

The explicit expression for $B(z,z')$ implies the following proposition.

\begin{PROP}{}\label{PROP:22389} The canonical state is given by the formula
$$
\Psi ^{\,-\lambda }(g) = |\delta |^{- \lambda }
\quad \text{for $g= \begin{pmatrix} \alpha & \beta \\ \gamma & \delta \end{pmatrix}$}.
$$
\end{PROP}

\begin{COR*}{}
In matrix model $b$ of the groups $O(n,1)$ and $U(n,1)$,
the canonical states are given by the formula
\begin{equation}{}\label{2-3-4-5-6-2}
\Psi^{\, - \lambda } (g)
= \Bigl|\frac{g_{11} + g_{33} - g_{13} - g_{31}}{2} \Bigr| ^{- \lambda}.
\end{equation}
In particular,
$$
\Psi ^{\,- \lambda }  (\gamma )
= \left( 1 + \frac{ |  \gamma |^2 }{ 4 } \right)^{- \lambda }
\quad \text{and}\quad
\Psi ^{\,- \lambda }  (t,\gamma )
= \left( \frac{ t^2 }{ 4 }  +
\left(  1 +  \frac{ | \gamma | ^2 }{ 4 } \right)^2\, \right)^{-  \lambda /2}
$$
for elements $\gamma \in Z$ and $(t,\gamma ) \in H$
of the corresponding unipotent subgroups.
\end{COR*}

The canonical representations of the groups $O(n,1)$ and $U(n,1)$
associated with this canonical state act in the Hilbert spaces
$H^{\lambda }$ of functions on the real unit ball
$D \subset \RR^n$ and the complex unit unitary ball
$D \subset \CC^n$, respectively, with the norms
\begin{equation}{}\label{35-48}
\|f\|^2_{\lambda }=\int_{D \times D} \bigl|\frac{1-{z}{z'}^*}
{(1-{z}{z}^*)^{1/2}\,(1-{z'}{z'}^*)^{1/2}}\Bigr|^{- \lambda }\,
f(z_1)\OVER{f(z_2)}\,d\mu(z_1) \,d\mu(z_1).
\end{equation}
Here $d\mu(z)$ is the invariant measure on $D$ defined by 
$$
d\mu(z)=(1-zz^*)^{-k(n+1)}\,dz,
$$
where $dz$ is the standard Lebesgue measure, $k=\frac12$ in the case of 
$O(n,1)$, and $k=1$ in the case of $U(n,1)$. The operators of these
representations have the form
$$
(T_g f)(z)=f(zg).
$$

Passing from functions $f(z)$ to the functions
$\phi(z) = |1-zz^*|^{\, -\frac\lambda2 - k(n+1)}\,f(z)$,
we can realize these representations in the Hilbert spaces with the norms
\begin{equation}{}\label{35-45}
\|\phi \|^2_{\lambda }=\int_{D \times D} |1-{z}{z'}^*|^{- \lambda }\,
\phi (z_1)\OVER{\phi (z_2)}\,dz_1\,dz_2.
\end{equation}
In the new realization, the operators of the representation  are given by
$$
(T^{\lambda }_g \phi )(z)-\phi (zg)\,|z \beta + \delta |^{\lambda -2k(n+1)} \quad
\text{for} \quad
g= \begin{pmatrix} \alpha & \beta \\ \gamma & \delta \end{pmatrix}.
$$

\begin{REM*}{}  The definition of canonical states and canonical representations
can be extended to all groups $O(p,q)$, $U(p,q)$, $p\ge q\ge 1$. Elements
of these groups are written as block matrices
$g= \begin{pmatrix} \alpha & \beta \\ \gamma & \delta \end{pmatrix}$,
where $ \alpha $ and $ \delta$ are square matrices of orders $p$ and $q$.
The analogs of the unit balls are the unit matrix balls, that is,
the manifolds of real and complex  $q \times p$ matrices
$z$, respectively,  satisfying the condition $zz ^*< e_q$, where $e_q$ is the unit $q$-matrix.

In the general case, the canonical state $\Psi^{- \lambda }(g)$ and the kernel $B(z,z')$
are given by the formulas
$$
\Psi ^{\,-\lambda }(g) = |\det\delta |^{- \lambda }\quad
\text{for $g= \begin{pmatrix} \alpha & \beta \\ \gamma & \delta \end{pmatrix}$},
$$
$$
B({z},{z'})=\frac{|\det(e_q-{z}{z'}^*)|}
{\det^{1/2}(e_q-{z}{z}^*)  \,\det^{1/2}(e_q-{z'}{z'}^*)}.
$$
The case $q>1$ differs from the one considered here in that the canonical state
is positive definite only for $\lambda > q-1$.

Note that the functions $\det^{\,\lambda } (e_q - {z}{z'}^*)$ are the 
Berezin kernels. In the case of $U(n,1)$, instead of the Berezin kernels we
consider their absolute values.
\end{REM*}

\subsection{Relation to the complementary series representations}

For each of the groups $O(n,1)$ and $U(n,1)$ there exists an isometric embedding
$$
L^{\lambda }\to H^{\lambda }
$$
commuting with the action of the group, where $L^{\lambda }$  is the space
of an arbitrary complementary series representation and $H^{\lambda }$
is the space of the canonical representation with the same parameter
$\lambda$. Let us describe this embedding explicitly.

We will use the realization of complementary series representations 
in the space of functions $f(\omega )$ on the unit sphere $ S \subset \RR^n$
in the case of $O(n,1)$ and on the unitary unit sphere
$S \subset \CC^n$ in the case of $U(n,1)$.

\begin{PROP}{}\label{PROP:7:6} If $H^{\lambda }$ is realized as the space with 
the norm \eqref{35-45}, then for each of the groups $O(n,1)$ and $U(n,1)$
the map
$$
\tau : f(\omega ) \to \phi (z) = \delta (1- |z|^2) \,f(\tfrac{z}{|z|})
$$
(where $\delta (t)$ is the delta function on $\RR$) is an isometric embedding
of $L^{\lambda }$ into $H^{\lambda }$ commuting with the action
of the corresponding group.
\end{PROP}

\begin{proof}{} For definiteness, consider the case of the group $U(n,1)$.
The complementary series representation $T^{\lambda }$ of this group
acts in the Hilbert space $L^{\lambda }$ of functions on the unitary sphere
$S = S^{2n-1}$ in $\CC^n$ with the norm
\begin{equation*}{}
\|f\|^2=\int_{S \times S}
|1-  \omega \omega^{\prime *}|^{- \lambda }\,
f(\omega )\,\OVER{f(\omega ')}\,d \omega \,d \omega ',
\end{equation*}
where $d \omega $ is the invariant measure on $S$. The 
operators of this representation have the form
\begin{equation*}{}
T^{\lambda }_g f(\omega )=f(\omega g)\,|\omega \beta + \delta |^{\lambda - 2n}.
\end{equation*}

Let us prove that  $\tau$ commutes with the action of $U(n,1)$. We have
$$
(T_g \tau f)(z) = \delta (1 - |z{g}|^2)\, f(\tfrac{zg}{|z{g}|})
\,|z\,\beta + \delta |^{\lambda - 2(n+1)}.
$$
Since $\delta (1- |z{g}|^2) = \delta (1-|z|^2) \,|z \beta + \delta |^2$,
it follows that
$$
(T_g \tau f)(z) = \delta (1-|z|^2) \, f(\tfrac{zg}{|z{g}|})
\,|z\,\beta + \delta |^{\lambda - 2n}
= (\tau T_g f)(z).
$$
Let us prove that the embedding $L^{\lambda }\to H^{\lambda }$ is isometric.
Let $\phi (z) = (\tau f)(z) = \delta (1-|z|^2)\,f(\tfrac{z}{|z|}) $.
Passing to the spherical coordinates on $D=K \backslash  U(n,1)$,
i.e., the coordinates $z = r \omega $, $r>0$, $\omega \in S$, we obtain
\begin{align*}{}
\| T_g \phi \|^2_{\lambda }
&
= \int |1 - rr' \omega \omega ^{\prime *}|^{ - \lambda }
\, \delta (1 - r^2)\, \delta (1 - r^{\prime 2})\,(rr')^{2n-1}\,f(\omega )\,\OVER{f(\omega ')}
\, dr \,dr'\,d \omega \,d \omega '
\\ &
= \int _{S \times S } |1 - \omega \omega ^{\prime *}|^{ - \lambda }
\,f(\omega )\,\OVER{f(\omega ')} \,d \omega \,d \omega '
= \|f\|^2.
\end{align*}
The proposition is proved.
\end{proof}

\subsection{The structure of the limiting space $H^0$; 
the $1$-cocycles and special representations associated with $H^0$}

According to general definitions, the space $H^0$, which is the $\lambda \to 0$ limit
of the spaces $H^{\lambda }$ of the canonical representations,
contains an invariant subspace $H \subset H^0$ 
of codimension $1$, namely, the space of functions satisfying the relation
$$
\int_D f(z)\,d\mu(z)=0;
$$
the quotient representation in $H^0/H$ is the identity representation.
Thus, as
$\lambda \to0$, the representations $T^{\lambda }$ in the spaces $H^{\lambda }$
tend to the identity representation in the quotient space $H^0/H$.

Let us study the structure of the subspace $H$ glued to $H^0/H$.
On this subspace there is a seminorm $\|f\|$ invariant under
the operators $(T_g f)(z)=f(zg)$:
$$
\|f\|^2=-\int_{D  \times D}\,\log B(z,z')\,
f(z) \,\OVER{f(z')}\,d\mu(z)\,d\mu(z').
$$
The condition $f\in H$ implies that this seminorm can be written in the form
\begin{equation}{}\label{`8-35}
\|f\|^2=-\int_{D  \times D}\,\log |1-zz'^*|\,
f(z) \,\OVER{f(z')}\,d\mu(z)\,d\mu(z').
\end{equation}
Passing from functions $f(z)$ to the functions $\phi(z)=|1-zz^*|^{-k(n+1)}$, where
$k=1/2$ and $k=1$ for the groups $O(n,1)$ and $U(n,1)$, respectively,
we obtain the following expression for the norm:
\begin{equation}{}\label{8-38}
\|\phi \|^2=-\int_{D  \times D}\,\log |1-zz'^*|\,
\phi(z) \,\OVER{\phi(z')}\,dz\,dz',
\end{equation}
where $dz$ is the Lebesgue measure on $D$. In the new realization, 
the operators of the representation are given by the formula
$$
(T_g\phi)(z) = \phi(zg)\,|z \beta + \delta |^{ -2 k(n+1)}
\quad \text{
for $g= \begin{pmatrix} \alpha & \beta \\ \gamma & \delta \end{pmatrix}$}.
$$

In the case of $O(n,1)$, the norm \eqref{8-38} is nondegenerate; in the case
of $U(n,1)$, it degenerates on some infinite-dimensional subspace
$H_0$. Consider these cases separately.

\subsubsection{The case of the group $O(n,1)$}

Expanding $\log |1-zz'^*|=\log (1-zz'^*)$ into a power series, we obtain the
following expression for the norm $\|\phi\|$:
$$
\|\phi\|^2=\sum_{|k|>0}\frac{(|k|-1)!}{k!}\,\Bigl|\int z^k\,
\phi(z)\,dz\Bigr|^2,
$$
where, as usual, $k!= \prod k_i!$, $|k|= \sum k_i$, and $z^k= \prod z_i^{k_i}$.
It follows that this norm is strictly positive. Denote by $\Cal H$
the Hilbert space obtained by completing $H$ with respect to  the norm $\|\phi \|$.

\begin{PROP}{}\label{PROP:734-34} The space $\Cal H$ contains the space
$\Cal L$ of the special representation of $O(n,1)$ as a simple
summand. If  $\Cal L$ is realized as the space of functions 
$f(\omega )$ on the unit sphere $S^{n-1} \subset \RR^n$, then the embedding
$\tau:\, \Cal L\to\Cal H$ has the form
$$
\tau: \quad f(\omega)  \to \phi(z)= \delta (1-zz^*)\,f(\frac{z}{|z|}),
$$
where $\delta (t)$ is the delta function on $\RR$.
\end{PROP}

The proof is the same as in Proposition~\ref{PROP:7:6}.

\begin{PROP}{}\label{PROP:45-84}  For every vector $\xi \in H^0 \setminus H$,
the function $a(g)= T_g\xi-\xi$ is a nontrivial $1$-cocycle $G\to \Cal H$.
If $\xi=\delta (1-|z|^2)$, then $ a(g) \in \Cal L$,
i.e., it is a $1$-cocycle $G\to \Cal L $.
\end{PROP}

\begin{COR*}{} The subspace $\Cal H'=\Cal H \ominus \Cal L$
does not contain nontrivial $1$-cocycles. (It decomposes into a direct integral
of principal series representations.)
\end{COR*}

The structure of the special representation of the group 
$O(n,1)$ was considered in \S\,2 
of this paper.

\subsubsection{The case of the group $U(n,1)$}

We have
$$
\log |1-zz'^*|=\frac12\log (1-zz'^*)+\frac12\log (1-z'z^*).
$$
Expanding each term of this sum into a power series, we obtain
the following proposition.

\begin{PROP}{}\label{PROP:29-48} The norm $\|\phi\|$ on $H$ 
can be written as the sum
$$
\|\phi\|^2=\frac12\,\|\phi\|^2_++\frac12\,\|\phi\|^2_-,
$$
where
$$
\|\phi\|^2_+=\sum_{k> 0}\frac{(|k|-1)!}{k!}\,
\Bigl|\int \OVER z^k\,\phi(z)\,dz\Bigr|^2,
\qquad
\|\phi\|^2_-\sum_{k> 0}\frac{(|k|-1)!}{k!}\,
\Bigl|\int z^k\,\phi(z)\,dz\Bigr|^2.
$$
These norms, as well as the original norm $\|\phi\|$, are invariant under
the action of $U(n,1)$.
\end{PROP}

It follows from the explicit expression for $\|\phi\|$ that this norm is positive;
however, in contrast to the case of the group $O(n,1)$, it degenerates on some
infinite-dimensional subspace $H_0$. Let $\Cal H$ be the completion of 
the quotient space $H/H_0$ with respect to this norm.

The invariance of the norms $\|\phi\|_{\pm}$ implies the following proposition.

\begin{PROP}{}\label{PROP:895-4} The space $\Cal H$ is the direct sum of two
invariant subspaces:
$$
\Cal H=\Cal H_+ \oplus   \Cal H_-,
$$
where $\Cal H_+ = \{\phi  \in\Cal H \mid \|\phi  \|_- =0\} $,
$\Cal H_- = \{\phi  \in\Cal H \mid \|\phi  \|_+=0\}$.
\end{PROP}

The representations of the group $U(n,1)$ in these subspaces are irreducible,
and the projections of every nontrivial $1$-cocycle $U(n,1)\to H$
to these subspaces do not vanish. Thus the representations of
$U(n,1)$ in the spaces $\Cal H_{\pm}$ are special.

The structure of these two special representations of the group
$U(n,1)$ was considered in \S\,6 of the paper.

\section{Appendix~2. Derivation of some formulas}\label{sect:8}

\subsection{The Fourier transform of the functions
$(a^2+|\gamma|^2)^{-\lambda/2}$ and $|\gamma|^{-\lambda}$ in  $\RR^{n-1}$, where
$|\gamma|=(\sum_{i=1}^{n-1}\gamma_i^2)^{1/2}$,} is given by
\begin{equation}{}\label{veron6:(2):17-3}
\int_{\RR^{n-1}} (a^2+|\gamma|^2)^{- \lambda /2}\,
        e^{\,i\,\langle \xi,\gamma\rangle}\,d\gamma
= c_n\frac{a^{\frac{\lambda }{2}-1}}{2^{\frac{\lambda }{2}-1}\,\Gamma (\lambda /2)}\,
  |\xi|^{\frac{\lambda -n+1}{2}}\,K_{\frac{n-1+ \lambda }{2}}(a|\xi|),
\end{equation}
\begin{equation}{}\label{veron6:(2):17-4}
\int_{\RR^{n-1}} |\gamma|^{-\lambda}\,e^{\,i\,\langle \xi,\gamma\rangle}\,d\gamma =
c_n\frac{2^{\,-\lambda}\, \Gamma ((n-1-\lambda)/2) }{\Gamma (\lambda /2)}.
\,|\xi|^{\lambda -n+1},
\end{equation}
where the coefficient $c_n$, the same in
\eqref{veron6:(2):17-3} and \eqref{veron6:(2):17-4}, depends only on $n$.

The Bessel function  $K_{\rho}$ is given by the formula (see \cite[Vol.~2]{BE})
\begin{equation}{}\label{veron6:17-5}
K_{\rho}(2z) = \frac{\pi}{2\sin(\pi\rho)}\,(I_{-\rho}(2z)-I_{\rho}(2z)),
\end{equation}
where
\begin{equation}{}\label{veron6:17-6}
I_{\rho}(2z) = \sum_{m=0}^\infty \frac{z^{2m+\rho}}{m! \,\Gamma (m+\rho+1)}.
\end{equation}

\begin{REM*}{}
For integer values of $\rho$, the series for $K_{\rho}$ contains terms with
$\log z$; for half-integer values of $\rho$, the expression for $K_{\rho}$ 
can be reduced to a simpler form:
$$
K_{n+\frac12}(z)=(\frac{\pi}{2z})^{1/2}\,e^{-z}\,
\sum_{k=0}^n \frac{(n+k)!}{k!(n-k)!(2z)^k}.
$$
\end{REM*}

Let us derive formulas \eqref{veron6:(2):17-3} and
\eqref{veron6:(2):17-4} for $n>2$. Passing to spherical coordinates,
we can bring the first integral to the form
$$
J_1=c\,\int_0^{\infty }\int_0^{\phi} (a^2+r^2)^{- \lambda /2}\,
e^{\,i|\xi|\cos\phi}\,r^{n-2}\,\sin^{n-3}\phi\,d\phi\,dr.
$$
Integrating with respect to $\phi$, we obtain (see \cite[formula 3.915.5]{G-R})
$$
J_1=c\,|\xi|^{-\frac{n-3}{2}}\int_0^{\infty }
(a^2+r^2)^{- \lambda /2}\,r^{\frac{n-1}{2}}\,
J_{\frac{n-3}{2}}(|\xi|r)\,dr,
$$
where $J_{\rho}(\gamma)$ is a Bessel function of the first kind
(see \cite[Vol.~2]{BE}). Analogously, for the integral 
\eqref{veron6:(2):17-4} we obtain the expression
$$
J_2=c\,|\xi|^{-\frac{n-3}{2}}\int_0^{\infty }
r^{\frac{n-1}{2}- \lambda }\,
J_{\frac{n-3}{2}}(|\xi|r)\,dr.
$$
Integrating with respect to $r$, we obtain for $J_1$ and $J_2$, respectively,
formulas \eqref{veron6:(2):17-3} and \eqref{veron6:(2):17-4},
see \cite[formulas 6.565.4 and 6.561.14]{G-R}.

Formula \eqref{veron6:(2):17-3} immediately implies its complex analog:
\begin{equation*}{}\label{veron6:(2):17-33}
\int_{\CC^{n-1}}
(a^2+|\gamma|^2)^{-\lambda/2}\,e^{\,i\,\langle \xi,\gamma\rangle}\,d\gamma =
c_n\frac{2^{- \lambda /2}}{\Gamma (\lambda /2)}\,
|\xi|^{(\lambda -2n + 2)/2}\,
K_{(2n-2+ \lambda )/2}(a |\xi|).
\end{equation*}

\subsection{Formulas for integrals over the unitary unit sphere $S \subset\CC^n$}

The following formula holds:
\begin{equation}{}\label{953--71}
I(\alpha )\equiv \int_S \prod_{i=1}^n
|\omega _i|^{\alpha _i}\,d \omega =
3\pi^n\,\frac{\prod_{i=1}^n \Gamma (\frac{\alpha _i}{2}+1)}
{\Gamma (\frac12 \sum \alpha _i+n )}.
\end{equation}
In particular,
\begin{equation}{}\label{953--7}
\int_S |\omega ^k|^2\,d \omega =\frac{2\,\pi^n\,k!}{(|k|+n-1)!},\quad
k\in\ZZ_+^n,
\end{equation}
where $\omega ^k=\prod \omega _i^{k_i},$ $k!=\prod k_i!$,
$|k|=\sum k_i$.

\begin{proof}
Consider the following integral over the interior of the unit ball in $\CC^n$:
$$
I = \int_{|z|\le 1}\prod_{i=1}^n |z_i|^{\alpha _i} \,d\mu(z),
$$
where $d\mu(z)$ is the Lebesgue measure on $\CC^n$. Passing to the spherical coordinates
$z=r\, \omega $, $r>0$, $\omega \in S$, we obtain
$$
I=\int_0^1 r^{\sum \alpha_i+2n-1}\, dr\,I(\alpha )=
\frac{1}{2(\frac12\,\sum \alpha_i +n)}\,I(\alpha ).
$$
On the other hand, passing to the polar coordinates $z_j=r_je^{i\phi_j}$,
$j=1, \dots ,n$, we obtain
\begin{gather*}{}
I=(2\pi)^n\,\int_{|r|\le 1} \prod_{j=1}^n r_j^{\alpha _j+1}\,dr_j=
\pi^n\,\int_{\sum s_j\le 1}\prod_{j=1}^n s_j^{\alpha_i/2}\,ds_j
\\
= \pi^n\,\frac{\prod_{i=1}^n \Gamma  (\frac{\alpha _i}{2}+1)}
{\Gamma  (\frac12\sum\alpha _i+n+1)}.
\end{gather*}
Comparing these expressions for $I$, we obtain \eqref{953--71}.
\end{proof}

\subsection{The multidimensional analog of the Dougall formula}

For every $n>1$,
\begin{equation}{}
\begin{gathered}{}\label{22-33-5}
\sum_{|l|=0} \Bigl(\prod_{i=1}^n  \Gamma (l_i+a_i+1)\Gamma (-l_i+b_i+1) \Bigr)^{-1}
\\
=\frac{\Gamma (\sum_{\mathstrut}(a_i+b_i)+1)}{
\Gamma (\sum^{\mathstrut} a_i+1)\,\Gamma (\sum b_i+1)\,\prod_{i=1}^n \Gamma (a_i+b_i+1)},
\end{gathered}
\end{equation}
where $|l|=l_1+ \ldots +l_n$, $l_i\in\ZZ$, or, in a nonsymmetric form,
\begin{equation}\label{22-33-44}
\begin{gathered}{}
\sum_{l'\in\ZZ^{n-1}} \Bigl( \Gamma (-|l'|+a_n+1)\,\Gamma (|'l|+b_n+1)\,
\prod_{i=1}^{n-1} \Gamma (l_i+a_i+1)\Gamma (-l_i+b_i+1) \Bigr)^{-1}
\\
=\frac{\Gamma (\sum_{\mathstrut}(a_i+b_i)+1)}{
\Gamma (\sum^{\mathstrut} a_i+1)\,\Gamma (\sum b_i+1)\,\prod_{i=1}^n \Gamma (a_i+b_i+1)},
\end{gathered}
\end{equation}
where $|l'|=\sum_{i=1}^{n-1} l_i$.

\begin{proof}{} For $n=2$, formula \eqref{22-33-44} coincides with the classical
Dougall formula, see \cite{BE}:
$$
\begin{gathered}{}
\sum_{l=- \infty }^\infty \Bigl(\Gamma (l+ \alpha _1+1)\,
\Gamma (l+ \alpha _2+1)\,\Gamma (-l+ \beta_1+1)\,
\Gamma (-l+ \beta _2+1) \Bigr)^{-1}
\\
=\frac{\Gamma (\alpha _1+ \alpha _2+ \beta _1+ \beta _2+1)}{
\Gamma (\alpha _1+ \beta _1+1)\,\Gamma (\alpha _1+ \beta _2+1)\,
\Gamma (\alpha _2+ \beta _1+1)\,\Gamma (\alpha _2+ \beta _3+1)}.
\end{gathered}
$$
Let us prove \eqref{22-33-44} for every $n>2$ by induction on $n$. 
By the Dougall formula, 
\begin{gather*}{}
\sum_{l_{n-1}=- \infty }^\infty \Bigl(\Gamma (-|l'|+a_n+1)\,
\Gamma (|l'|+b_n+1)\,\Gamma (l_{n-1}+a_{n-1}+1)\,
\Gamma (-l_{n-1}+b_{n-1}+1) \Bigr)^{-1}
\\
=\frac{ \Gamma (a_{n-1}+b_{n-1}+a_n+b_n+1)
        }{
        \Gamma (a_{n-1}+b_{n-1}+1)\,\Gamma (a_n+b_n+1)\,
        \Gamma (-|l''|+a_{n-1}+a_n+1)\,\Gamma ( |l''|+b_{n-1}+b_n+1)},
\end{gather*}
where $l''=l_1+ \ldots +l_{n-2}$. Substituting this equation into
\eqref{22-33-44} yields
$$
I_n = \frac{\Gamma (a_{n-1}+b_{n-1}+a_n+b_n+1)}{\Gamma (a_{n-1}+b_{n-1}+1)\,
\Gamma (a_n+b_n+1)}\, I_{n-1},
$$
where
$$
I_{n-1}=\sum{}'\Bigl(
        \Gamma (-|l''|+a_{n-1}+a_n+1)\,
        \Gamma (|l''|+b_{n-1}+b_n+1)\,
\prod_{i=1}^{n-2} \Gamma (l_i+a_i+1)\Gamma (-l_i+b_i+1) \Bigr)^{-1};
$$
the prime means that the sum ranges only over $l_1, \dots ,l_{n-2}$.
By the induction hypothesis,
$$
I_{n-1}=\frac{\Gamma (\sum(a_i+b_i)+1)}{\Gamma (a_{n-1}+b_{n-1}+a_n+b_n+1)
\Gamma (\sum a_i+1)\Gamma (\sum b_i+1)
\prod_{i=1}^{n-2}\Gamma (a_i+b_i+1)}.
$$
This implies \eqref{22-33-44}.
\end{proof}

\begin{REM*}{}
The following generalization of \eqref{22-33-5} is obvious:
\begin{equation}{}
\begin{gathered}{}\label{22-33-455}
\sum_{|l|=m}
\Bigl(\prod_{i=1}^n  \Gamma (l_i+a_i+1)\,\Gamma (-l_i+b_i+1) \Bigr)^{-1}
\\
=\frac{\Gamma (\sum_{\mathstrut}(a_i+b_i)+1)
        }{
        \Gamma (\sum^{\mathstrut} a_i+m+1)\,\Gamma (\sum b_i-m+1)\,
        \prod_{i=1}^n \Gamma (a_i+b_i+1)}
\end{gathered}
\end{equation}
for every $m\in\ZZ$.
\end{REM*}

\subsection{Calculation of the norm of the vector $f_p$}

The norm of the vector $f_p$ in the space
$L^{\lambda }(\rho )$ equals
\begin{equation}{}\label{44-44-11}
\|\, f_p \, \|^{\,2}_{L} = const\,\rho ^{-|p|}\,p!\,
        \frac{
        \Gamma (n - \lambda )\,\Gamma (|p| + \frac{\lambda }{2})
        }{
        \Gamma (|p| + n - \frac{\lambda }{2})\,\Gamma^{\,2} (\frac{\lambda }{2})
        }.
\end{equation}

\begin{proof}{}
According to the original definition,
$$
\|f_p\|^2=c(\lambda )\,\rho ^{\,2n-2} \int_{\CC^{n-1} \times \CC^{n-1} }
\OVER u^p v^p\,Q(\rho |u-v|^2)\,
e^{\,\rho (vu^* - uu^* - vv^*)}\,d\mu(u)\,d\mu(v),
$$
where
$$
Q(t)=t^{-\frac{\lambda -1}{2}}\,e^{t/2}\,K_{\frac{\lambda -1}{2}}(t/2).
$$

\subsubsection{}
First let us reduce this expression for the norm to
\begin{equation}{}\label{0101-21}
\|f_p\|^2=c(\lambda )\, \rho^{\,-|p|}\,(p\,!)^2\,
\sum_{l\le p} \frac{(-1)^{|l|}}{l!\,(m-l)!\,(|l|+n-2)!}\,
\int_0^\infty s^{|l|+n-2}\,Q(s)\,e^{\,-s}\,ds,
\end{equation}
where $c(\lambda ) = \frac{2\pi^{1/2}}{\Gamma (\lambda /2)}$.

For this, make the change of variables $ u = \rho ^{\,-1/2 }\,u'$,
$ v = \rho ^{\,-1/2 }\,v'$, $u'=z+v'$ and use the identity
$$
\rho ^{\,-1/2 }\,(vu^* - uu^* - vv^*)
= -zv^{\prime *} - v'v^{\prime *} - zz^*.
$$
In the new variables $z$ and $v'$, $ v'\to v$, we obtain
\begin{gather*}{}
\|f_p \|^2=c(\lambda )\, \rho ^{-|p|}\, \int_{\CC^{n-1} \times \CC^{n-1} }
(\OVER z+\OVER v)^p v^p\,Q(|z|^2)\,
e^{\,-(zv^* + vv^* + zz^*)}\,d\mu(z)\,d\mu(v)
\\
= c(\lambda )\,\rho ^{-|p|}\,
\sum_{l\le p} \frac{p\,!}{l!(p-l)!} \int_{\CC^{n-1} \times \CC^{n-1} }
\OVER z^l\OVER v^{p-l} v^p\,Q(|z|^2)\,
e^{\,-(zv^* + vv^* + zz^*)}\,d\mu(z)\,d\mu(v).
\end{gather*}
Let us substitute into this formula the power expansion 
$e^{\,-zv^*}=\sum_{\,k\in\ZZ_+^{n-1}}
\frac{(-1)^{|k|}}{k!} \, z^k\OVER v^k$.
Throwing off the terms that integrate to $0$, we obtain
\begin{gather*}{}
\|f_p \|^2 = c(\lambda )\,\rho ^{-|p|}\,p\,!\,
\sum_{l\le p} \frac{(- 1 )^{|l|}}{(l!)^2(m-l)!}
\int_{\CC^{n-1} \times \CC^{n-1} }
\OVER z^l z^l\OVER v^p v^p\,Q(|z|^2)\,
e^{\,- ( vv^*+zz^*)}\,d\mu(z)\,d\mu(v)
\\
=c(\lambda )\,\rho ^{-|p|}\,p\,!\,\,
\Bigl(\int_{\CC^{n-1}} \,v^p\OVER v^p\,e^{\,- vv^*}\,d\mu(v)\,\Bigr)\,
\sum_{l\le p} \frac{(- 1)^{|l|}}{(l!)^2\,(m-l)!}\,
\int_{\CC^{n-1}} z^l\OVER z^l\,Q(|z|^2)\,e^{\,- |z|^2}\,d\mu(z).
\end{gather*}
The first integral in the obtained equation equals $p\,!$.
The second one can be reduced to the form
\begin{gather*}{}
I=\int_{\RR^{n-1}_+}s^l\,Q({\textstyle\sum\,} s_i)\,e^{\,- \sum s_i}\,\prod ds_i
\\
\int_{\RR^{n-1}_+}(s-\sum_{i=2}^{n-1} s_i)^{l_1}\,
(\prod_{i=2}^{n-1} s_i)\,Q(s)\,e^{\,-s}\,ds
\,\prod_{i=2}^{n-1}ds_i
\\
\frac{l!}{(|l|+n-2)!}\,\int_0^\infty
s^{|l|+n-2}\,Q(s)\,e^{\,-s}\,ds.
\end{gather*}
Thus
$$
\|f_p\|^2=c(\lambda )\, \rho^{\,-|p|}\,(p\,!)^2\,
\sum_{l\le p} \frac{(-1)^{|l|}}{l!\,(m-l)!\,(|l|+n-2)!}\,
\int_0^\infty s^{|l|+n-2}\,Q(s)\,e^{\,-s}\,ds.
$$

\subsubsection{}
Now let us check that substituting the expression for $Q(t)$ into
\eqref{0101-21} yields
\begin{equation}{}\label{3-6}
\begin{gathered}{}
\|f_p\|^2=\frac{2\pi}{\Gamma (\lambda /2)}\rho^{\,-|p|}\,(p\,!)^2\,
\sum_{l\le p}
\frac{(-1)^{|l|}\,\Gamma (|l|+n- \lambda)
        }{
        l\,!\,(p-l)\,!\,\Gamma (|l|+n-\frac{\lambda }{2})
        }
=\pi^{3/2}\,\frac{\rho^{\,-|p|-2n+2}\,(p\,!)^2}{\sin\pi(n- \lambda  )}\,
\times
\\
\times
\sum_{l\in\ZZ^{n-1}}\Bigl( \Gamma (|l|+n-\frac{\lambda }{2})\,
\Gamma (-|l|-n+\lambda +1)\,\prod_{i=1}^{n-1} \Gamma (l_i+1)
\Gamma (p_i-l_i+1) \Bigr)^{-1}.
\end{gathered}
\end{equation}
Let us use the formula
$$
\int_0^\infty e^{\,-as}\,K_{\nu}(as)\,s^{\mu-1}\,ds=\pi^{1/2}
\frac{\Gamma (\mu+\nu)\,\Gamma (\mu-\nu) }{(2a)^{\mu}\, \Gamma (\mu+1/2)},
$$
see  \cite[Vol.~1, p.~331, (28)]{Bateman-Tables}.
For $\nu=\frac{\lambda -1}{2}$,
$\mu=|l|+n-\frac{\lambda +1}{2}$, $a=1/2$, we have
\begin{gather*}{}
\int_0^\infty s^{|l|+n-2}\,Q(s)\,e^{\,-s}\,ds
= \int_0^\infty s^{|l|+n-2-\frac{\lambda -1}{2}}\,
K_{\frac{\lambda -1}{2}}(s/2)\,e^{\,-s/2}\,ds
\\
=\pi^{1/2}\,
\frac{\Gamma (|l|+n-1)\, \Gamma (|l|+n- \lambda)}
{\Gamma (|l|+n-\frac{\lambda }{2})}.
\end{gather*}
This implies \eqref{3-6}.

\subsubsection{}
Finally, apply the Dougall formula \eqref{22-33-44} to the sum \eqref{3-6}.
In our case,  $a_i=0$, $b_i=p_i$ for $i\le n$, $a_n=-n+ \lambda $, and
$a=n-\frac{\lambda }{2}-1$, $b_i=p_i$. Therefore
$$
\|f_p\|^2=c(\lambda )\,\pi^{3/2}\,
\frac{\rho^{\,-|p|}\,(p\,!)^2}{\sin\pi(n- \lambda)}\,
\frac{\Gamma (|p|+\frac{\lambda }{2})
        }{
        \Gamma (\frac{\lambda }{2})\,
        \Gamma (|p|+n-\frac{\lambda }{2})\,
        \Gamma (\lambda -n+1)\,p\,!}.
$$
Substituting the equality $\sin\pi(n- \lambda  )\,\Gamma (\lambda -n+1)
= \pi\, \Gamma^{-1}(n- \lambda )$ yields \eqref{44-44-11}.
\end{proof}

\subsection{Derivation of formula \eqref{7-4} for the operator $T_s$ 
in the Bargmann model:}

$$
T_s^+(\rho_+ ^{\alpha }\,z^k)=(-1)^{|k|}\,\frac{\Gamma (\alpha )}
{\Gamma (- \alpha +|k|)}\rho_+^{\,-\alpha+|k| }\,z^k.
$$
\begin{proof}{}
It follows from Proposition~\ref{PROP:L_k,T_s} that
$$
T_s^+(\rho_+ ^{\alpha }\,z^k)=f(\rho )\,z^k=\tau\,\psi_1(\rho ,z),
$$
where, according to \eqref{249-34}, the function $\psi_1(\rho ,z)$
is given by the formula
\begin{gather*}{}
\psi_i(\rho ,z)= \rho ^{-n+1}\,
\Bigl(\int\limits_0^{ \infty }
\int\limits_{- \infty }^{ \infty }
e^{\,-\rho \zeta +\frac{\,\rho '\mathstrut}{\zeta \mathstrut}}\,
|\zeta |^{-2n}\,\zeta^{-|k|}\,\rho ^{\prime \alpha+n-1 }\,dt\,
d \rho '\Bigr)\,z^k
\\
=(-1)^{|k|}\, \Gamma (\alpha +n)\, \rho ^{-n+1}\,
\Bigl(
\int\limits_{- \infty }^{ \infty }
e^{\,-\rho \zeta }\,
(-\bar\zeta )^{-n}\,(-\zeta)^{-|k|+ \alpha }\,dt\Bigr)\,z^k,
\quad \zeta =it-\frac{|z|^2}{2}.
\end{gather*}

According to \eqref{377-70}, the equality $f(\rho )\,z^k=\tau\,\psi_1(\rho ,z)$
implies that the function  $f(\rho  )$ is given by the formula
$$
f(\rho )=\frac{\rho ^{|k|+n-1}}{k!}\,
\int_{\CC^{n-1}}\psi_1(\rho ,z)\,\bar z^k\,
e^{\,- \rho |z|^2}\,d\mu(z).
$$
Let us proceed to the calculation of this function. Taking into account that
$-\rho \zeta- \rho |z|^2= \rho \bar \zeta$, we have
$$
f(\rho )=\frac{(-1)^{|k|} \Gamma(\alpha +n)  }{k!}\,\rho ^{|k|}\,
\int_{- \infty }^\infty  \int_{\CC^{n-1}} e^{\,\rho \bar \zeta }\,
(-\bar \zeta )^{-n}\,(- \zeta )^{-|k|+ \alpha }\,|z^k|^2\,d \mu (z) \,dt.
$$
Making the change of variables $t\to \rho^{\,-1}t$, 
$z\to \rho^{\,-1/2}z$, we obtain
$$
f(\rho )=\frac{(-1)^{|k|} \Gamma(\alpha +n)  }{k!}\,c(\alpha ,k)\,
\rho^{\,- \alpha + |k|},
$$
where 
$$
c(\alpha ,k)=
\int_{- \infty }^\infty  \int_{\CC^{n-1}} e^{ \bar \zeta }\,
(-\bar \zeta )^{-n}\,(- \zeta )^{-|k| +\alpha }\,|z^k|^2\,d\mu(z)\,dt.
$$
To calculate $c(\alpha ,k)$, replace the complex variables
$z_i$ by the variables $r_i=\frac{|z|^2}{2}$ and $\phi_i=\arg z_i$.
Since the integrand does not depend on $\phi_i$, we obtain, setting 
$r=\sum r_i$,
\begin{gather*}{}
c(\alpha ,k)=c\,
\int_{- \infty }^\infty  \int_{\RR_+^{n-1}} e^{\,-(r+it)}\,
(r+it)^{-n}\,(r-it)^{-|k| + \alpha }\,\prod r_i^{k_i}dr_i\,dt
\\
=c'\,\int_{- \infty }^\infty  \int_0^\infty  e^{\,-(r+it)}\,
(r+it)^{-n}\,(r-it)^{-|k| + \alpha }\, r^{|k|+n-2}\,dr\,dt.
\end{gather*}

Now, setting $r+it=u\,e^{i\phi}$, $u>0$,
$-\frac{\pi}{2}<\phi< \frac{\pi}{2}$, we bring the integral to the form
$$
c(\alpha ,k)=c'\,
        \int_{-\frac{\pi}{2}}^{\frac{\pi}{2}}
        \int_0^\infty  e^{\,-u e^{ i\phi}}\,
        u^{\alpha -1}\,
e^{\,-(-|k|+ \alpha +n)\phi}\,\cos^{|k|+n-2}\phi\,d\mu(u)\,d\phi.
$$
Integrating with respect to $u$ yields
\begin{gather*}{}
c(\alpha ,k)=c'\, \Gamma (\alpha )\,
\int_{-\frac{\pi}{2}}^{\frac{\pi}{2}}
e^{\,-(-|k|+ 2\alpha +n)\phi}\, \cos^{|k|+n-2}\phi\,d\phi
\\
=c''\,\frac{ \Gamma (\alpha )}
{\Gamma (\alpha +n) \Gamma (- \alpha +|k|) },
\end{gather*}
see \cite[p.~386, formula 9, expression for the last integral]{G-R}.

As a result we have
$$
T_s(\rho ^{\alpha }\,z^k)=c\,(-1)^{|k|}\,\frac{\Gamma (\alpha )}
{\Gamma (- \alpha +|k|)}\rho^{\,-\alpha+|k| }\,z^k, \quad c>0.
$$
Since $T^2_s=\opn{id}$, it follows that $c=1$.
\end{proof}

        \let\define\def \let\redefine\def
        \let\aw\bibitem

        \def\No{ No } 

        \define\FA#1#2#3#4{//Функциональный анализ и его приложения.
                          -- #2 -- т.~ #1 ~,   вып. ~#3 ~. \ С. ~#4 ~}
        \def\DAN#1#2#3#4{//ДАН СССР. 
                           -- #2  -- т.~ #1 ~,   \No #3 ~. \ С. ~#4 ~}
        \define\UMN#1#2#3#4{//Успехи Математических Наук.
                            -- #2  -- т.~ #1 ~,   вып.~#3 ~. \ С. ~#4 ~}
        \define\TrudyMMO#1#2#3{//Труды Московского Математического Об\-щест\-ва.
                           -- #2  -- т.~ #1 ~. \  С. ~#3 ~}
        \define\IPM#1#2#3{//Пре\-принт ИПМ АН СССР
                            \No #1 за #2~. M: #2 ~. #3 ~с~}
        \define\MSBOR#1#2#3#4{//Математический сборник.
                            -- ~#1 ( #2 ) за #3 ~.\  С. ~#4 ~}

\makeatletter
\makeatother

\let\sectionname\null
\end{document}